\documentclass[final]{siamltex}
\usepackage[utf8]{inputenc}
\usepackage{color}
\usepackage{array}
\usepackage{float}
\usepackage{rotfloat}
\usepackage{epstopdf}
\usepackage{graphicx}
\usepackage{algorithmicx}
\usepackage{booktabs}
\usepackage{amsfonts,amssymb, latexsym}
\usepackage{mathtools}
\usepackage{algpseudocode}
\usepackage{amsmath}
\usepackage{xcolor}
\usepackage{epstopdf}
\usepackage{url}

\floatstyle{ruled}
\newfloat{algorithm}{tbp}{loa}
\floatname{algorithm}{Algorithm}

\newtheorem{thm}{Theorem}[section]
\newtheorem{cor}[thm]{Corollary}
\newtheorem{lem}[thm]{Lemma}
\newcommand{\pl}{\hspace{-2pt}+\hspace{-2pt}}
\newcommand{\Y}{\mathcal{Y}}
\newcommand{\Yhk}{\hat{\mathcal{Y}}_k}
\newcommand{\Ghk}{\hat{G}_k}
\newcommand{\whkj}{\hat{w}'_{k,j}}
\newcommand{\whkjt}{\hat{w}'^T_{k,j}}
\newcommand{\vhkj}{\hat{v}'_{k,j}}
\newcommand{\vhkjt}{\hat{v}'^T_{k,j}}
\newcommand{\uhkj}{\hat{u}'_{k,j}}

\newcommand{\vhskj}{\hat{v}_{sk+j}}
\newcommand{\uhskj}{\hat{u}_{sk+j}}
\newcommand{\vhskjt}{\hat{v}^T_{sk+j}}

\newcommand{\eps}{\varepsilon}
\newcommand{\epss}{{\varepsilon^2}}
\newcommand{\mn}{\hspace{-2pt}-\hspace{-2pt}}
\newcommand{\lvvert}{\hspace{1pt}\vert}
\newcommand{\rvvert}{\vert\hspace{1pt}}
\newcommand{\op}{\hspace{1pt}}
\newcommand{\tp}{\hspace{2pt}}

\title{Mixed Precision $s$-step Lanczos and Conjugate Gradient Algorithms}
\author{Erin Carson\thanks{Charles University, Faculty of Mathematics and Physics, \{carson,gergelits\}@karlin.mff.cuni.cz. Both authors were funded by Charles University PRIMUS project
No. PRIMUS/19/SCI/11 and Lawrence Livermore National Security, LLC Subcontract Award B639388 under Prime Contract No. DE-AC52-07NA27344. The first author was additionally supported by Charles University Research Program No. UNCE/SCI/023.} \and Tom{\'a}{\v s} Gergelits\footnotemark[1]}
\date{December 2020}

\begin{document}

\maketitle

\begin{abstract}
Compared to the classical Lanczos algorithm, the $s$-step Lanczos variant has the potential to improve performance by asymptotically decreasing the synchronization cost per iteration. However, this comes at a cost. Despite being mathematically equivalent, the $s$-step variant  is known to behave quite differently in finite precision, with potential for greater loss of accuracy and a decrease in the convergence rate relative to the classical algorithm. It has previously been shown that the errors that occur in the $s$-step version follow the same structure as the errors in the classical algorithm, but with the addition of an amplification factor that depends on the square of the condition number of the $O(s)-$dimensional Krylov bases computed in each outer loop. As the condition number of these $s$-step bases grows (in some cases very quickly) with $s$, this limits the parameter $s$ that can be chosen and thus limits the performance that can be achieved. 
In this work we show that if a select few computations in $s$-step Lanczos are performed in double the working precision, the error terms then depend only linearly on the conditioning of the $s$-step bases. This has the potential for drastically improving the numerical behavior of the algorithm with little impact on per-iteration performance. Our numerical experiments demonstrate the improved numerical behavior possible with the mixed precision approach, and also show that this improved behavior extends to the $s$-step CG algorithm in mixed precision. 
\end{abstract}



\pagestyle{myheadings}
\thispagestyle{plain}
\markboth{CARSON AND GERGELITS}{MIXED PRECISION S-STEP LANCZOS}

\section{Introduction}
The Lanczos algorithm \cite{lanczos1950iteration, lanczos1952solution} is a popular approach to solving large, sparse, symmetric eigenvalue problems. Given an $n\times n$ matrix $A$ and a starting vector $v_1$, after $n$ iterations the Lanczos algorithm constructs an orthogonal matrix $V=[v_1,\ldots v_n]$ and a tridiagonal matrix $T$ such that $AV=VT$. The matrices $A$ and $T$ thus have the same eigenvalues, and the eigenvectors of $A$ can easily be obtained from those of $T$. Because $T$ is tridiagonal, its eigendecomposition can be easily computed. The columns of $V$ form an orthogonal basis for the Krylov subspace 
\[
\mathcal{K}_n(A,v_1)=\text{span}\{v_1, Av_1, A^2v_1, \ldots, A^{n-1}v_1\}.
\]

Even after $m<n$ iterations, the eigenvalues and eigenvectors of $T$ can give close approximations to some of those of $A$, and this is how the method is most often used in practice. We show the classical Lanczos algorithm in Algorithm~\ref{alg:Lanczos}. Note that this variant uses auxiliary vectors $u_m$ in the computation; this 2-term recurrence formulation is equivalent to modified Gram-Schmidt and is thus numerically superior to the alternative 3-term recurrence variant. The convergence criteria used to stop the iterations is not important in the current setting; we only assume that we execute some number $m\leq n$ iterations. 
We will also assume in the current work that no breakdown occurs during the algorithm, i.e., we do not encounter a 0 value of $\beta_{m+1}$ in line \ref{lanczos:beta}.

The conjugate gradient (CG) method \cite{hestenes1952methods} is the method of choice for solving symmetric positive definite linear systems $Ax=b$, and is based on an underlying Lanczos process. Given an initial approximate solution $x_1$ with corresponding residual $r_1=b-Ax_1$, in iteration $m$ the CG method selects updates an approximate solution $x_m\in x_1+\mathcal{K}_m(A,r_1)$ such that $r_k \perp \mathcal{K}_m(A,r_1)$. This is equivalent to solving $T_m y_m = \Vert r_1 \Vert e_1$ and then taking $x_m = x_1 + V_m y_m$, where $V_m$ forms a basis for $\mathcal{K}_m(A,r_1)$, where we note that this $T$ and $V$ are the same as those from Lanczos. Mathematically, this means that $x_m$ is chosen in each iteration to be the vector in $x_1+\mathcal{K}_m(A,r_1)$ that minimizes the $A$-norm of the error, that is, $\Vert x_m - x \Vert_A = (x_m-x)^TA(x_m-x)$, and we are guaranteed to converge to the exact solution after $n$ iterations. 

In practice, of course, both Lanczos and CG, which are based on short recurrences, suffer from the effects of finite precision computation. Namely, the basis vectors $V$ are subject to loss of orthogonality, which can cause deviation from the exact processes. 

Implementations of the classical formulations of Lanczos and CG are often limited by communication. This is because each iteration consists of interdependent sparse matrix-vector multiplications (SpMVs) and inner products, which both suffer from low computational intensity and create a communication bottleneck. For example, the High-Performance CG Benchmark (HPCG), used as a complement to LINPACK for ranking supercomputers, only achieves about 3\% of the peak performance on the top machine as of November 2020 \cite{top500} (compared to 82\% for the LINPACK benchmark).

This has inspired a number of approaches for alleviating the performance bottleneck in CG, Lanczos, and other Krylov subspace methods. One particular approach is termed $s$-step Krylov subspace algorithms. This approach has a long history, with much early work by Chronopoulos and Gear \cite{chronopoulos1989s} and others. See \cite{hoemmen2010communication} and \cite{carson2015communication} for overviews of historical references. 

The idea here is that iterations are computed in blocks of $s$. First, one computes an $O(s)$-dimensional Krylov subspace basis with the current iteration vectors, such that all vectors that will be computed in the next $s$ iterations can be obtained from this subspace. Under certain constraints on matrix sparsity and partitioning, this can be accomplished with $O(1)$ messages in parallel \cite{demmel2008avoiding, ballard2014communication}. This ``matrix powers computation'' is followed by a block orthogonalization; in short recurrence methods like CG this involves the computation of a Gram matrix, and in long recurrence methods such as GMRES this is done via a block Gram-Schmidt procedure. See \cite[Chapters 7 \& 8]{ballard2014communication} for details. In this way, the per-iteration latency cost in the parallel case is reduced by an asymptotic factor $O(s)$. We note that there are other approaches that hide synchronization cost, such as pipelined Krylov subspace methods; see, e.g., \cite{ghysels2014hiding}. 

It has long been known that these $s$-step formulations are potentially unstable. In \cite{carson2015communication}, the maximum attainable accuracy of $s$-step CG is analyzed, and the results of Paige for classical Lanczos were extended to the $s$-step case. The main source of instability is in the potentially ill-conditioned $O(s)$-dimensional bases computed via repeated SpMVs in each outer loop iteration. One simple solution is to use Chebyshev or Newton polynomials to reduce the condition number of the bases. Other approaches to improving numerical behavior include residual replacement \cite{carson2014residual} and using adaptive $s$ values \cite{carson2018adaptive, carson2020adaptive}.

A recent development is the inclusion of multiprecision capabilities in hardware. Modern GPUs such as the NVIDIA V100 and A100, a crucial part of exascale architectures, offer half precision (16-bit) to double precision (64-bit), as well as special TensorCore instructions in various precisions, which perform a $4\times 4$ matrix multiply in one clock cycle with extended precision accumulation of inner products. It is expected that the range of available precisions will only expand going forward, and there are many efforts to develop new algorithms that exploit this hardware; see, e.g., \cite{abdelfattah2020survey}. With this new multiprecision ecosystem available, it is natural to ask whether we exploit multiple precisions to improve the numerical behavior of $s$-step Krylov subspace methods without significant performance overhead.

In this work, we develop mixed precision $s$-step Lanczos algorithm. In this algorithm, a given working precision is used for all computations and storage of all quantities, with the exception that the Gram matrix is computed and applied in double the working precision. By extending the work of Paige for classical Lanczos \cite{paige1976error, paige1980accuracy}, as well as the work in \cite{carson2015accuracy} for the uniform precision $s$-step Lanczos algorithm, we give theoretical error bounds on the loss of orthogonality and other quantities in the mixed precision $s$-step Lanczos algorithm, and also state implications for eigenvalue and eigenvector convergence. 

The main result is that in this mixed precision regime, compared to the uniform precision case, the loss of orthogonality and loss of normality of the Lanczos basis vectors is reduced by a factor relating to the condition number of the computed $s$-step bases. This means that we expect numerical behavior closer to that of the classical Lanczos algorithm. We present a few small numerical experiments that support the theoretical analysis. 

We also extend this approach to the $s$-step CG algorithm and develop a mixed precision $s$-step CG algorithm that makes use of extended precision in the same way. Numerical results confirm that behavior is improved compared to the uniform precision case, and is closer to that of the classical CG algorithm. 

We argue that, especially in latency-bound cases, the performance overhead per iteration of using double the working precision in these select computations will be minimal. When the improved convergence rate is accounted for, we expect that the mixed precision approach will, in many cases, provide an improved time-to-solution. We also argue that the approach of using mixed precision is complementary to other approaches, such as residual replacement and the use of better polynomial bases. 

The remainder of the paper is outlined as follows. In Section \ref{sec:relatedwork}, we discuss the related work on which the present work is based. In Section \ref{sec:algs}, we give a brief derivation of $s$-step Lanczos from the classical Lanczos algorithm in order to establish notation. Our mixed precision $s$-step Lanczos algorithm is represented and rounding errors are analyzed in Section \ref{sec:vtv}, and a few numerical experiments are presented in Section \ref{sec:num}. The extension of Paige's results \cite{paige1980accuracy} to the mixed precision $s$-step Lanczos case is discussed in Section \ref{sec:evals}. Section \ref{sec:cg} details the corresponding mixed precision $s$-step CG algorithm along with numerical experiments, and we discuss open problems and future work in Section \ref{sec:conclusion}.

\section{Related work}
\label{sec:relatedwork}

Our analysis relies heavily on existing results, in particular those of Paige for the classical Lanczos algorithm. In \cite{paige1976error}, Paige performed a rounding error analysis of the Lanczos algorithm, giving bounds on loss of orthogonality and other important quantities. This analysis was subsequently used in his seminal analysis \cite{paige1980accuracy}, which rigorously showed, amongst other fundamental results, that the loss of orthogonality is linked with the convergence of eigenvalue approximations. We recommend \cite{meurant2006lanczos} for an insightful overview of the significance of Paige's work. 

In \cite{carson2015accuracy}, the results of Paige in \cite{paige1976error, paige1980accuracy} are extended to the (uniform precision) $s$-step Lanczos algorithm. The key result is that, under a constraint on the condition numbers of the $s$-step bases, the main structure of the analysis still applies, but with almost all bounds being amplified by quantities related to the basis condition numbers (or their squares). We heavily reuse these results here, with only small modifications to the analysis to account for the use of mixed precision. We do, however, give new insights into the meaning of Paige's analysis for the (uniform and mixed precision) $s$-step case; see Section \ref{sec:evals}. 

The seminal work of Greenbaum \cite{greenbaum1989behavior} expanded upon the results of Paige, and developed a ``backward-like'' error analysis for the CG algorithm (in particular, a variant sometimes referred to as ``CG/Lanczos''). Greenbaum's analysis says that finite precision CG for a linear system with matrix $A$ behaves like exact CG with a larger matrix $\tilde{A}$ whose eigenvalues lie in tight intervals around the eigenvalues of $A$. The analysis of Greebaum has not yet been extended to the $s$-step variants of CG, but we make some conjectures on what the results would look like in Section \ref{sec:cg}.

\section{The $s$-step Lanczos method}
\label{sec:algs}
The classical Lanczos method is shown in Algorithm~\ref{alg:Lanczos}. For simplicity, we will assume no breakdown occurs and thus breakdown conditions are not discussed here. We present a brief derivation of $s$-step Lanczos from the classical Lanczos recurrences in order to establish notation.

\begin{algorithm}[h]
\caption{Lanczos}\label{alg:Lanczos}
\begin{algorithmic}[1]
\Require {$n$-by-$n$ real symmetric matrix $A$ and length-$n$ starting vector $v_1$ such that $\Vert v_1 \Vert_2=1$}
\State {$u_1 = Av_1$}
\For {$m=1,2,\ldots$ until convergence}
\State {$\alpha_m = v_m^Tu_m$} \label{lanczos:alpha}
\State {$w_m = u_m - \alpha_m v_m$} \label{lanczos:w}
\State {$\beta_{m+1} = \Vert w_m \Vert_2$} \label{lanczos:beta}
\State {$v_{m+1} = w_m/\beta_{m+1}$}\label{lanczos:v}
\State {$u_{m+1} = Av_{m+1} - \beta_{m+1} v_m$}\label{lanczos:u}
\EndFor
\end{algorithmic}
\end{algorithm}

Suppose we are beginning iteration $m=sk+1$ where $k\in\mathbb{N}$ is the number of outer loops we have performed so far and $0<s\in\mathbb{N}$ is the number of inner iterations within each outer loop. By induction on lines~\ref{lanczos:v} and~\ref{lanczos:u} of Algorithm~\ref{alg:Lanczos}, we can write
\begin{equation}
v_{sk+j}, u_{sk+j} \in \mathcal{K}_{s+1}(A,v_{sk+1}) + \mathcal{K}_{s+1}(A,u_{sk+1}),
\label{eq:vudeps}
\end{equation}
for $j\in\{1,\ldots,s+1\}$, where $\mathcal{K}_{i}(A,x)=\text{span}\{x,Ax,\ldots,A^{i-1}x\}$ denotes the Krylov subspace of dimension $i$ of matrix $A$ with respect to vector $x$. Note in the case $k=0$ (the first outer loop iteration), we have
\begin{equation*}
v_{j}, u_{j} \in \mathcal{K}_{s+2}(A,v_{1}),
\end{equation*}
for $j\in\{1,\ldots,s+1\}$, since $u_1=Av_1$. 

Thus to perform the next $s$ inner loop iterations, we need bases for the Krylov subspaces in \eqref{eq:vudeps}. We define $\mathcal{Y}_k=[\mathcal{V}_k, \mathcal{U}_k]$, where $\mathcal{V}_k$ and $\mathcal{U}_k$ are size $n$-by-$(s+1)$ matrices whose columns form bases for $\mathcal{K}_{s+1}(A,v_{sk+1})$ and $\mathcal{K}_{s+1}(A,u_{sk+1})$, respectively. For the case $k=0$, we can save on computation and define $\mathcal{Y}_0$ to be a size $n$-by-$(s+2)$ matrix whose columns span $\mathcal{K}_{s+2}(A,v_{1})$. We refer to these $\mathcal{Y}_k$'s as `$s$-step basis matrices'. 

Given these computed bases,~\eqref{eq:vudeps} says that we can represent the length-$n$ vectors $v_{sk+j}$ and $u_{sk+j}$ by their coordinates (denoted with primes) in $\mathcal{Y}_k$, i.e.,
\begin{equation}
v_{sk+j}=\mathcal{Y}_k v'_{k,j}, \qquad u_{sk+j}=\mathcal{Y}_k u'_{k,j},
\label{eq:basischange}
\end{equation}
for $j\in \{1,\ldots,s+1\}$. Note that we define the coordinate vectors to be length $s+2$ for $k=0$ and length $2s+2$ for $k>0$. Because the vector $w_{sk+j}$ is a linear combination of $u_{sk+j}$ and $v_{sk+j}$, we can also represent it by its coordinates in $\mathcal{Y}_k$, i.e., $w_{sk+j}=\mathcal{Y}_k w'_{k,j}$ for $j\in\{1,\ldots,s\}$. 

Another advantage of having the $s$-step basis matrix is that we can precompute all the inner products between basis vectors, and this only need be done once per outer loop. We define the Gram matrix $G_k = \mathcal{Y}_k^T \mathcal{Y}_k$, which is size $(s+2)$-by-$(s+2)$ for $k=0$ and $(2s+2)$-by-$(2s+2)$ for $k>0$. In this way, we can write the inner products in lines~\ref{lanczos:alpha} and~\ref{lanczos:beta} as 
\begin{align}
\alpha_{sk+j} &= v_{sk+j}^T u_{sk+j} = v'^{T}_{k,j}\mathcal{Y}_k^T \mathcal{Y}_k u'_{k,j} = v'^{T}_{k,j}G_k u'_{k,j} \quad \text{and}\label{eq:slanczos:alpha}\\
\beta_{sk+j+1} &= (w_{sk+j}^T w_{sk+j})^{1/2} = (w'^{T}_{k,j}\mathcal{Y}_k^T \mathcal{Y}_k w'_{k,j})^{1/2} = (w'^{T}_{k,j}G_k w'_{k,j})^{1/2}.\label{eq:slanczos:beta}
\end{align} 

Another question is what polynomials are used for generating the basis vectors in $\mathcal{Y}_k$. This will turn out to have a large impact numerically. In practice, one could use any polynomial basis desired. We define matrix $\mathcal{B}_k$ which stores the coefficients of the polynomials, i.e., $\mathcal{B}_k$ is a matrix such that 
\begin{equation*}
A \underline{\hat{\mathcal{Y}}}_k = \hat{\mathcal{Y}}_k \mathcal{B}_k
\label{eq:b}
\end{equation*}
where $\mathcal{B}_k$ is size $(s+2)$-by-$(s+2)$ for $k=0$ and size $(2s+2)$-by-$(2s+2)$ for $k>0$, and $\underline{\hat{\mathcal{Y}}}_k = \big[\hat{\mathcal{V}}_k [I_s,0_{s,1}]^T, 0_{n,1},\hat{\mathcal{U}}_k [I_s,0_{s,1}]^T, 0_{n,1} \big]$. We note that $\mathcal{B}_k$ is in general upper Hessenberg but often tridiagonal in practice since we have a symmetric $A$. 

Therefore, for $j\in\{1,\ldots,s\}$, we can write the sparse matrix-vector products as 
\begin{equation}
Av_{sk+j+1} = A\mathcal{Y}_k v'_{k,j+1} = A \underline{\hat{\mathcal{Y}}}_k v'_{k,j+1} =\mathcal{Y}_k\mathcal{B}_k v'_{k,j+1}.
\label{eq:AVVB}
\end{equation}

Given these quantities, the $s$-step variant of Lanczos works as follows. In some outer loop $k$, starting at global iteration $sk+1$, we will compute global iterations $sk+2$ through $sk+s+1$. In the outer loop, we will generate the $s$-step basis matrix $\mathcal{Y}_k$ such that~\eqref{eq:b} holds. Under certain assumptions on the sparsity structure of the matrix and its partitioning in the parallel case, this can be accomplished with $O(1)$ messages; see \cite{hoemmen2010communication}. We can then compute the Gram matrix $G_k$, which in parallel requires one Allreduce collective. This is all the communication that need happen within outer loop $k$.  

The inner loop iterations can then proceed without communication. Updates to the length-$n$ vectors are performed implicitly by updating instead the length-$(2s+2)$ coordinates for those vectors in the basis $\mathcal{Y}_k$. The matrices $\mathcal{B}_k$ and $G_k$ are small (square with dimension $O(s)$), and so they can be stored locally on each processor (or kept in cache in the sequential case). Thus the matrix-vector products can be done locally (as in \eqref{eq:AVVB}) and the inner products can also be computed locally using the Gram matrix as in~\eqref{eq:slanczos:alpha} and~\eqref{eq:slanczos:beta}. After the inner loop iterations finish, the length-$n$ vectors $v_{sk+j+1}$ and $u_{sk+j+1}$, for $j\in\{1,\ldots, s\}$ can be recovered via a single block multiplication of the corresponding coordinate vectors with the basis matrix $\mathcal{Y}_k$.  

We show the resulting $s$-step Lanczos algorithm, as it appears in \cite{carson2015accuracy}, in Algorithm~\ref{alg:sLanczos}. 
As described above, in practice one would only perform the basis change operation~\eqref{eq:basischange} on a block of coordinate vectors at the end of each outer loop, however, in Algorithm~\ref{alg:sLanczos} we have shown the length-$n$ vector updates in each inner iteration (lines~\ref{slanczos:v} and~\ref{slanczos:u}) for clarity.

\begin{algorithm}[h]
\caption{$s$-step Lanczos}\label{alg:sLanczos}
\begin{algorithmic}[1]
\Require {$n$-by-$n$ real symmetric matrix $A$ and length-$n$ starting vector $v_1$ such that $\Vert v_1 \Vert_2=1$}
\State {$u_1 = Av_1$}
\For {$k=0,1,\ldots$ until convergence}
\State {Compute $\mathcal{Y}_k$ with change of basis matrix $\mathcal{B}_k$}\label{slanczos:basis}
\State {Compute $G_k= \mathcal{Y}_k^T\mathcal{Y}_k$} \label{slanczos:G}
\State {$v'_{k,1}=e_{1}$}
\If {$k=0$} \State{$u'_{0,1}=\mathcal{B}_0 e_{1}$}
\Else 
\State{$u'_{k,1}=e_{s+2}$}
\EndIf
\For {$j=1,2,\ldots,s$}
\State {$\alpha_{sk+j} = v'^{T}_{k,j}G_k u'_{k,j}$} \label{slanczos:alpha}
\State {$w'_{k,j} = u'_{k,j} - \alpha_{sk+j} v'_{k,j}$} \label{slanczos:wcoeff}
\State {$\beta_{sk+j+1} = ({w'^{T}_{k,j} G_k w'_{k,j}})^{1/2}$} \label{slanczos:beta}
\State {$v'_{k,j+1} = w'_{k,j}/\beta_{sk+j+1}$} \label{slanczos:vcoeff}
\State {$v_{sk+j+1} =  \mathcal{Y}_kv'_{k,j+1} $} \label{slanczos:v}
\State {$u'_{k,j+1} = \mathcal{B}_k v'_{k,j+1} - \beta_{sk+j+1} v'_{k,j}$} \label{slanczos:ucoeff}
\State {$u_{sk+j+1} =  \mathcal{Y}_k u'_{k,j+1} $} \label{slanczos:u}
\EndFor
\EndFor
\end{algorithmic}
\end{algorithm}

\section{The $s$-step Lanczos method in mixed precision}
\label{sec:vtv}

Throughout our analysis, we use a standard model of floating point arithmetic where we assume the computations are carried out on a machine with relative working precision $\eps$ (see~\cite{golub1996matrix}). We ignore underflow and overflow. Following Paige~\cite{paige1976error}, we use the $\eps$ symbol to represent the relative precision as well as terms whose absolute values are bounded by the relative precision. 

The analysis in \cite{carson2015accuracy} suggests that a mixed precision approach can potentially have significant benefit. It is clear that the square of the $s$-step basis condition number enters the bounds through the formation of the Gram matrix $G_k$. We therefore suggest a mixed precision approach as follows. In each outer loop iteration, the Gram matrix $G_k$ should be computed and stored in precision $\epss$ (double the working precision). This will double the number of bits moved, but since $G_k$ is modestly-sized (square with dimension $2s+2$), this will not cause significant overhead, particularly in latency-bound cases. This will not affect the number of synchronizations, and further, only needs to occur every $s$ iterations. Within each inner loop, $G_k$ is applied twice to a single vector. The inner products involved in this matrix-vector multiplication should be accumulated in precision $\epss$, but the result can be stored in the working precision $\eps$. Again, since $G_k$ is of small dimension, and these matrix-vector multiplies are done locally on each processor, this is an insignificant cost to performance. All other computations are performed in the working precision $\eps$. We summarize this approach in Algorithm~\ref{alg:msLanczos}.

\begin{algorithm}[h]
\caption{Mixed precision $s$-step Lanczos}\label{alg:msLanczos}
\begin{algorithmic}[1]
\Require {$n$-by-$n$ real symmetric matrix $A$ and length-$n$ starting vector $v_1$ such that $\Vert v_1 \Vert_2=1$}
\State {$u_1 = Av_1$ (precision $\eps$)}
\For {$k=0,1,\ldots$ until convergence}
\State {Compute $\mathcal{Y}_k$ with change of basis matrix  $\mathcal{B}_k$ (precision $\eps$).}\label{mslanczos:basis}
\State {Compute and store $G_k= \mathcal{Y}_k^T\mathcal{Y}_k$ in precision $\epss$.}  \label{mslanczos:G}
\State {$v'_{k,1}=e_{1}$}
\If {$k=0$} \State{$u'_{0,1}=\mathcal{B}_k e_{1}$}
\Else 
\State{$u'_{k,1}=e_{s+2}$}
\EndIf
\For {$j=1,2,\ldots,s$}
\State {Compute $g = G_k u'_{k,j}$ in precision $\epss$, store in precision $\eps$.}
\State {$\alpha_{sk+j} = v'^{T}_{k,j}g$ (precision $\eps$)} \label{mslanczos:alpha}
\State {$w'_{k,j} = u'_{k,j} - \alpha_{sk+j} v'_{k,j}$ (precision $\eps$)} \label{mslanczos:wcoeff}
\State {Compute $c=G_k w'_{k,j}$ in precision $\epss$, store in precision $\eps$.}
\State {$\beta_{sk+j+1} = ({w'^{T}_{k,j} c})^{1/2}$ (precision $\eps$)} \label{mslanczos:beta}
\State {$v'_{k,j+1} = w'_{k,j}/\beta_{sk+j+1}$ (precision $\eps$)} \label{mslanczos:vcoeff}
\State {$v_{sk+j+1} =  \mathcal{Y}_kv'_{k,j+1} $ (precision $\eps$)} \label{mslanczos:v}
\State {$u'_{k,j+1} = \mathcal{B}_k v'_{k,j+1} - \beta_{sk+j+1} v'_{k,j}$ (precision $\eps$)} \label{mslanczos:ucoeff}
\State {$u_{sk+j+1} =  \mathcal{Y}_k u'_{k,j+1} $ (precision $\eps$)} \label{mslanczos:u}
\EndFor
\EndFor
\end{algorithmic}
\end{algorithm}

We will model floating point computation in a precision $\eps$ using the following standard conventions (see, e.g., ~\cite{golub1996matrix}): for vectors $u,v \in \mathbb{R}^{n}$, matrices $A\in\mathbb{R}^{n\times m}$ and $G\in\mathbb{R}^{n\times n}$, and scalar $\alpha$, 
\begin{align*}
fl_\eps(u-\alpha v) =& u - \alpha v - \delta w,  &&| \delta w | \leq (| u| + 2 |\alpha v |)\eps,\\
fl_\eps(v^T u) =& (v+\delta v)^Tu,  &&| \delta v | \leq n \eps | v |, \\
fl_\eps(Au) =& (A+\delta A)u, &&|\delta A| \leq m\eps | A |, \quad\text{and}\\
fl_\eps(A^TA) =& A^TA + \delta E,  &&| \delta E | \leq n\eps | A^T | |A|.
\end{align*} 
where $fl_\eps()$ represents the evaluation of the given expression in floating point arithmetic with unit roundoff $\eps$ and terms with $\delta$ denote error terms. We decorate quantities computed in finite precision arithmetic with hats, e.g., if we are to compute the expression $\alpha=v^Tu$ in finite precision, we get $\hat{\alpha}=fl_\eps(v^T u)$.

We will use the following lemma from \cite{carson2015accuracy}, which will be useful in our analysis. The proof is trivial and thus omitted. 
\begin{lem}
\label{lem:condnum}
Assume we have rank-$r$ matrix $Y \in \mathbb{R}^{n\times r}$, where $n\geq r$. Let $Y^{+}$ denote the pseudoinverse of $Y$, i.e., $Y^{+} = (Y^{T}Y)^{-1}Y^{T}$. Then for any vector $x\in\mathbb{R}^{r}$, we can bound 
\begin{equation*}
\Vert \lvvert Y \rvvert \lvvert x \rvvert \Vert_2 \leq \Vert \lvvert Y \rvvert \Vert_2 \Vert  x  \Vert_2 \leq \Gamma \Vert Y x \Vert_2.
\end{equation*}
where $\Gamma = \big\Vert Y^{+} \big\Vert_2 \op \big\Vert \lvvert Y \rvvert \big\Vert_2\leq \sqrt{r}\op\big\Vert Y^{+} \big\Vert_2 \op\big\Vert Y \big\Vert_2$.
\end{lem}

We note that the term $\Gamma$ can be thought of as a type of condition number for the matrix $Y$. In the analysis, we will apply the above lemma to the computed `basis matrix' $\hat{\mathcal{Y}}_k$; in particular, we will use the definition $\Gamma_k = \big\Vert \Y_k^{+} \big\Vert_2 \op \big\Vert \lvvert \Y_k \rvvert \big\Vert_2$. We assume throughout that the generated bases $\hat{\mathcal{U}}_k$ and $\hat{\mathcal{V}}_k$ are numerically full rank. That is, all singular values of $\hat{\mathcal{U}}_k$ and $\hat{\mathcal{V}}_k$ are greater than $\epsilon n \cdot 2^{\lfloor{\log_2{\theta_1}}\rfloor}$ where $\theta_1$ is the largest singular value of $A$. We further make the key assumption that $\eps n \Gamma_k \ll 1$ for all $k$. Based on this assumption, throughout the analysis we will drop terms of order $(\eps n \Gamma_k)^2$. Including these terms in the analysis would not change the fundamental structure of the results, but would result in a serious overestimate of the bounds due to growing constant terms. 

The results of this section are summarized in the following theorem, which has the same structure as that of the uniform precision $s$-step Lanczos algorithm in \cite{carson2015accuracy}, which in turn has the same structure as Paige's results for classical Lanczos in \cite{paige1976error}.

 \begin{thm}
\label{thm:mainthm}
Assume that Algorithm~\ref{alg:sLanczos} is implemented in floating point with working precision $\epsilon$ and applied for $m=sk+j$ steps to the $n$-by-$n$ real symmetric matrix $A$ \color{black} with at most $N$ nonzeros per row, \color{black} starting with vector $v_1$ with $\| v_1 \|_2 = 1$. \color{black} Let $\sigma\equiv \| A \|_2$, $\theta\sigma= \| |A| \|_2 $ and $\tau_k \sigma= \| |\mathcal{B}_k| \|_2 $, \color{black} where $\mathcal{B}_k$ is defined in~\eqref{eq:AVVB}, and let 
\begin{equation*}
\bar{\Gamma}_k =\hspace{-1mm} \max_{i\in \{0,\ldots,k\} } \Vert \hat{\mathcal{Y}}_i^{+} \Vert_2 \Vert \lvvert \hat{\mathcal{Y}}_i \rvvert \Vert_2\geq 1  \quad\text{and} \quad \bar{\tau}_k =\hspace{-1mm} \max_{i\in\{0,\ldots,k\} }  \tau_i,
\end{equation*}
where above the superscript `+' denotes the Moore-Penrose pseudoinverse, i.e., $\hat{\mathcal{Y}}_i^{+} = (\hat{\mathcal{Y}}_i^T \hat{\mathcal{Y}}_i)^{-1} \hat{\mathcal{Y}}_i^T$. Further, assume that 
\begin{equation}
2m\eps\big( (N\pl 2s\pl 5)\theta + (4s\pl 9)\bar{\tau}_k  \pl 34s\pl 55\big)\bar\Gamma_k\ll 1
    \label{ass1}
\end{equation}
and that $\eps n \bar\Gamma_k \ll 1$. Then $\hat\alpha_{sk+j}$, $\hat\beta_{sk+j+1}$, and $\hat{v}_{sk+j+1}$ will be computed such that, for $i \in \{1,\ldots, m\}$,
\begin{equation*}
A\hat{V}_{m} = \hat{V}_{m} \hat{T}_{m} + \hat\beta_{m+1}\hat{v}_{m+1}e^T_{m} - \delta \hat{V}_{m},\\
\end{equation*}
with
\begin{align*}
\hat{V}_{m} & = [\hat{v}_1,\hat{v}_2,\ldots, \hat{v}_{m}]\\
\delta\hat{V}_{m} & = [\delta\hat{v}_1,\delta\hat{v}_2,\ldots, \delta\hat{v}_{m}]\\
\hat{T}_{m} & = \left[\begin{array}{cccc}
\hat\alpha_1 & \hat{\beta}_2 & & \\
\hat{\beta}_2 & \ddots & \ddots & \\
 & \ddots & \ddots & \hat\beta_{m} \\
 & & \hat\beta_{m} & \hat\alpha_{m} 
\end{array}\right]
\end{align*}
and
\color{black}
\begin{align}
\Vert  \delta\hat{v}_{i} \Vert_2 \leq& \hspace{1mm} \eps_1 \sigma, \label{mainthm:deltav}\\ \vspace{2mm}
\hat\beta_{i+1}\vert \hat{v}_{i}^T\hat{v}_{i+1} \vert \leq& \hspace{1mm}  \eps_0 \sigma,  \label{mainthm:orth}\\ \vspace{2mm}
\vert \hat{v}_{i+1}^T \hat{v}_{i+1} -1 \vert \leq& \hspace{1mm}  \eps_0/2,\quad\text{and} \label{mainthm:norm} \\
\Big\vert \hat\beta_{i+1}^2  + \hat\alpha_{i}^2 + \hat{\beta}_{i}^2 - \Vert A \hat{v}_{i} \Vert_2^2 \Big\vert \leq& 
2  i (3\eps_0 + 2\eps_1) \sigma^2, \label{mainthm:cols}
\end{align}
\color{black}

\color{black}
where
\begin{equation}
\eps_0 \equiv 2\eps(9s\pl 14)\bar{\Gamma}_k \quad\text{and}\quad
\eps_1 \equiv \eps\big( (N\pl 2s\pl 5)\theta + (4s\pl 9)\bar{\tau}_k + (10s \pl 16) \big) \bar{\Gamma}_k.
\label{eq:e0e1}
\end{equation}
\color{black}

Furthermore, if $R_{m}$ is the strictly upper triangular matrix such that 
\begin{equation}
\hat{V}_{m}^T \hat{V}_{m} = {R}_{m}^T + \text{diag}(\hat{V}_{m}^T \hat{V
}_{m}) + {R}_{m},
\label{eq:Rm}
\end{equation}
then 
\begin{equation}
\hat{T}_{m}R_{m} - R_{m}\hat{T}_{m} = \hat\beta_{m+1}\hat{V}_{m}^T \hat{v}_{m+1}e_{m}^T + H_{m},
\label{mainthm:orthrec}
\end{equation}
where $H_{m}$ is upper triangular with elements $\eta$ such that 
\color{black}
\begin{equation}\label{mainthm:etas}
\begin{split}
\vert\eta_{1,1}\vert \leq & \eps_0\sigma, \hspace{2mm} \text{and, for} \hspace{2mm} i\in\{2,\ldots,m\}, \\
\vert \eta_{i,i} \vert \leq &  2\eps_0 \sigma , \\
\vert \eta_{i-1,i} \vert \leq & (\eps_0 + 2\eps_1)\sigma,  \hspace{2mm}\text{and}\\
\vert \eta_{\ell,i} \vert  \leq & 2 \eps_1 \sigma, \hspace{1mm}\text{for}\hspace{2mm}\ell\in\{1,\ldots,i\mn 2\}. 
\end{split}
\end{equation}
\color{black}

\end{thm}

\paragraph{Remarks} 
We reiterate that this theorem has the same structure as that for uniform precision $s$-step Lanczos appearing in \cite{carson2015accuracy}, but with the notable exception that \textit{the term $\eps_0$ now contains only a factor of $\bar{\Gamma}_k$ rather than $\bar{\Gamma}_k^2$}. As $\Gamma_k$ can potentially grow very quickly with $s$, this is a significant improvement, and indicates that, among other things, the Lanczos basis vectors will maintain significantly better orthogonality and normality due to the selective use of higher precision. We again note that this theorem also has the same structure as the theorem given by Paige~\cite{paige1976error} for classical Lanczos. 

We briefly discuss the meaning of the bounds in Theorem~\ref{thm:mainthm}, which give insight into how orthogonality is lost in the mixed precision $s$-step Lanczos algorithm. Equation~\eqref{mainthm:deltav} bounds the error in the columns of the perturbed Lanczos recurrence. Equation ~\eqref{mainthm:norm} bounds how far the Lanczos vectors deviate from normality, and~\eqref{mainthm:orth} bounds loss of orthogonality between adjacent vectors. The bound~\eqref{mainthm:cols} describes the deviation of the columns of $A\hat{V}_{m}$ and $\hat{T}_{m}$. Finally,~\eqref{mainthm:orthrec} gives a recurrence for the loss of orthogonality between Lanczos vectors and shows how errors propagate through the iterations. Note that $\| |\mathcal{B}_k|\|_2$ depends on the polynomial basis used in generating the $s$-step basis vectors, and should be $\lesssim\| |A|\|_2$ in practice. We also note that the assumption~\eqref{ass1} is needed for the analysis but is likely to be overly strict in practice. 

\subsection{Proof of Theorem~\ref{thm:mainthm}}
The remainder of this section is focused on the proof of Theorem~\ref{thm:mainthm}. We stress that this analysis follows almost exactly the same structure as that for the uniform precision $s$-step Lanczos in \cite{carson2015accuracy} which in turns follows almost exactly the same structure as the proof of Paige \cite{paige1976error}. The proof is trivially adapted from those in \cite{carson2015accuracy} and \cite{paige1976error}, but we include the entire derivation here for posterity. 

Recall that we assume that $\eps n \bar{\Gamma}_k\ll 1$. As previously stated, for presentation purposes, we will thus exclude terms of order $(\eps n \bar{\Gamma}_k)^2$ and higher from the analysis. We repeat that this does not affect the fundamental structure or spirit of the results. The main difference between this analysis and the analysis in \cite{carson2015accuracy} is that here, we perform certain computations in double the working precision, leading to terms of order $O(\eps^2)$ which, by our assumptions, we can ignore.

We first proceed toward proving~\eqref{mainthm:norm}.
We construct the Gram matrix in line~\ref{slanczos:G} of Algorithm~\ref{alg:sLanczos} in double the working precision. This gives
\begin{equation}
\hat{G}_k = fl_\epss(\hat{\mathcal{Y}}_k^T\hat{\mathcal{Y}}_k) = \hat{\mathcal{Y}}_k^T\hat{\mathcal{Y}}_k + \delta G_k, \quad \text{where}\quad  | \delta G_k | \leq \epss n | \hat{\mathcal{Y}}_k^T | | \hat{\mathcal{Y}}_k |.
\label{eq:Gk}
\end{equation}
We must also apply the Gram matrix in double the working precision within the inner loop. Note that the computed Gram matrix and vector to be multiplied are very small and likely fit in cache, so this computation in double the working precision is not expected to add significant cost in terms of bits moved. Thus we now turn to the computation of $\hat{\beta}_{sk+j+1}$ in 
line~\ref{slanczos:beta} of Algorithm~\ref{alg:sLanczos}.

We first apply the Gram matrix in double the working precision, computing
\begin{equation}
    c = fl_\epss(\hat{G}_k \hat{w}'_{k,j}) = (\hat{G}_k + \delta\hat{G}_{k,w_k})\hat{w}'_{k,j}, \quad |\delta \hat{G}_{k,w_j}| \leq \epss(2s+2)|\hat{G}_k|.  \label{eq:Gw}
\end{equation}
This result is then rounded to working precision, giving 
\begin{equation*}
    \hat{c}=fl_\eps(c) = c+\delta c, 
\end{equation*}
where 
\begin{equation*}
    |\delta c| \leq \eps|c| \leq \eps\left(|\hat{G}_k \hat{w}'_{k,j}| + \epss(2s\pl 2)|\hat{G}_k||\hat{w}'_{k,j}| \right).
\end{equation*}
Then let
\begin{align}
    d = fl_\eps(\whkjt \hat{c}) &= (\whkj+\delta\whkj)^T\hat{c}, \quad |\delta\whkj|\leq \eps (2s\pl 2)|\whkj| \label{eq:deltawhat}\\
    &=\whkjt \Yhk^T \Yhk \whkj + \delta d \nonumber\\
    &=\Vert \Yhk \whkj \Vert_2^2 + \delta d,\nonumber 
\end{align}
where
\begin{align*}
    \delta d &= \whkjt \delta G_{k} \whkj + \whkjt \delta \hat{G}_{k,w_j} \whkj+ \whkjt \delta c \\
    &\phantom{=}+ \delta \whkjt \hat{G}_k \whkj     + \delta \whkjt \delta \hat{G}_{k,w_j} \whkj + \delta \whkjt \delta c.
\end{align*}
We can then write the bound 
\begin{align*}
    |\delta d| &\leq \epss n |\whkjt||\Yhk^T||\Yhk||\whkj| + \epss(2s \pl 2) |\whkjt||\Yhk^T||\Yhk||\whkj| + \eps |\whkjt||\Yhk^T||\Yhk\whkj| \\
    &\phantom{=}+ \eps^3(2s\pl 2)|\whkjt||\Yhk^T||\Yhk||\whkj| 
    +\eps(2s\pl 2)|\whkjt||\Yhk^T||\Yhk\whkj| \\
    &\phantom{=}+ \eps^3(2s\pl 2)^2 |\whkjt||\Yhk^T||\Yhk||\whkj| +\epss(2s\pl 2) |\whkjt||\Yhk^T||\Yhk\whkj|\\
    &\phantom{=} + \eps^4(2s\pl 2)^2 |\whkjt||\Yhk^T||\Yhk||\whkj| \\
    &\leq \epss n \Gamma_k^2 \Vert \Yhk \whkj\Vert_2^2 + \epss(2s\pl 2)\Gamma_k^2 \Vert \Yhk \whkj\Vert_2^2 + \eps \Gamma_k \Vert \Yhk \whkj\Vert_2^2  \\
    &\phantom{\leq} + \eps^3(2s\pl 2)\Gamma_k^2 \Vert \Yhk \whkj\Vert_2^2 +\eps(2s\pl 2)\Gamma_k \Vert \Yhk \whkj\Vert_2^2 + \eps^3 (2s\pl 2)^2 \Gamma_k^2 \Vert \Yhk \whkj\Vert_2^2\\
    &\phantom{\leq} +\epss(2s\pl 2) \Gamma_k \Vert \Yhk \whkj\Vert_2^2 + \eps^4 (2s\pl 2)^2 \Gamma_k^2 \Vert \Yhk \whkj\Vert_2^2.
\end{align*}
Now, using the assumption that $\eps n \Gamma_k$ is sufficiently less than 1 and omitting higher order terms, we have
\begin{align}
    |\delta d| &\leq \eps(2s\pl 3)\Gamma_k \Vert \Yhk \whkj\Vert_2^2. \label{eq:cbound}
\end{align}
Then writing
\begin{equation*}
d = \big\Vert \Yhk \whkj\big\Vert_2^2 + \delta d=\big\Vert \Yhk \whkj\big\Vert_2^2 + \delta d \cdot \frac{\big\Vert \Yhk \whkj \big\Vert_2^2}{\big\Vert \Yhk \whkj \big\Vert_2^2}
= \big\Vert \Yhk \whkj \big\Vert_2^2 \left( 1 + \frac{\delta d}{\big\Vert \Yhk \whkj \big\Vert_2^2} \right),
\end{equation*}
the computation of $\hat\beta_{sk+j+1}$ can be written
\begin{equation}
\hat\beta_{sk+j+1} \hspace{-1.5pt}=\hspace{-1.5pt} fl_\eps( \sqrt{d}) \hspace{-1pt}=\hspace{-2pt} \sqrt{d} +\delta \beta_{sk+j+1} \hspace{-1pt}\approx\hspace{-1pt} \|\Yhk \whkj \|_2 \hspace{-1mm}\left(1\pl \frac{\delta d}{2\|\Yhk \whkj\|_2^2}\hspace{-1pt}\right) \hspace{-1pt}\pl \delta \beta_{sk+j+1}, 
\label{eq:betahat} 
\end{equation}
where
\begin{equation}
\quad| \delta \beta_{sk+j+1} | \leq \eps \sqrt{d} \lesssim \eps \| \Yhk \whkj\|_2 \left(1+ \frac{\eps(2s\pl 3)\Gamma_k}{2} \right).
\label{eq:deltabeta}
\end{equation}

For the computed $\hat{v}'_{k,j+1}$ we have
\begin{equation}
\hat{v}'_{k,j+1} = fl_\eps (\whkj/\hat\beta_{sk+j+1})= (\whkj + \delta \tilde{w}'_{k,j})/\hat{\beta}_{sk+j+1},
\label{eq:hatvcoeff}
\end{equation}
where
\begin{equation}
\quad | \delta \tilde{w}'_{k,j} | \leq \eps | \whkj |.
\label{eq:deltavp}
\end{equation}

The corresponding (length-$n$) Lanczos vector $\hat{v}_{sk+j+1}$, as well as $\hat{u}_{sk+j+1}$,  are recovered via a change of basis using the basis matrix $\Y_k$. In finite precision, we have 
\begin{equation}
\hat{v}_{sk+j+1} = fl_\eps (\Yhk \hat{v}'_{k,j+1})=\Big(\Yhk + \delta\hat{\mathcal{Y}}_{k,v_{j+1}}\Big) \hat{v}'_{k,j+1},\quad \vert \delta\hat{\mathcal{Y}}_{k,v_{j+1}} \vert \leq \eps(2s+ 2)\vert \Yhk \vert, 
\label{eq:bchangev}
\end{equation}
and
\begin{equation}
\hat{u}_{sk+j+1} = fl_\eps (\Yhk \hat{u}'_{k,j+1})=\Big(\Yhk + \delta\hat{\mathcal{Y}}_{k,u_{j+1}} \Big) \hat{u}'_{k,j+1},\quad \vert \delta\hat{\mathcal{Y}}_{k,u_{j+1}} \vert \leq \eps(2s+ 2)\vert \Yhk \vert.
\label{eq:bchangeu}
\end{equation}

We can now prove~\eqref{mainthm:norm} in Theorem~\ref{thm:mainthm}. Using~\eqref{eq:betahat},~\eqref{eq:hatvcoeff}, and~\eqref{eq:bchangev},
%
\begin{align*}
\hat{v}_{sk+j+1}^T \hat{v}_{sk+j+1} &=  \hat{v}'^T_{k,j+1} (\hat{\mathcal{Y}}_k + \delta\hat{\mathcal{Y}}_{k,v_{j+1}})^T (\hat{\mathcal{Y}}_k + \delta\hat{\mathcal{Y}}_{k,v_{j+1}}) \hat{v}'_{k,j+1} \\
&= \left( \frac{\hat{w}'_{k,j} + \delta \tilde{w}'_{k,j}}{\hat{\beta}_{sk+j+1}} \right)^T (\hat{\mathcal{Y}}_k^T \hat{\mathcal{Y}}_k + 2\delta\hat{\mathcal{Y}}_{k,v_{j+1}}^T \hat{\mathcal{Y}}_k ) \left( \frac{\hat{w}'_{k,j} + \delta \tilde{w}'_{k,j}}{\hat{\beta}_{sk+j+1}} \right) \\
&= \frac{\| \hat{\mathcal{Y}}_k \hat{w}'_{k,j}\|_2^2 + 2\hat{w}'^T_{k,j} \delta\hat{\mathcal{Y}}_{k,v_{j+1}}^T\hat{\mathcal{Y}}_k \hat{w}'_{k,j} + 2\delta \tilde{w}'^T_{k,j} \hat{\mathcal{Y}}_k^T \hat{\mathcal{Y}}_k \hat{w}'_{k,j} }{\hat{\beta}_{sk+j+1}^2} \\
&= \frac{\| \hat{\mathcal{Y}}_k \hat{w}'_{k,j}\|_2^2 + 2\hat{w}'^T_{k,j} \delta\hat{\mathcal{Y}}_{k,v_{j+1}}^T\hat{\mathcal{Y}}_k \hat{w}'_{k,j} + 2\delta \tilde{w}'^T_{k,j} \hat{\mathcal{Y}}_k^T \hat{\mathcal{Y}}_k \hat{w}'_{k,j} }{\| \hat{\mathcal{Y}}_k \hat{w}'_{k,j}\|_2^2  + (\delta d+2\| \hat{\mathcal{Y}}_k \hat{w}'_{k,j}\|_2 \cdot \delta \beta_{sk+j+1})} \\
&= \frac{\| \hat{\mathcal{Y}}_k \hat{w}'_{k,j} \|_2^4}{\| \hat{\mathcal{Y}}_k \hat{w}'_{k,j} \|_2^4} - \frac{\| \hat{\mathcal{Y}}_k \hat{w}'_{k,j} \|_2^2(\delta d+2\| \hat{\mathcal{Y}}_k \hat{w}'_{k,j}\|_2 \cdot \delta \beta_{sk+j+1})}{\| \hat{\mathcal{Y}}_k \hat{w}'_{k,j} \|_2^4} \\
&\phantom{=} + \frac{2\| \hat{\mathcal{Y}}_k \hat{w}'_{k,j} \|_2^2 (\hat{w}'^T_{k,j} \delta\hat{\mathcal{Y}}_{k,v_{j+1}}^T\hat{\mathcal{Y}}_k \hat{w}'_{k,j} + \delta \tilde{w}'^T_{k,j} \hat{\mathcal{Y}}_k^T \hat{\mathcal{Y}}_k \hat{w}'_{k,j} )  }{\| \hat{\mathcal{Y}}_k \hat{w}'_{k,j} \|_2^4}\\
&= 1 - \frac{\delta d + 2\| \hat{\mathcal{Y}}_k \hat{w}'_{k,j}\|_2 \cdot \delta \beta_{sk+j+1}}{\| \hat{\mathcal{Y}}_k \hat{w}'_{k,j} \|_2^2}\\
&\phantom{=} + \frac{ 2 (\hat{w}'^T_{k,j} \delta\hat{\mathcal{Y}}_{k,v_{j+1}}^T\hat{\mathcal{Y}}_k \hat{w}'_{k,j} +  \delta \tilde{w}'^T_{k,j} \hat{\mathcal{Y}}_k^T \hat{\mathcal{Y}}_k \hat{w}'_{k,j} )}{\| \hat{\mathcal{Y}}_k \hat{w}'_{k,j} \|_2^2}.
\end{align*}
Using the bounds in~\eqref{eq:Gk},~\eqref{eq:deltawhat},~\eqref{eq:Gw},~\eqref{eq:cbound},~\eqref{eq:bchangev},~\eqref{eq:bchangeu}, and Lemma~\ref{lem:condnum}, we obtain the bound 
\begin{align*}
\vert \hat{v}_{sk+j+1}^T \hat{v}_{sk+j+1} -1 \vert \leq& \eps (2s\pl 3) \Gamma_k + 2\eps + 2\eps (2s\pl 2)\Gamma_k + 2\eps \Gamma_k\\
\leq& \eps(6s \pl 11) \Gamma_k.
\end{align*}
Notice that due to the use of mixed precision, this bound only depends linearly on the quantity $\Gamma_k$. This thus proves~\eqref{mainthm:norm}. We now proceed toward proving~\eqref{mainthm:orth}. 

Similarly as before, to compute $\hat{\alpha}_{sk+j}$ in line~\ref{slanczos:alpha} of Algorithm~\ref{alg:sLanczos}, we first apply the Gram matrix in double the working precision, round the result to the working precision, and then compute the inner product of the two vectors in the working precision. Letting 
\begin{equation*}
    g = fl_\epss(\Ghk + \delta \hat{G}_{k,u_j})\uhkj, \quad |\delta \hat{G}_{k,u_j}|\leq \epss(2s\pl 2)|\Ghk|,
\end{equation*}
we then round the result to working precision to obtain 
\begin{equation*}
    \hat{g}=fl_\eps(g)=g+\delta g, \quad |\delta g|\leq \eps |g| \leq \eps\left(|\Ghk \uhkj|+\epss(2s\pl 2)|\Ghk||\uhkj| \right).
\end{equation*}
Then to compute $\hat{\alpha}_{sk+j}$, we have 
\begin{equation*}
\hat\alpha_{sk+j} = fl_\eps(\vhkjt \hat{g}) = (\vhkj + \delta \vhkj)^T\hat{g}, \quad |\delta \vhkj|\leq\eps(2s\pl 2)|\vhkj|,
\end{equation*}
and expanding the above equation and using~\eqref{eq:Gk}, and~\eqref{eq:bchangev}, we obtain
\begin{align} 
\hat\alpha_{sk+j}&= (\vhkj + \delta \vhkj)^T\left( (\Ghk + \delta \hat{G}_{k,u_j})\uhkj + \delta g \right) \nonumber\\
&= \vhkjt \Ghk \uhkj + \vhkjt \delta \hat{G}_{k,u_j}\uhkj + \vhkjt \delta g + \delta\vhkjt\Ghk\uhkj \nonumber\\
&=\vhkj\Yhk^T\Yhk\uhkj + \vhkjt \delta G_k \uhkj + \vhkjt \delta \hat{G}_{k,u_j}\uhkj + \vhkjt \delta g + \delta\vhkjt\Ghk\uhkj \nonumber\\
&=(\vhskj-\delta\hat{\Y}_{k,v_j}\vhkj)^T (\uhskj-\delta\hat{\Y}_{k,u_j} \uhkj) \nonumber\\
&\phantom{=}+ \vhkjt \delta G_k \uhkj + \vhkjt \delta \hat{G}_{k,u_j}\uhkj + \vhkjt \delta g + \delta\vhkjt\Ghk\uhkj \nonumber\\
&=\vhskjt\uhskj - \vhskjt\delta\hat{\Y}_{k,u_j}\uhkj - \vhkjt\delta\hat{\Y}_{k,v_j}\uhskj \nonumber\\
&\phantom{=}+ \vhkjt \delta G_k \uhkj + \vhkjt \delta \hat{G}_{k,u_j}\uhkj + \vhkjt \delta g + \delta\vhkjt\Ghk\uhkj \nonumber \\
&= \vhskjt\uhskj + \delta\hat{\alpha}_{sk+j}, \label{eq:alphaexp}
\end{align}
with 
\[\delta\hat{\alpha}_{sk+j} =  \delta\hat{v}'^{T}_{k,j}\hat{G}_k \hat{u}'_{k,j} +
\vhkjt\delta g+\hat{v}'^{T}_{k,j} (\delta G_k + \delta\hat{G}_{k,u_j} - \hat{\mathcal{Y}}_k^T \delta \hat{\mathcal{Y}}_{k,u_j} - \delta \hat{\mathcal{Y}}^T_{k,v_j} \hat{\mathcal{Y}}_k ) \hat{u}'_{k,j}.
\]

Using bounds in~\eqref{mainthm:norm},~\eqref{eq:Gk},~\eqref{eq:deltawhat},~\eqref{eq:Gw},~\eqref{eq:bchangev}, and~\eqref{eq:bchangeu}, and applying Lemma~\ref{lem:condnum}, we can write (again, ignoring higher order terms),
\begin{align*}
\vert \delta\hat{\alpha}_{sk+j} \vert \leq& \eps (2s\pl 2) | \hat{v}'^{T}_{k,j} | | \hat{\mathcal{Y}}_k^T | | \hat{\mathcal{Y}}_k \hat{u}'_{k,j} | + \eps|\vhkjt||\Yhk^T||\Yhk\uhkj| \\ 
&\phantom{\leq} + \epss n | \hat{v}'^{T}_{k,j} | | \hat{\mathcal{Y}}_k^T | | \hat{\mathcal{Y}}_k | |\hat{u}'_{k,j} |  \\
&\phantom{\leq} + \epss (2s\pl 2) | \hat{v}'^{T}_{k,j} | | \hat{\mathcal{Y}}_k^T | | \hat{\mathcal{Y}}_k | |\hat{u}'_{k,j} |  \\
&\phantom{\leq}+\eps(2s\pl 2) | \hat{v}'^{T}_{k,j}  \hat{\mathcal{Y}}_k^T | | \hat{\mathcal{Y}}_k | |\hat{u}'_{k,j} |  \\
&\phantom{\leq}+\eps(2s\pl 2) | \hat{v}'^{T}_{k,j} | | \hat{\mathcal{Y}}_k^T | | \hat{\mathcal{Y}}_k \hat{u}'_{k,j} |  \\
&\leq \eps (6s\pl 7) \Gamma_k \|  \vhskj \|_2 \hspace{1mm} \| \uhskj \|_2,
\end{align*}
and using~\eqref{mainthm:norm}, we obtain
\begin{equation}
    \vert \delta\hat{\alpha}_{sk+j} \vert \leq \eps (6s\pl 7) \Gamma_k  \| \uhskj \|_2.\label{eq:deltaalpha}
\end{equation}
From~\eqref{eq:alphaexp}, and using the bounds in~\eqref{eq:deltaalpha} and~\eqref{mainthm:norm}, we can then write 
\begin{align}
| \hat{\alpha}_{sk+j} | \leq& \| \hat{v}_{sk+j} \|_2 \| \hat{u}_{sk+j} \|_2 + \vert \delta\hat{\alpha}_{sk+j} \vert \nonumber \\
\leq& \big(1+(\eps/2)(6s \pl 11) \Gamma_k\big) \| \hat{u}_{sk+j} \|_2 +\eps ( 6s\pl 7)\Gamma_k \| \hat{u}_{sk+j} \|_2  \nonumber \\
\leq& \Big(1 + (\eps/2) \big(18s\pl 25\big) \Gamma_k \Big) \| \hat{u}_{sk+j} \|_2.
\label{eq:alphabound}
\end{align}

In finite precision, line~\ref{slanczos:wcoeff} of Algorithm~\ref{alg:sLanczos} is computed as
\begin{equation}
\hat{w}'_{k,j} = \hat{u}'_{k,j} -\hat\alpha_{sk+j}\hat{v}'_{k,j} -\delta {w}'_{k,j}, \quad\text{where}\quad \vert \delta {w}'_{k,j} \vert \leq \eps (\vert \hat{u}'_{k,j} \vert + 2\vert \hat\alpha_{sk+j} \hat{v}'_{k,j}\vert ).
\label{eq:deltaw}
\end{equation}
Multiplying both sides of~\eqref{eq:deltaw} by $\hat{\mathcal{Y}}_k$ gives
\begin{equation*}
\hat{\mathcal{Y}}_k \hat{w}'_{k,j} = \hat{\mathcal{Y}}_k \hat{u}'_{k,j} - \hat{\alpha}_{sk+j} \hat{\mathcal{Y}}_k \hat{v}'_{k,j} - \hat{\mathcal{Y}}_k \delta w'_{k,j},
\end{equation*}
and multiplying each side by its own transpose, we get
\begin{align*}
\hat{w}'^{T}_{k,j} \hat{\mathcal{Y}}_k^T \hat{\mathcal{Y}}_k \hat{w}'_{k,j} &= (\hat{\mathcal{Y}}_k \hat{u}'_{k,j} \mn \hat{\alpha}_{sk+j} \hat{\mathcal{Y}}_k \hat{v}'_{k,j} \mn \hat{\mathcal{Y}}_k \delta w'_{k,j})^T(\hat{\mathcal{Y}}_k \hat{u}'_{k,j} \mn \hat{\alpha}_{sk+j} \hat{\mathcal{Y}}_k \hat{v}'_{k,j} \mn \hat{\mathcal{Y}}_k \delta w'_{k,j})\\
&= \hat{u}'^{T}_{k,j} \hat{\mathcal{Y}}_k^T \hat{\mathcal{Y}}_k \hat{u}'_{k,j} - 2\hat{\alpha}_{sk+j} \hat{u}'^{T}_{k,j} \hat{\mathcal{Y}}_k^T \hat{\mathcal{Y}}_k \hat{v}'_{k,j} + \hat{\alpha}^2_{sk+j} \hat{v}'^{T}_{k,j} \hat{\mathcal{Y}}_k^T \hat{\mathcal{Y}}_k \hat{v}'_{k,j} \\
&\phantom{=} - \delta w'^{T}_{k,j} \hat{\mathcal{Y}}_k^T (\hat{\mathcal{Y}}_k \hat{u}'_{k,j} \mn \hat{\alpha}_{sk+j}\hat{\mathcal{Y}}_k \hat{v}'_{k,j})  \mn (\hat{\mathcal{Y}}_k \hat{u}'_{k,j} \mn \hat{\alpha}_{sk+j}\hat{\mathcal{Y}}_k \hat{v}'_{k,j})^T \hat{\mathcal{Y}}_k \delta w'_{k,j} .
\end{align*}
Then using~\eqref{eq:bchangev} and~\eqref{eq:bchangeu}, we can write
\begin{align*}
\hat{w}'^{T}_{k,j} \hat{\mathcal{Y}}_k^T \hat{\mathcal{Y}}_k \hat{w}'_{k,j} &= 
(\hat{u}_{sk+j} - \delta \hat{\mathcal{Y}}_{k,u_j} \hat{u}'_{k,j} )^T(\hat{u}_{sk+j} - \delta \hat{\mathcal{Y}}_{k,u_j} \hat{u}'_{k,j}) \\
&\phantom{=} - 2\hat{\alpha}_{sk+j} (\hat{u}_{sk+j} - \delta \hat{\mathcal{Y}}_{k,u_j} \hat{u}'_{k,j} )^T (\hat{v}_{sk+j} - \delta \hat{\mathcal{Y}}_{k,v_j} \hat{v}'_{k,j} ) \\
&\phantom{=}+ \hat{\alpha}^2_{sk+j}  (\hat{v}_{sk+j} - \delta \hat{\mathcal{Y}}_{k,v_j} \hat{v}'_{k,j} )^T (\hat{v}_{sk+j} - \delta \hat{\mathcal{Y}}_{k,v_j} \hat{v}'_{k,j} ) \\
&\phantom{=}- 2 \delta w'^{T}_{k,j} \hat{\mathcal{Y}}_k^T (\hat{\mathcal{Y}}_k \hat{u}'_{k,j} - \hat{\alpha}_{sk+j} \hat{\mathcal{Y}}_k \hat{v}'_{k,j}) \\
&= \hat{u}^T_{sk+j} \hat{u}_{sk+j} - 2\hat{u}^T_{sk+j}  \delta \hat{\mathcal{Y}}_{k,u_j} \hat{u}'_{k,j} - 2\hat{\alpha}_{sk+j}\hat{u}^T_{sk+j} \hat{v}_{sk+j}\\
&\phantom{=}+ 2\hat{\alpha}_{sk+j} \hat{u}^T_{sk+j} \delta \hat{\mathcal{Y}}_{k,v_j} \hat{v}'_{k,j} + 2\hat{\alpha}_{sk+j} \hat{u}'^T_{k,j} \delta \hat{\mathcal{Y}}_{k,u_j}^T \hat{v}_{sk+j} \\
&\phantom{=} + \hat{\alpha}^2_{sk+j} \hat{v}^T_{sk+j} \hat{v}_{sk+j} - 2 \hat{\alpha}^2_{sk+j} \hat{v}^T_{sk+j} \delta \hat{\mathcal{Y}}_{k,v_j} \hat{v}'_{k,j}\\
&\phantom{=}-2\delta w'^{T}_{k,j} \hat{\mathcal{Y}}_k^T (\hat{\mathcal{Y}}_k \hat{u}'_{k,j} - \hat{\alpha}_{sk+j} \hat{\mathcal{Y}}_k \hat{v}'_{k,j})\\
&= \hat{u}^T_{sk+j} \hat{u}_{sk+j} - 2\hat{\alpha}_{sk+j}\hat{u}^T_{sk+j}\hat{v}_{sk+j} + \hat{\alpha}^2_{sk+j} \hat{v}^T_{sk+j} \hat{v}_{sk+j} \\
&\phantom{=} -2(\delta \hat{\mathcal{Y}}_{k,u_j} \hat{u}'_{k,j} - \hat{\alpha}_{sk+j} \delta \hat{\mathcal{Y}}_{k,v_j} \hat{v}'_{k,j})^T (\hat{u}_{sk+j} - \hat{\alpha}_{sk+j} \hat{v}_{sk+j}) \\
&\phantom{=} -2 \delta w'^{T}_{k,j} \hat{\mathcal{Y}}_k^T (\hat{\mathcal{Y}}_k \hat{u}'_{k,j} - \hat{\alpha}_{sk+j} \hat{\mathcal{Y}}_k \hat{v}'_{k,j}).
\end{align*}
This can be written
\begin{align*}
\Vert \hat{\mathcal{Y}}_k \hat{w}'_{k,j} \Vert_2^2 &= \Vert \hat{u}_{sk+j} \Vert_2^2 - 2\hat{\alpha}_{sk+j}\hat{u}^T_{sk+j} \hat{v}_{sk+j} + \hat{\alpha}^2_{sk+j} \Vert \hat{v}_{sk+j} \Vert_2^2 \\
&\phantom{=} - 2 (\delta \hat{\mathcal{Y}}_{k,u_j} \hat{u}'_{k,j} - \hat{\alpha}_{sk+j} \delta \hat{\mathcal{Y}}_{k, v_j} \hat{v}'_{k,j} + \hat{\mathcal{Y}}_k \delta w'_{k,j} )^T (\hat{u}_{sk+j} - \hat{\alpha}_{sk+j} \hat{v}_{sk+j}),
\end{align*}
where we have used $\hat{\mathcal{Y}}_k \hat{u}'_{k,j} - \hat{\alpha}_{sk+j} \hat{\mathcal{Y}}_k \hat{v}'_{k,j} = \hat{u}_{sk+j} - \hat{\alpha}_{sk+j} \hat{v}_{sk+j} + O(\eps)$.
Now, using~\eqref{eq:alphaexp}, 
\begin{align*}
\Vert \hat{\mathcal{Y}}_k \hat{w}'_{k,j} \Vert_2^2 &= \Vert \hat{u}_{sk+j} \Vert_2^2 - 2\hat{\alpha}_{sk+j}(\hat{\alpha}_{sk+j} - \delta\hat{\alpha}_{sk+j}) + \hat{\alpha}^2_{sk+j} \Vert \hat{v}_{sk+j} \Vert_2^2 \\
&\phantom{=} - 2 (\delta \hat{\mathcal{Y}}_{k,u_j} \hat{u}'_{k,j} - \hat{\alpha}_{sk+j} \delta \hat{\mathcal{Y}}_{k, v_j} \hat{v}'_{k,j} + \hat{\mathcal{Y}}_k \delta w'_{k,j} )^T (\hat{u}_{sk+j} - \hat{\alpha}_{sk+j} \hat{v}_{sk+j}) \\
&= \Vert \hat{u}_{sk+j} \Vert_2^2 + \hat{\alpha}^2_{sk+j} (\Vert \hat{v}_{sk+j} \Vert_2^2 - 2) + 2\hat{\alpha} _{sk+j} \delta\hat{\alpha}_{sk+j} \\
&\phantom{=} - 2 (\delta \hat{\mathcal{Y}}_{k,u_j} \hat{u}'_{k,j} - \hat{\alpha}_{sk+j} \delta \hat{\mathcal{Y}}_{k, v_j} \hat{v}'_{k,j} + \hat{\mathcal{Y}}_k \delta w'_{k,j} )^T (\hat{u}_{sk+j} - \hat{\alpha}_{sk+j} \hat{v}_{sk+j}).
\end{align*}

Rearranging the above equation, we obtain
\begin{multline*}
\Vert \hat{\mathcal{Y}}_k \hat{w}'_{k,j} \Vert_2^2 + \hat{\alpha}^2_{sk+j} - \Vert \hat{u}_{sk+j} \Vert_2^2 = \hat{\alpha}^2_{sk+j} (\Vert \hat{v}_{sk+j} \Vert_2^2 - 1) + 2\hat{\alpha} _{sk+j} \delta\hat{\alpha}_{sk+j} \\
- 2 (\delta \hat{\mathcal{Y}}_{k,u_j} \hat{u}'_{k,j} - \hat{\alpha}_{sk+j} \delta \hat{\mathcal{Y}}_{k, v_j} \hat{v}'_{k,j} + \hat{\mathcal{Y}}_k \delta w'_{k,j} )^T (\hat{u}_{sk+j} - \hat{\alpha}_{sk+j} \hat{v}_{sk+j}).
\end{multline*}
Using Lemma~\ref{lem:condnum} and bounds in~\eqref{mainthm:norm},~\eqref{eq:bchangev},~\eqref{eq:bchangeu},~\eqref{eq:deltaalpha},~\eqref{eq:alphabound}, and~\eqref{eq:deltaw}, along with the assumption $\eps(n\pl 8s \pl 13)\Gamma_k\leq 1$, we can then write
\begin{align}
\Vert \hat{\mathcal{Y}}_k \hat{w}'_{k,j} \Vert_2^2 \pl \hat{\alpha}^2_{sk+j} \mn \Vert \hat{u}_{sk+j} \Vert_2^2 
&\leq \Big(1\pl (\eps/2) \big(18s\pl 25\big)\Gamma_k\Big)^2 \hspace{-1mm} \cdot \eps (6s\pl 11)\Gamma_k \|\hat{u}_{sk+j} \|_2^2 \nonumber\\
&\phantom{=} + 2\Big(\hspace{-2pt}1\pl (\eps/2) \big(18s\pl 25\big)\Gamma_k\Big) \hspace{-1mm} \cdot \eps(6s\pl 7)\Gamma_k  \|\hat{u}_{sk+j} \|_2^2 \nonumber\\
&\phantom{=} + \eps \big( 16s\pl 28\big)\Gamma_k \|\hat{u}_{sk+j} \|_2^2  \nonumber \\
&\leq \eps (34s \pl 53) \Gamma_k \Vert \hat{u}_{sk+j} \Vert_2^2.
\label{eq:wtf}
\end{align}

Given the above, we can also write the bound
\begin{equation}
\Vert \hat{\mathcal{Y}}_k \hat{w}'_{k,j} \Vert_2^2 \leq  \Vert \hat{\mathcal{Y}}_k \hat{w}'_{k,j} \Vert_2^2 + \hat{\alpha}^2_{sk+j} \leq \big(1+ \eps(34s\pl 53) \Gamma_k \big) \|\hat{u}_{sk+j} \|_2^2,
\label{eq:Yw}
\end{equation}
and using~\eqref{eq:cbound},~\eqref{eq:betahat}, and~\eqref{eq:deltabeta},
\begin{align*}
| \hat{\beta}_{sk+j+1} | &\leq \Vert \hat{\mathcal{Y}}_k \hat{w}'_{k,j} \Vert_2 \left( 1+ \eps + \frac{\delta d}{2\Vert \hat{\mathcal{Y}}_k \hat{w}'_{k,j} \Vert_2^2} \right)\\
&\phantom{=}\lesssim \big(1+ \eps(34s\pl 53) \Gamma_k \big) \|\hat{u}_{sk+j} \|_2 \left( 1+ \eps + \frac{\delta d}{2\Vert \hat{\mathcal{Y}}_k \hat{w}'_{k,j} \Vert_2^2} \right)\\
&\phantom{=}\leq \big(1+\eps + (\eps/2)(2s\pl 3) \Gamma_k + \eps(34s\pl 53) \Gamma_k \big)\|\hat{u}_{sk+j} \|_2.
\end{align*}
Combining terms, the above becomes
%
%
%
\begin{align}
\vert \hat{\beta}_{sk+j+1} \vert \leq& 
 \big(1 + \eps(35s\pl 56)\Gamma_k \big) \big\Vert \hat{u}_{sk+j} \big\Vert_2.
\label{eq:betabound}
\end{align}


Now, rearranging~\eqref{eq:hatvcoeff}, we can write
\begin{equation*}
\hat\beta_{sk+j+1}\hat{v}'_{k,j+1} = \hat{w}'_{k,j} + \delta \tilde{w}'_{k,j},
\end{equation*}
and premultiplying by $\hat{\mathcal{Y}}_k$, we obtain
\begin{equation*}
\hat\beta_{sk+j+1}\hat{\mathcal{Y}}_k\hat{v}'_{k,j+1} = \hat{\mathcal{Y}}_k\hat{w}'_{k,j} + \hat{\mathcal{Y}}_k\delta \tilde{w}'_{k,j}.
\end{equation*}
Using~\eqref{eq:bchangev}, this can be written
\begin{equation*}
\hat\beta_{sk+j+1}(\hat{v}_{sk+j+1} - \delta \hat{\mathcal{Y}}_{k,v_{j+1}}\hat{v}'_{k,j+1}) = \hat{\mathcal{Y}}_k\hat{w}'_{k,j} + \hat{\mathcal{Y}}_k\delta \tilde{w}'_{k,j}.
\end{equation*}
Rearranging and using~\eqref{eq:hatvcoeff}, 
\begin{align}
\hat\beta_{sk+j+1}\hat{v}_{sk+j+1}  &= \hat{\mathcal{Y}}_k\hat{w}'_{k,j} + \hat{\mathcal{Y}}_k\delta \tilde{w}'_{k,j}+  \hat\beta_{sk+j+1}\delta \hat{\mathcal{Y}}_{k,v_{j+1}}\hat{v}'_{k,j+1}\nonumber \\
&= \hat{\mathcal{Y}}_k\hat{w}'_{k,j} + \hat{\mathcal{Y}}_k\delta \tilde{w}'_{k,j}+  \delta \hat{\mathcal{Y}}_{k,v_{j+1}}(\hat{w}'_{k,j} + \delta \tilde{w}'_{k,j}) \nonumber \\
&= \hat{\mathcal{Y}}_k\hat{w}'_{k,j} + \hat{\mathcal{Y}}_k\delta \tilde{w}'_{k,j}+  \delta \hat{\mathcal{Y}}_{k,v_{j+1}}\hat{w}'_{k,j}\nonumber \\
&\equiv \hat{\mathcal{Y}}_k\hat{w}'_{k,j} + \delta w_{sk+j}, 
\label{eq:betav}
\end{align}
where $\delta w_{sk+j}=\hat{\mathcal{Y}}_k\delta \tilde{w}'_{k,j}+  \delta \hat{\mathcal{Y}}_{k,v_{j+1}}\hat{w}'_{k,j}$. Using Lemma~\ref{lem:condnum} and bounds in~\eqref{eq:deltavp}, ~\eqref{eq:bchangev}, and~\eqref{eq:Yw},
\begin{align}
\Vert \delta w_{sk+j} \Vert_2 &\leq \eps \| |\hat{\mathcal{Y}}_k | |\hat{w}'_{k,j}| \|_2 + \eps (2s\pl 2)  \| |\hat{\mathcal{Y}}_k | |\hat{w}'_{k,j}| \|_2 \nonumber \\
&\leq \eps (2s\pl 3)\Gamma_k \|\hat{\mathcal{Y}}_k \hat{w}'_{k,j} \|_2  \nonumber \\
&\leq \eps (2s\pl 3)\Gamma_k \Vert \hat{u}_{sk+j} \Vert_2.
\label{eq:deltwvec}
\end{align}
%
%

We premultiply~\eqref{eq:betav} by $\hat{v}_{sk+j}^T$ and use~\eqref{eq:bchangev},~\eqref{eq:bchangeu},~\eqref{eq:alphaexp}, and~\eqref{eq:deltaw}  to obtain 
\begin{align*}
\hat{\beta}_{sk+j+1}\hat{v}_{sk+j}^T\hat{v}_{sk+j+1} &= \hat{v}^T_{sk+j} ( \hat{\mathcal{Y}}_k\hat{w}'_{k,j} +\delta {w}_{sk+j} )\\
&= \hat{v}^T_{sk+j} ( \hat{\mathcal{Y}}_k\hat{u}'_{k,j}-\hat{\alpha}_{sk+j} \hat{\mathcal{Y}}_k\hat{v}'_{k,j} - \hat{\mathcal{Y}}_k\delta \hat{w}'_{k,j} ) + \hat{v}^T_{sk+j}\delta {w}_{sk+j} \\
&= \hat{v}^T_{sk+j} \big( \hat{\mathcal{Y}}_k \hat{u}'_{k,j} - \hat{\alpha}_{sk+j} \hat{\mathcal{Y}}_k\hat{v}'_{k,j}\big)- \hat{v}^T_{sk+j}  \big(\hat{\mathcal{Y}}_k\delta w'_{k,j} - \delta w_{sk+j} \big) \\
&= \hat{v}^T_{sk+j} \big( (\hat{u}_{sk+j} -  \delta \hat{\mathcal{Y}}_{k,u_j} \hat{u}'_{k,j})  - \hat{\alpha}_{sk+j} (\hat{v}_{sk+j} -  \delta \hat{\mathcal{Y}}_{k,v_j} \hat{v}'_{k,j}) \big) \\
&\phantom{=}-\hat{v}^T_{sk+j} (\hat{\mathcal{Y}}_k \delta w'_{k,j} -\delta {w}_{sk+j} )\\
&= \hat{v}^T_{sk+j} \hat{u}_{sk+j} - \hat{\alpha}_{sk+j} \hat{v}^T_{sk+j} \hat{v}_{sk+j} \\
&\phantom{=}-\hat{v}^T_{sk+j} ( \delta \hat{\mathcal{Y}}_{k,u_j} \hat{u}'_{k,j} -  \hat{\alpha}_{sk+j}\delta \hat{\mathcal{Y}}_{k,v_j} \hat{v}'_{k,j} +  \hat{\mathcal{Y}}_k \delta w'_{k,j} -\delta {w}_{sk+j} )\\
&= (\hat{\alpha}_{sk+j} - \delta\hat{\alpha}_{sk+j} ) - \hat{\alpha}_{sk+j} \| \hat{v}_{sk+j}\|_2^2 \\
&\phantom{=}-\hat{v}^T_{sk+j} ( \delta \hat{\mathcal{Y}}_{k,u_j} \hat{u}'_{k,j} -  \hat{\alpha}_{sk+j}\delta \hat{\mathcal{Y}}_{k,v_j} \hat{v}'_{k,j} +  \hat{\mathcal{Y}}_k \delta w'_{k,j} -\delta {w}_{sk+j} )\\
&=  -\delta\hat{\alpha}_{sk+j}  - \hat{\alpha}_{sk+j} (\| \hat{v}_{sk+j}\|_2^2-1) \\
&\phantom{=}-\hat{v}^T_{sk+j} ( \delta \hat{\mathcal{Y}}_{k,u_j} \hat{u}'_{k,j} -  \hat{\alpha}_{sk+j}\delta \hat{\mathcal{Y}}_{k,v_j} \hat{v}'_{k,j} +  \hat{\mathcal{Y}}_k \delta w'_{k,j} -\delta {w}_{sk+j} ),
%
%
%
%
%
%
%
%
%
%
\end{align*}
%
%
and using Lemma~\ref{lem:condnum} and bounds in~\eqref{mainthm:norm},~\eqref{eq:hatvcoeff},~\eqref{eq:bchangev},~\eqref{eq:bchangeu},~\eqref{eq:deltaalpha},~\eqref{eq:alphabound},~\eqref{eq:deltaw}, and~\eqref{eq:deltwvec}, we can write the bound
\begin{align}
\Big\vert \hat\beta_{sk+j+1}\cdot\hat{v}_{sk+j}^T \hat{v}_{sk+j+1}  \Big\vert &\leq | \delta \hat{\alpha}_{sk+j} | + | \hat{\alpha}_{sk+j} | |\hat{v}^T_{sk+j}\hat{v}_{sk+j} -1| \nonumber \\
&\phantom{\leq} + \|\hat{v}_{sk+j} \|_2 \big( \| \lvvert \delta \hat{\mathcal{Y}}_{k,u_j} \vert  \vert\hat{u}'_{k,j} \rvvert \|_2 + | \hat{\alpha}_{sk+j} |  \| \lvvert  \delta \hat{\mathcal{Y}}_{k,v_j} \vert \vert \hat{v}'_{k,j} \rvvert  \|_2  \big) \nonumber \\
&\phantom{\leq} + \|\hat{v}_{sk+j} \|_2 \big(\| \lvvert \hat{\mathcal{Y}}_k \vert \vert \delta w'_{k,j} \rvvert \|_2 + \|\delta w_{sk+j} \|_2 \big) \nonumber \\
%
%
&\leq  \eps (18s\pl 28) \Gamma_k \| \hat{u}_{sk+j}\|_2. \label{eq:orthbound1}
\end{align}

This is a start toward proving~\eqref{mainthm:orth}. We will return to the above bound once we later prove a bound on $\Vert \hat{u}_{sk+j} \Vert_2$. Our next step is to analyze the error in each column of the finite precision $s$-step Lanczos recurrence. First, we note that we can write the error in computing the $s$-step bases (line~\ref{slanczos:basis} in Algorithm~\ref{alg:sLanczos}) by
\begin{equation}
A \underline{\hat{\mathcal{Y}}}_k = \hat{\mathcal{Y}}_k \mathcal{B}_k + \delta E_k
\label{eq:berror}
\end{equation}
where $\underline{\hat{\mathcal{Y}}}_k = \big[\hat{\mathcal{V}}_k [I_s,0_{s,1}]^T, 0_{n,1},\hat{\mathcal{U}}_k [I_s,0_{s,1}]^T, 0_{n,1} \big]$. It can be shown (see, e.g.,~\cite{carson2014residual}) that if the basis is computed in the usual way by repeated SpMVs,
\color{black}
\begin{equation}
\vert \delta E_k \vert \leq \eps \big( (3\pl N)\vert A \vert \vert \underline{\hat{\mathcal{Y}}}_k \vert + 4 \vert \hat{\mathcal{Y}}_k \vert \vert \mathcal{B}_k \vert \big),
\label{eq:Ebound}
\end{equation} 
where $N$ is the maximum number of nonzeros per row over all rows of $A$. 
\color{black}

In finite precision, line~\ref{slanczos:ucoeff} in Algorithm~\ref{alg:sLanczos} is computed as
\begin{equation}
\hat{u}'_{k,j}\hspace{-1pt} =\hspace{-1pt} \mathcal{B}_k \hat{v}'_{k,j} \mn \hat{\beta}_{sk+j}\hat{v}'_{k,j-1} \pl \delta{u}'_{k,j}, \hspace{1mm} \vert \delta{u}'_{k,j} \vert \hspace{-1pt}\leq \hspace{-1pt} \eps \big((2s\pl 3)\vert \mathcal{B}_k \vert \vert \hat{v}'_{k,j} \vert \pl 2 \vert \hat{\beta}_{sk+j} \hat{v}'_{k,j-1}\vert \big),
\label{eq:deltau}
\end{equation}
and then, with Lemma~\ref{lem:condnum},~\eqref{eq:bchangev},~\eqref{eq:bchangeu},~\eqref{eq:berror}, and~\eqref{eq:deltau}, we can write
\begin{align}
 \hat{u}_{sk+j} &= (\hat{\mathcal{Y}}_k + \delta \hat{\mathcal{Y}}_{k,u_j})\hat{u}'_{k,j} \nonumber \\
 &=(\hat{\mathcal{Y}}_k + \delta \hat{\mathcal{Y}}_{k,u_j}) (\mathcal{B}_k \hat{v}'_{k,j} -\hat{\beta}_{sk+j}\hat{v}'_{k,j-1} + \delta{u}'_{k,j})\nonumber \\
&= \hat{\mathcal{Y}}_k \mathcal{B}_k \hat{v}'_{k,j} - \hat{\beta}_{sk+j} \hat{\mathcal{Y}}_k \hat{v}'_{k,j-1} + \hat{\mathcal{Y}}_k \delta u'_{k,j} + \delta \hat{\mathcal{Y}}_{k,u_j} \mathcal{B}_k \hat{v}'_{k,j} - \hat{\beta}_{sk+j} \delta \hat{\mathcal{Y}}_{k,u_j} \hat{v}'_{k,j-1}\nonumber \\
&= (A\underline{\hat{\mathcal{Y}}}_k - \delta E_k )\hat{v}'_{k,j} - \hat{\beta}_{sk+j} (\hat{v}_{sk+j-1} - \delta \hat{\mathcal{Y}}_{k,v_{j-1}}\hat{v}'_{k,j-1}) + \hat{\mathcal{Y}}_k \delta u'_{k,j}  \nonumber \\
&\phantom{=} + \delta \hat{\mathcal{Y}}_{k,u_j} \mathcal{B}_k \hat{v}'_{k,j} - \hat{\beta}_{sk+j} \delta \hat{\mathcal{Y}}_{k,u_j} \hat{v}'_{k,j-1}\nonumber \\
&= A\underline{\hat{\mathcal{Y}}}_k\hat{v}'_{k,j} - \delta E_k \hat{v}'_{k,j} - \hat{\beta}_{sk+j} \hat{v}_{sk+j-1} + \hat{\beta}_{sk+j} \delta \hat{\mathcal{Y}}_{k,v_{j-1}}\hat{v}'_{k,j-1} + \hat{\mathcal{Y}}_k \delta u'_{k,j}  \nonumber \\
&\phantom{=} + \delta \hat{\mathcal{Y}}_{k,u_j} \mathcal{B}_k \hat{v}'_{k,j} - \hat{\beta}_{sk+j} \delta \hat{\mathcal{Y}}_{k,u_j} \hat{v}'_{k,j-1}\nonumber \\
&= A ( \hat{v}_{sk+j} - \delta{\hat{\mathcal{Y}}}_{k,v_j} \hat{v}'_{k,j}) - \delta E_k \hat{v}'_{k,j} - \hat{\beta}_{sk+j} \hat{v}_{sk+j-1} + \hat{\beta}_{sk+j} \delta \hat{\mathcal{Y}}_{k,v_{j-1}}\hat{v}'_{k,j-1}  \nonumber \\
&\phantom{=} + \hat{\mathcal{Y}}_k \delta u'_{k,j} + \delta \hat{\mathcal{Y}}_{k,u_j} \mathcal{B}_k \hat{v}'_{k,j} - \hat{\beta}_{sk+j} \delta \hat{\mathcal{Y}}_{k,u_j} \hat{v}'_{k,j-1}\nonumber \\
&= A\hat{v}_{sk+j} - A\delta \hat{\mathcal{Y}}_{k,v_j}\hat{v}'_{k,j}  - \delta E_k \hat{v}'_{k,j} - \hat{\beta}_{sk+j} \hat{v}_{sk+j-1} + \hat{\beta}_{sk+j} \delta \hat{\mathcal{Y}}_{k,v_{j-1}}\hat{v}'_{k,{j-1}} \nonumber \\
&\phantom{=} + \hat{\mathcal{Y}}_k \delta u'_{k,j} + \delta \hat{\mathcal{Y}}_{k,u_j} \mathcal{B}_k \hat{v}'_{k,j} - \hat{\beta}_{sk+j} \delta \hat{\mathcal{Y}}_{k,u_j} \hat{v}'_{k,j-1}\nonumber \\
&\equiv A\hat{v}_{sk+j}- \hat{\beta}_{sk+j} \hat{v}_{sk+j-1} +\delta u_{sk+j},
\label{eq:uhat}
\end{align}
where 
\begin{equation*}
\delta u_{sk+j} = \hat{\mathcal{Y}}_k \delta u'_{k,j} - (A\delta \hat{\mathcal{Y}}_{k,v_j} - \delta \hat{\mathcal{Y}}_{k,u_j} \mathcal{B}_k + \delta E_k) \hat{v}'_{k,j} + \hat{\beta}_{sk+j} (\delta \hat{\mathcal{Y}}_{k,v_{j-1}} - \delta \hat{\mathcal{Y}}_{k,u_j}  )\hat{v}'_{k,j-1}.
\end{equation*}
%
Using the bounds in~\eqref{eq:bchangev},~\eqref{eq:bchangeu},~\eqref{eq:betabound},~\eqref{eq:Ebound}, and~\eqref{eq:deltau}, we can write 
\color{black}
\begin{align*}
\vert \delta u_{sk+j} \vert &\leq \eps\big( (2s\pl 3) |\hat{\mathcal{Y}}_k | |\mathcal{B}_k | |\hat{v}'_{k,j} | + 2 |\hat{\beta}_{sk+j}| |\hat{\mathcal{Y}}_k ||\hat{v}'_{k,j-1} |    \big)\\
&\phantom{\leq} + \eps (2s\pl 2) |A||\hat{\mathcal{Y}}_k | |\hat{v}'_{k,j} | + \eps (2s\pl 2)|\hat{\mathcal{Y}}_k | |\mathcal{B}_k | |\hat{v}'_{k,j} |\\
&\phantom{\leq} + \eps \big( (3\pl N) |A||\hat{\mathcal{Y}}_k | |\hat{v}'_{k,j} | + 4|\hat{\mathcal{Y}}_k | |\mathcal{B}_k | |\hat{v}'_{k,j} | \big)\\
&\phantom{\leq} + 2\eps (2s\pl 2) |\hat{\beta}_{sk+j}| |\hat{\mathcal{Y}}_k ||\hat{v}'_{k,j-1} | \\
&\leq \eps (N\pl 2s\pl 5)|A||\hat{\mathcal{Y}}_k | |\hat{v}'_{k,j} | + \eps(4s\pl 9)|\hat{\mathcal{Y}}_k | |\mathcal{B}_k | |\hat{v}'_{k,j} |  \\
&\phantom{\leq} + \eps(4s\pl 6) \big( 1+\eps (35s\pl 56)\Gamma_k \big) \| \hat{u}_{sk+j-1} \|_2 \cdot |\hat{\mathcal{Y}}_k ||\hat{v}'_{k,j-1} |. 
\end{align*}
\color{black}
and from this we obtain
\color{black}
\begin{align*}
\Vert \delta u_{sk+j} \Vert_2 &\leq \eps (N\pl 2s\pl 5) \hspace{2pt} \|\lvvert A \rvvert  \|_2 \hspace{2pt} \|\lvvert \hat{\mathcal{Y}}_k  \rvvert  \|_2 \hspace{2pt} \| \hat{v}'_{k,j} \|_2 + \eps (4s\pl 9) \hspace{2pt} \|\lvvert \hat{\mathcal{Y}}_k  \rvvert  \|_2 \hspace{2pt} \|\lvvert {\mathcal{B}}_k  \rvvert  \|_2 \hspace{2pt} \| \hat{v}'_{k,j} \|_2 \\
&\phantom{\leq}+ \eps (4s\pl 6) \hspace{2pt} \|\lvvert \hat{\mathcal{Y}}_k  \rvvert  \|_2 \hspace{2pt}  \| \hat{v}'_{k,j-1} \|_2 \hspace{2pt} \| \hat{u}_{sk+j-1} \|_2\\
&\leq \eps(N\pl 2s\pl 5) \Gamma_k \hspace{2pt} \|\lvvert A \rvvert  \|_2 \| \vhskj\|_2 \hspace{2pt} + \eps (4s\pl 9) \Gamma_k  \hspace{2pt} \|\lvvert \mathcal{B}_k \rvvert  \|_2  \hspace{2pt}\| \vhskj\|_2 \\
&\phantom{\leq}+ \eps (4s\pl 6) \Gamma_k\hspace{2pt} \| \hat{u}_{sk+j-1} \|_2 \| \hat{v}_{sk+j-1} \|_2 \hspace{2pt}  \\
&\leq  \eps(N\pl 2s\pl 5) \Gamma_k \hspace{2pt} \|\lvvert A \rvvert  \|_2  \hspace{2pt} + \eps (4s\pl 9) \Gamma_k  \hspace{2pt} \|\lvvert \mathcal{B}_k \rvvert  \|_2  + \eps (4s\pl 6) \Gamma_k\hspace{2pt} \| \hat{u}_{sk+j-1} \|_2 .  \\
\end{align*}
\color{black}
We will now introduce and make use of the quantities
\color{black}
$\sigma=\|A\|_2$, $\theta\equiv \| \lvvert A\rvvert \|_2/\sigma$ and $\tau_k \equiv \| \lvvert \mathcal{B}_k\rvvert \|_2/\sigma$. 
\color{black}
Note that the quantity $\| \lvvert \mathcal{B}_k \rvvert \|_2$ is in some sense controlled by the user, and as previously stated, for the usual basis choices, including monomial, Newton, or Chebyshev bases, it should be the case that $\| \lvvert \mathcal{B}_k \rvvert \|_2 \lesssim\| \lvvert A\rvvert \|_2$. Using these quantities, the bound above can be written
\color{black}
\begin{equation}
\big\Vert \delta {u}_{sk+j} \big\Vert_2 \leq  \eps \Big( \big( (N\pl 2s \pl 5) \theta + (4s\pl 9) \tau_k \big)\sigma +(4s\pl 6) \| \hat{u}_{sk+j-1} \|_2  \Big)\Gamma_k.
\label{eq:ubound}
\end{equation}
\color{black}
%

%
Manipulating~\eqref{eq:betav}, and using~\eqref{eq:bchangev},~\eqref{eq:bchangeu}, and~\eqref{eq:deltaw}, we have 
\begin{align*}
\hat{\beta}_{sk+j+1} \hat{v}_{sk+j+1} &= \hat{\mathcal{Y}}_k \hat{w}'_{k,j} + \delta w_{sk+j} \\
&= \hat{\mathcal{Y}}_k \hat{u}'_{k,j} - \hat{\alpha}_{sk+j} \hat{\mathcal{Y}}_k \hat{v}'_{k,j} - \hat{\mathcal{Y}}_k \delta w'_{k,j} + \delta w_{sk+j} \\
&= (\hat{u}_{sk+j} \mn \delta \hat{\mathcal{Y}}_{k,u_j} \hat{u}'_{k,j}) \mn \hat{\alpha}_{sk+j} (\hat{v}_{sk+j} \mn \delta \hat{\mathcal{Y}}_{k,v_j} \hat{v}'_{k,j}) \mn \hat{\mathcal{Y}}_k \delta w'_{k,j} \pl \delta w_{sk+j}\\
&= \hat{u}_{sk+j} \mn \hat{\alpha}_{sk+j} \hat{v}_{sk+j}  \mn \delta \hat{\mathcal{Y}}_{k,u_j} \hat{u}'_{k,j} \pl \hat{\alpha}_{sk+j} \delta \hat{\mathcal{Y}}_{k,v_j} \hat{v}'_{k,j} \mn \hat{\mathcal{Y}}_k \delta w'_{k,j} \\
&\phantom{=} + \delta w_{sk+j},
\end{align*}
and substituting in the expression for $\hat{u}_{sk+j}$ in~\eqref{eq:uhat} on the right, we obtain
\begin{equation}
\hat{\beta}_{sk+j+1} \hat{v}_{sk+j+1} \equiv A\hat{v}_{sk+j}- \hat{\alpha}_{sk+j} \hat{v}_{sk+j} -\hat{\beta}_{sk+j} \hat{v}_{sk+j-1} + \delta \hat{v}_{sk+j},
\label{eq:betavhat}
\end{equation}
where
\begin{equation*}
\delta \hat{v}_{sk+j} = \delta u_{sk+j}  - \delta \hat{\mathcal{Y}}_{k,u_j} \hat{u}'_{k,j} + \hat{\alpha}_{sk+j} \delta \hat{\mathcal{Y}}_{k,v_j} \hat{v}'_{k,j} - \hat{\mathcal{Y}}_k \delta w'_{k,j} + \delta w_{sk+j}.
\end{equation*}
%
%
%
From this we can write the componentwise bound
\begin{align*}
\vert \delta \hat{v}_{sk+j} \vert &\leq |\delta {u}_{sk+j} | + |\delta \hat{\mathcal{Y}}_{k,u_j}|\hspace{2pt} |\hat{u}'_{k,j} | + |\hat{\alpha}_{sk+j} | \hspace{2pt} |\delta \hat{\mathcal{Y}}_{k,v_j} | \hspace{2pt} | \hat{v}'_{k,j} | + |\hat{\mathcal{Y}}_k | \hspace{2pt} |\delta w'_{k,j} | + |\delta w_{sk+j}|,
\end{align*}
and using Lemma~\ref{lem:condnum},~\eqref{eq:deltavp},~\eqref{eq:bchangev},~\eqref{eq:bchangeu},~\eqref{eq:alphabound},~\eqref{eq:deltaw}, and~\eqref{eq:deltwvec} we obtain 
\begin{align*}
\Vert \delta \hat{v}_{sk+j} \Vert_2 &\leq 
\Vert \delta u_{sk+j} \Vert_2 + \eps(2s\pl 2) \Gamma_k \| \uhskj\|_2   + \eps (2s\pl 2) \Gamma_k \| \hat{u}_{sk+j} \|_2  \\
&\phantom{\leq} + \eps \Gamma_k \|\uhskj\|_2 
+ 2\eps\Gamma_k \|\uhskj\|_2 + \epsilon (2s\pl 3) \Gamma_k \| \hat{u}_{sk+j} \|_2\\
&\leq \| \delta {u}_{sk+j} \|_2 + \epsilon (6s\pl 10) \Gamma_k \| \hat{u}_{sk+j} \|_2.
%
\end{align*}
Using~\eqref{eq:ubound}, this gives the bound
\color{black}
\begin{multline}
\|\delta\hat{v}_{sk+j} \|_2 \leq \\
\eps \Big(\big((N\pl 2s\pl 5)\theta \pl (4s\pl 9) \tau_k \big)\sigma \pl (6s\pl 10) \| \hat{u}_{sk+j}\|_2  \pl (4s\pl 6) \| \hat{u}_{sk+j-1} \|_2  \Big) \Gamma_k.
\label{eq:deltavnorm}
\end{multline}
\color{black}

We now have everything we need to write the mixed precision $s$-step Lanczos recurrence in matrix form. Let
\begin{align*}
\hat{T}_{sk+j} & = \left[\begin{array}{cccc}
\hat\alpha_1 & \hat{\beta}_2 & & \\
\hat{\beta}_2 & \ddots & \ddots & \\
 & \ddots & \ddots & \hat\beta_{sk+j} \\
 & & \hat\beta_{sk+j} & \hat\alpha_{sk+j} 
\end{array}\right]
\end{align*}
be the tridiagonal Jacobi matrix composed of the computed $\hat{\alpha}$'s and $\hat{\beta}$'s, and let $\hat{V}_{sk+j} = [\hat{v}_1,\hat{v}_2,\ldots, \hat{v}_{sk+j}]$ and $\delta\hat{V}_{sk+j} = [\delta\hat{v}_1,\delta\hat{v}_2,\ldots, \delta\hat{v}_{sk+j}]$. Note that $\hat{T}_{sk+j}$ has dimension $(sk+j)$-by-$(sk+j)$, and $\hat{V}_{sk+j}$ and $\delta \hat{V}_{sk+j}$ have dimension $n$-by-$(sk+j)$. Then~\eqref{eq:betavhat} in matrix form gives
\begin{equation}
A\hat{V}_{sk+j} = \hat{V}_{sk+j} \hat{T}_{sk+j}+ \hat\beta_{sk+j+1}\hat{v}_{sk+j+1}e^T_{sk+j}  - \delta \hat{V}_{sk+j}.
\label{eq:matrixrecurrence}
\end{equation}
Thus~\eqref{eq:deltavnorm} gives a bound on the error in the columns of the mixed precision $s$-step Lanczos recurrence. Again, we will return to~\eqref{eq:deltavnorm} to prove~\eqref{mainthm:deltav} once we bound $\Vert \hat{u}_{sk+j} \Vert_2$. 

We now turn to the loss of orthogonality in the vectors $\hat{v}_1,\ldots,\hat{v}_{sk+j+1}$. 
Let ${R}_{sk+j}$ be the strictly upper triangular matrix of dimension $(sk+j)$-by-$(sk+j)$ with elements $\rho_{i,j}$, for $i,j\in \{1,\ldots,sk+j\}$, such that 
\begin{equation*}
\hat{V}_{sk+j}^T \hat{V}_{sk+j} = {R}_{sk+j}^T + \text{diag}(\hat{V}_{sk+j}^T \hat{V}_{sk+j}) + {R}_{sk+j}.
\end{equation*}
We also define $\rho_{sk+j,sk+j+1}\equiv \hat{v}_{sk+j}^T \hat{v}_{sk+j+1}$, noting that $\rho_{sk+j,sk+j+1}$ is not an element of $R_{sk+j}$ but would be an element of $R_{sk+j+1}$). 
Multiplying~\eqref{eq:matrixrecurrence} on the left by $\hat{V}_{sk+j}^T$, we get
\begin{equation*}
\hat{V}_{sk+j}^T A\hat{V}_{sk+j} = \hat{V}_{sk+j}^T\hat{V}_{sk+j} \hat{T}_{sk+j}+ \hat\beta_{sk+j+1}\hat{V}_{sk+j}^T\hat{v}_{sk+j+1}e^T_{sk+j}  - \hat{V}_{sk+j}^T\delta \hat{V}_{sk+j}.
\end{equation*}
Using symmetry, we can equate the right hand side with its own transpose to obtain
\begin{align*}
\hat{T}_{sk+j} (R_{sk+j}^T + R_{sk+j}) &- (R_{sk+j}^T + R_{sk+j})\hat{T}_{sk+j} \\
&=\hat\beta_{sk+j+1}(\hat{V}_{sk+j}^T\hat{v}_{sk+j+1}e^T_{sk+j}-e_{sk+j}\hat{v}^T_{sk+j+1}\hat{V}_{sk+j} )\\
  &\phantom{=}+ \hat{V}_{sk+j}^T\delta \hat{V}_{sk+j} - \delta \hat{V}_{sk+j}^T\hat{V}_{sk+j} \\
  &\phantom{=}+ \text{diag}(\hat{V}_{sk+j}^T \hat{V}_{sk+j})\cdot \hat{T}_{sk+j} -\hat{T}_{sk+j} \cdot \text{diag}(\hat{V}_{sk+j}^T \hat{V}_{sk+j}).
\end{align*}

Now, let $M_{sk+j}\equiv \hat{T}_{sk+j}R_{sk+j} - R_{sk+j}\hat{T}_{sk+j}$, which is upper triangular and has dimension $(sk+j)$-by-$(sk+j)$. Then the left-hand side above can be written as $M_{sk+j}-M_{sk+j}^T$, and we can equate the strictly upper triangular part of $M_{sk+j}$ with the strictly upper triangular part of the right-hand side above. The diagonal elements can be obtained from the definition $M_{sk+j}\equiv \hat{T}_{sk+j}R_{sk+j} - R_{sk+j}\hat{T}_{sk+j}$:
\begin{align*}
m_{1,1} =& -\hat\beta_{2}\rho_{1,2}, \qquad m_{sk+j,sk+j}= \hat\beta_{sk+j}\rho_{sk+j-1,sk+j},\quad\text{and}\\
m_{i,i} =& \hat\beta_{i}\rho_{i-1,i} -\hat\beta_{i+1}\rho_{i,i+1},\quad\text{for}\quad  i \in \{2,\ldots,sk+j-1\}.
\end{align*}
Therefore, we can write
\begin{equation*}
M_{sk+j} = \hat{T}_{sk+j}R_{sk+j} - R_{sk+j}\hat{T}_{sk+j} = \hat\beta_{sk+j+1}\hat{V}_{sk+j}^T \hat{v}_{sk+j+1}e_{sk+j}^T + H_{sk+j},
\end{equation*}
where $H_{sk+j}$ has elements satisfying
\begin{equation}\label{eq:etas}
\begin{split}
\eta_{1,1} = & - \hat\beta_{2}\rho_{1,2},\\
\eta_{i,i} = & \hat\beta_{i} \rho_{i-1,i} - \hat\beta_{i+1} \rho_{i,i+1}, \quad\text{for}\quad i\in\{2,\ldots,sk\pl j\},\\
\eta_{i-1,i}= & \hat{v}_{i-1}^T \delta\hat{v}_{i} - \delta\hat{v}_{i-1}^T \hat{v}_{i} + \hat{\beta}_{i} (\hat{v}_{i-1}^T\hat{v}_{i-1} - \hat{v}_{i}^T \hat{v}_{i}),\hspace{2mm}\text{and} \\
\eta_{\ell,i} = & \hat{v}_{\ell}^T \delta \hat{v}_{i} - \delta\hat{v}_{\ell}^T \hat{v}_{i}, \quad \text{for}\quad \ell\in\{1,\ldots,i-2\}.
\end{split}
\end{equation}
To simplify notation, we introduce the quantities
\begin{equation*}
\bar{u}_{sk+j} = \hspace{-1mm}\max_{i\in\{1,\ldots,sk+j\}}\| \hat{u}_{i} \|_2, \quad \bar{\Gamma}_k =\hspace{-1mm} \max_{i\in\{0,\ldots,k\} } \Gamma_i ,\quad\text{and} \quad \bar{\tau}_k =\hspace{-1mm} \max_{i\in\{0,\ldots,k\} }  \tau_i .
\end{equation*}
Using this notation and~\eqref{mainthm:norm},~\eqref{eq:betabound},~\eqref{eq:orthbound1}, and~\eqref{eq:deltavnorm}, the quantities in~\eqref{eq:etas} can be bounded by
\begin{equation}\label{eq:etasbounds}
\begin{split}
\vert\eta_{1,1}\vert \leq & \tp \eps (18s\pl 28)\op \bar\Gamma_k\op \bar{u}_{sk+j}, \hspace{2mm}\text{and, for } i\in\{2,\ldots,sk\pl j\}, \\
\vert \eta_{i,i} \vert \leq & \tp 2\eps (18s\pl 28)\op \bar\Gamma_k \op \bar{u}_{sk+j}, \\
\color{black}\vert \eta_{i-1,i} \vert \leq & \color{black}\tp 2\eps \Big( \hspace{-2pt} \big( (N\pl 2s\pl 5)\theta \pl (4s \pl 9)\bar{\tau}_k \big) \sigma \pl (16s\pl 27)\bar{u}_{sk+j}  \hspace{-2pt}\Big) \bar{\Gamma}_k ,\\
%
%
\color{black}\vert \eta_{\ell,i} \vert  \leq & \color{black}\tp 2\eps \Big( \hspace{-2pt} \big( (N\pl 2s\pl 5)\theta \pl (4s \pl 9)\bar{\tau}_k \big) \sigma \pl (10s\pl 16)\bar{u}_{sk+j}  \hspace{-2pt}\Big) \bar{\Gamma}_k,
\end{split}
\end{equation}
for $\ell \in\{1,\ldots,i\mn 2\}$.

This is start toward proving~\eqref{mainthm:etas}, but to proceed we first need a bound on $\Vert \hat{u}_{sk+j} \Vert_2$.
We must first find a bound for $\vert \rho_{sk+j-2,sk+j} \vert=\vert \hat{v}_{sk+j}^T\hat{v}_{sk+j-2} \vert$. From the definition of $M_{sk+j}$, we know the ($1,2$) element of $M_{sk+j}$ is 
\[
\hat\alpha_1 \rho_{1,2} - \hat\alpha_2 \rho_{1,2} - \hat\beta_3 \rho_{1,3}=\eta_{1,2},
\]
and for $i>2$, the ($i \mn 1, i$) element is 
\begin{equation*}
\hat\beta_{i-1}\rho_{i-2,i} + (\hat\alpha_{i-1} -\hat\alpha_{i})\rho_{i-1,i} - \hat\beta_{i+1}\rho_{i-1,i+1}=\eta_{i-1,i}.
\end{equation*}
Then, defining
\[
\xi_{i} \equiv (\hat\alpha_{i-1}-\hat\alpha_{i})\hat\beta_{i}\rho_{i-1,i} - \hat\beta_{i}\eta_{i-1,i}
\]
for $i\in\{2,\ldots,sk\pl j\}$, we have
\[
\hat\beta_{i} \hat\beta_{i+1} \rho_{i-1,i+1} = \hat\beta_{i-1} \hat\beta_{i} \rho_{i-2,i} + \xi_{i} = \xi_{i}+\xi_{i-1}+\ldots+\xi_{2}. 
\]
This, along with~\eqref{eq:alphabound},~\eqref{eq:betabound},~\eqref{eq:orthbound1}, and~\eqref{eq:etasbounds} gives
\begin{align}
\hat\beta_{sk+j} &\hat\beta_{sk+j+1} \vert \rho_{sk+j-1,sk+j+1} \vert = \hat\beta_{sk+j} \hat\beta_{sk+j+1} |\hat{v}^T_{sk+j-1} \hat{v}_{sk+j+1}  | \nonumber \\
\leq& \sum_{i=2}^{sk+j} | \xi_{i} | \leq \sum_{i=2}^{sk+j} (| \hat{\alpha}_{i-1} | + |\hat{\alpha}_{i}|)|\hat{\beta}_{i} \rho_{i-1,i}| + |\hat{\beta}_{i} | |\eta_{i-1,i}| \nonumber \\
%
\color{black}\leq& \color{black}2 \eps \sum_{i=2}^{sk+j} \Big( \big( (N\pl 2s\pl 5)\theta + (4s \pl 9)\bar{\tau}_k \big) \sigma + (34s\pl 55)\bar{u}_{sk+j} \Big) \bar{\Gamma}_k \bar{u}_{sk+j}\nonumber \\
\color{black}\leq& \color{black}2\eps(sk\pl j \mn 1) \Big( \big( (N\pl 2s\pl 5)\theta + (4s \pl 9)\bar{\tau}_k \big) \sigma + (34s\pl 55)\bar{u}_{sk+j} \Big) \bar{\Gamma}_k \bar{u}_{sk+j}.
\label{eq:vvjm2}
\end{align}

Rearranging~\eqref{eq:uhat} we have 
\begin{equation*}
\hat{u}_{sk+j} - \delta u_{sk+j} = A\hat{v}_{sk+j} -\hat{\beta}_{sk+j}\hat{v}_{sk+j-1},
\end{equation*}
and multiplying each side by its own transpose (and ignoring $\eps^2$ terms), we obtain
\begin{equation}
\hat{u}_{sk+j}^T \hat{u}_{sk+j} - 2\hat{u}_{sk+j}^T\delta u_{sk+j}   = \Vert A \hat{v}_{sk+j} \Vert_2^2 \pl \hat\beta_{sk+j}^2 \Vert \hat{v}_{sk+j-1} \Vert_2^2 \mn 2\hat\beta_{sk+j}\hat{v}_{sk+j}^T A \hat{v}_{sk+j-1}.
\label{eq:updu}
\end{equation}
Rearranging~\eqref{eq:betavhat} gives
\begin{equation*}
A\hat{v}_{sk+j-1}= \hat{\beta}_{sk+j} \hat{v}_{sk+j} + \hat{\alpha}_{sk+j-1}\hat{v}_{sk+j-1} + \hat{\beta}_{sk+j-1} \hat{v}_{sk+j-2} - \delta \hat{v}_{sk+j-1}, 
\end{equation*}
and premultiplying this expression by $\hat{\beta}_{sk+j} \hat{v}^T_{sk+j}$, we get
\begin{align}
\hat{\beta}_{sk+j}& \hat{v}_{sk+j}^T A \hat{v}_{sk+j-1} \nonumber \\
=& \hat{\beta}_{sk+j}\hat{v}_{sk+j}^T \big( \hat{\beta}_{sk+j}\hat{v}_{sk+j} + \hat{\alpha}_{sk+j-1}\hat{v}_{sk+j-1} + \hat{\beta}_{sk+j-1} \hat{v}_{sk+j-2} -\delta\hat{v}_{sk+j-1} \big) \nonumber \\
=& \hat{\beta}_{sk+j}^2 \Vert \hat{v}_{sk+j} \Vert_2^2 + \hat{\alpha}_{sk+j-1} (\hat{\beta}_{sk+j}\hat{v}_{sk+j}^T \hat{v}_{sk+j-1}) + \hat{\beta}_{sk+j}\hat{\beta}_{sk+j-1} \hat{v}_{sk+j}^T \hat{v}_{sk+j-2} \nonumber \\
&- \hat{\beta}_{sk+j} \hat{v}_{sk+j}^T \delta\hat{v}_{sk+j-1} \nonumber \\
\equiv& \hat{\beta}_{sk+j}^2 + \delta\hat{\beta}_{sk+j}, 
\label{eq:bvAv}
\end{align}
where, using bounds in~\eqref{mainthm:norm},~\eqref{eq:alphabound},~\eqref{eq:betabound},~\eqref{eq:orthbound1},~\eqref{eq:deltavnorm}, and~\eqref{eq:vvjm2},
\color{black}
\begin{align}
\vert \delta\hat{\beta}_{sk+j} \vert &\leq \eps  \Big( \big( (N\pl 2s\pl 5)\theta + (4s \pl 9)\bar{\tau}_k \big) \sigma + (34s\pl 55)\bar{u}_{sk+j} \Big) \bar{\Gamma}_k   \bar{u}_{sk+j} \nonumber \\
&\phantom{\leq} \pl 2\eps (sk\pl j\mn 1) \Big( \big( (N\pl 2s\pl 5)\theta + (4s \pl 9)\bar{\tau}_k \big) \sigma + (34s\pl 55)\bar{u}_{sk+j} \Big) \bar{\Gamma}_k  \bar{u}_{sk+j}  \nonumber \\
&\leq \eps \big(2(sk\pl j) \mn 1 \big) \Big( \big( (N\pl 2s\pl 5)\theta + (4s \pl 9)\bar{\tau}_k \big) \sigma + (34s\pl 55)\bar{u}_{sk+j} \Big) \bar{\Gamma}_k  \bar{u}_{sk+j}. \label{eq:betahatbound}
\end{align}
\color{black}
%
We note that in contrast to the uniform precision case, this bound contains only a factor of $\bar{\Gamma}_k$ rather than $\bar{\Gamma}_k^2$. 
Adding $2\hat{u}_{sk+j}^T\delta u_{sk+j}$ to both sides of~\eqref{eq:updu} and using the bound~\eqref{eq:bvAv}, 
\begin{align}
\Vert \hat{u}_{sk+j} \Vert_2^2 &= \Vert A \hat{v}_{sk+j} \Vert_2^2 + \hat\beta_{sk+j}^2 \big( \Vert \hat{v}_{sk+j-1} \Vert_2^2 -2 \big) -2\delta\hat\beta_{sk+j} + 2\hat{u}^T_{sk+j} \delta u_{sk+j} \nonumber \\
&\equiv \Vert A \hat{v}_{sk+j} \Vert_2^2 + \hat\beta_{sk+j}^2 \big( \Vert \hat{v}_{sk+j-1} \Vert_2^2 -2 \big) + \delta\tilde\beta_{sk+j},
\label{eq:uhat2}
\end{align}
where $\delta\tilde{\beta}_{sk+j} = -2\delta\hat\beta_{sk+j} + 2\hat{u}^T_{sk+j}\delta u_{sk+j}$. Now, using~\eqref{eq:ubound} and~\eqref{eq:betahatbound},
\begin{align}
\vert \delta\tilde\beta_{sk+j} \vert  &\leq  2\vert \delta\hat\beta_{sk+j}\vert + 2 \| \hat{u}_{sk+j}^T \|_2 \| \delta{u}_{sk+j}\|_2 \nonumber\\
&\color{black}\leq  2  \eps \big(2(sk\pl j) \mn 1 \big) \Big( \big( (N\pl 2s\pl 5)\theta + (4s \pl 9)\bar{\tau}_k \big) \sigma + (34s\pl 55)\bar{u}_{sk+j} \Big) \bar{\Gamma}_k  \bar{u}_{sk+j} \nonumber \\
&\color{black}\phantom{\leq} + 2\eps \Big( \big( (N\pl 2s \pl 5) \theta + (4s\pl 9) \bar{\tau}_k \big)\sigma+ (4s\pl 6)\bar{u}_{sk+j}  \Big) \bar{\Gamma}_k \bar{u}_{sk+j} \nonumber \\
&\color{black}\leq  4\eps (sk\pl j) \big( (N\pl 2s \pl 5) \theta+ (4s\pl 9) \bar{\tau}_k \big)\sigma \bar{\Gamma}_k \bar{u}_{sk+j} \nonumber \\
&\color{black}\phantom{\leq} + 2\eps \Big( \big( 2(sk\pl j)\mn 1 \big) \big( 34s\pl 55 \big) + (4s \pl 6) \Big) \bar{\Gamma}_k  \bar{u}^2_{sk+j}. \label{eq:betatilde}
\end{align}

Using that $\hat\beta_{sk+j}^2\geq 0$ and~\eqref{eq:uhat2}, we can write
\begin{equation}
\Vert \hat{u}_{sk+j} \Vert_2^2 \leq \Vert \hat{u}_{sk+j} \Vert_2^2 \pl \hat{\beta}_{sk+j}^2\leq \sigma^2 \Vert \hat{v}_{sk+j} \Vert_2^2 \pl \hat\beta_{sk+j}^2 \big( \Vert \hat{v}_{sk+j-1} \Vert_2^2 \mn 1 \big) \pl  | \delta\tilde\beta_{sk+j} |. \label{eq:usqr}
\end{equation}
Let $\mu \equiv \max{\big\{\bar{u}_{sk+j}, \|A\|_2\big\}}$. Then~\eqref{eq:usqr} along with bounds in~\eqref{mainthm:norm},~\eqref{eq:betabound}, and~\eqref{eq:betatilde} gives
\color{black}
\begin{equation}
\Vert \hat{u}_{sk+j} \Vert_2^2 \leq  \sigma^2 +  4\eps(sk\pl j)\big( (N\pl 2s\pl 5)\theta + (4s\pl 9)\bar{\tau}_k  + (34s\pl 55) \big) \bar{\Gamma}_k  \mu^2. 
%
\label{eq:uboundmu}
\end{equation}
\color{black}
We consider the two possible cases for $\mu$. First, if $\mu=\|A\|_2$, then
\color{black} 
\begin{align*}
\Vert \hat{u}_{sk+j} \Vert_2^2 &\leq \sigma^2 \pl  4\eps(sk\pl j)\big( (N\pl 2s\pl 5)\theta + (4s\pl 9)\bar{\tau}_k  + (34s\pl 55) \big) \bar{\Gamma}_k  \sigma^2 \\
&\leq \sigma^2 \bigg(1+ 4\eps(sk\pl j)\big( (N\pl 2s\pl 5)\theta + (4s\pl 9)\bar{\tau}_k  + (34s\pl 55) \big) \bar{\Gamma}_k   \bigg).
\end{align*}
\color{black}
If instead $\mu=\bar{u}_{sk+j}$, since the bound in~\eqref{eq:uboundmu} holds for all $\| \hat{u}_{sk+j} \|_2^2$,  it also holds for $\bar{u}_{sk+j}^2=\mu^2$, and thus ignoring terms of order $\eps^2$,
\color{black}
\begin{align*}
\mu^2 &\leq  \sigma^2 \pl  4\eps(sk\pl j)\big( (N\pl 2s\pl 5)\theta + (4s\pl 9)\bar{\tau}_k  + (34s\pl 55) \big) \bar{\Gamma}_k  \mu^2\\
&\leq   \sigma^2 \pl  4\eps(sk\pl j)\big( (N\pl 2s\pl 5)\theta + (4s\pl 9)\bar{\tau}_k  + (34s\pl 55) \big) \bar{\Gamma}_k  \sigma^2\\
&\leq  \sigma^2 \bigg(1+ 4\eps(sk\pl j)\big( (N\pl 2s\pl 5)\theta + (4s\pl 9)\bar{\tau}_k  + (34s\pl 55) \big) \bar{\Gamma}_k   \bigg).
\end{align*}
\color{black}
Substituting this into~\eqref{eq:uboundmu}, we have
\color{black}
\begin{equation}
\| \hat{u}_{sk+j} \|_2^2\leq \sigma^2 \bigg(1+ 4\eps(sk\pl j)\big( (N\pl 2s\pl 5)\theta + (4s\pl 9)\bar{\tau}_k  + (34s\pl 55) \big) \bar{\Gamma}_k   \bigg).
\label{eq:mu2}
\end{equation}
\color{black}
Thus in both cases we obtain the same bound on $\|\hat{u}_{sk+j}\|_2^2$  so~\eqref{eq:mu2} holds. 

Taking the square root of~\eqref{eq:mu2}, we have
\color{black}
\begin{equation}
\|\hat{u}_{sk+j} \|_2 \leq \sigma \bigg(1+ 2\eps(sk\pl j)\big( (N\pl 2s\pl 5)\theta + (4s\pl 9)\bar{\tau}_k  + (34s\pl 55) \big) \bar{\Gamma}_k   \bigg),
\label{eq:yayu}
\end{equation}
\color{black}
and substituting~\eqref{eq:yayu} into~\eqref{eq:orthbound1},~\eqref{eq:deltavnorm}, and~\eqref{eq:etasbounds}, we obtain~\eqref{mainthm:orth},~\eqref{mainthm:deltav}, and~\eqref{mainthm:etas} in Theorem~\ref{thm:mainthm}, respectively, under the constraint 
\color{black}
\begin{equation*}
2\eps(sk\pl j)\big( (N\pl 2s\pl 5)\theta + (4s\pl 9)\bar{\tau}_k  + (34s\pl 55) \big) \bar{\Gamma}_k \ll 1.
\end{equation*}
We note that this constraint is obtained from what are almost surely very loose upper bounds, and thus this requirement is overly strict in practice. We also note that in constrast to the uniform precision case, there is again only a factor of $\bar{\Gamma}_k$ instead of $\bar{\Gamma}_k^2$. 
\color{black}

We now only need to prove~\eqref{mainthm:cols}. We first multiply both sides of~\eqref{eq:betav} by their own transposes to obtain
\begin{align*}
\hat{\beta}_{sk+j+1}^2 \| \hat{v}_{sk+j+1} \|_2^2 &= \| \hat{\mathcal{Y}}_k \hat{w}'_{k,j} \|_2^2 + 2 \delta w_{sk+j}^T \hat{\mathcal{Y}}_k \hat{w}'_{k,j}.
\end{align*}
Adding $\hat{\alpha}_{sk+j}^2 -\| \hat{u}_{sk+j} \|_2^2$ to both sides, 
\begin{align*}
\hat{\beta}_{sk+j+1}^2 \| \hat{v}_{sk+j+1} \|_2^2 + \hat{\alpha}_{sk+j}^2 -\| \hat{u}_{sk+j} \|_2^2 &= \| \hat{\mathcal{Y}}_k \hat{w}'_{k,j} \|_2^2 + \hat{\alpha}_{sk+j}^2 -\| \hat{u}_{sk+j} \|_2^2 \\
&\phantom{=}+ 2 \delta w_{sk+j}^T \hat{\mathcal{Y}}_k \hat{w}'_{k,j}.
\end{align*}
Substituting in~\eqref{eq:uhat2} on the left-hand side above,
\begin{multline*}
\hat{\beta}_{sk+j+1}^2 \| \hat{v}_{sk+j+1} \|_2^2 + \hat{\alpha}_{sk+j}^2 -\| A \hat{v}_{sk+j} \|_2^2 - \hat{\beta}_{sk+j}^2 (\| \hat{v}_{sk+j-1}\|_2^2 -2) - \delta \tilde{\beta}_{sk+j} =\\ \| \hat{\mathcal{Y}}_k \hat{w}'_{k,j} \|_2^2 + \hat{\alpha}_{sk+j}^2 -\| \hat{u}_{sk+j} \|_2^2 + 2 \delta w_{sk+j}^T \hat{\mathcal{Y}}_k \hat{w}'_{k,j},
\end{multline*}
and then subtracting $\hat{\beta}_{sk+j+1}^2$ from both sides gives
\begin{multline*}
\hat{\beta}_{sk+j+1}^2 (\| \hat{v}_{sk+j+1} \|_2^2 - 1)+ \hat{\alpha}_{sk+j}^2 -\| A \hat{v}_{sk+j} \|_2^2 - \hat{\beta}_{sk+j}^2 (\| \hat{v}_{sk+j-1}\|_2^2 -2) - \delta \tilde{\beta}_{sk+j} =\\ \| \hat{\mathcal{Y}}_k \hat{w}'_{k,j} \|_2^2 + \hat{\alpha}_{sk+j}^2 -\| \hat{u}_{sk+j} \|_2^2 + 2 \delta w_{sk+j}^T \hat{\mathcal{Y}}_k \hat{w}'_{k,j} - \hat{\beta}_{sk+j+1}^2.
\end{multline*}
This can be rearranged to give
\begin{align*}
\hat{\beta}_{sk+j+1}^2  + \hat{\alpha}_{sk+j}^2 + \hat{\beta}_{sk+j}^2 -\| A \hat{v}_{sk+j} \|_2^2 &= \| \hat{\mathcal{Y}}_k \hat{w}'_{k,j} \|_2^2 + \hat{\alpha}_{sk+j}^2 -\| \hat{u}_{sk+j} \|_2^2 \\
&\phantom{=} + 2 \delta w_{sk+j}^T \hat{\mathcal{Y}}_k \hat{w}'_{k,j}+\hat{\beta}_{sk+j}^2 (\| \hat{v}_{sk+j-1}\|_2^2 -1)\\
&\phantom{=} -\hat{\beta}_{sk+j+1}^2 (\| \hat{v}_{sk+j+1} \|_2^2 -1) + \delta \tilde{\beta}_{sk+j} ,
\end{align*}
and finally, using~\eqref{mainthm:norm},~\eqref{eq:wtf},~\eqref{eq:betabound},~\eqref{eq:deltwvec}, and~\eqref{eq:betatilde} gives the bound
\begin{multline*}
\Big\vert \hat\beta_{sk+j+1}^2  + \hat\alpha_{sk+j}^2 + \hat{\beta}_{sk+j}^2 - \Vert A \hat{v}_{sk+j} \Vert_2^2 \Big\vert \leq \\
\color{black} 
2\eps (sk\pl j) \Big(  2((N\pl 2s\pl 5)\theta \pl (4s\pl 9)\bar{\tau}_k  + (10s\pl 16)) \pl 3\cdot 2(9s\pl 14)  \Big) \bar{\Gamma}_k\sigma^2.
\end{multline*}
\color{black}
This proves~\eqref{mainthm:cols} and thus completes the proof of Theorem~\ref{thm:mainthm}.

\section{Numerical experiments}
\label{sec:num}

We perform a few numerical experiments in MATLAB (R2020a) in order to demonstrate the validity of the bounds and enable comparison between the uniform finite precision $s$-step Lanczos algorithm in double precision and the variant that uses a mixed precision approach, with quad precision used for the Gram matrix computations as described and double precision otherwise. For double precision, we use the built-in MATLAB datatype and for quadruple precision we use the Advanpix Toolbox \cite{ad06}. In all tests we use starting vector of unit 2-norm with equal components. For each test problem, we show the both the measured/computed values (solid lines) and the bounds on the values (dashed lines) given in \eqref{mainthm:deltav}, \eqref{mainthm:orth}, \eqref{mainthm:norm}, and \eqref{mainthm:cols} for the mixed precision case (red) and the equivalent bounds from \cite[Theorem 4.2]{carson2015accuracy} for the uniform precision case (blue). We additionally plot the quantity $\bar\Gamma_k$ in black. 

\begin{figure}
    \centering
    \includegraphics[trim={4cm 8cm 4cm 8cm},clip,width=6cm]{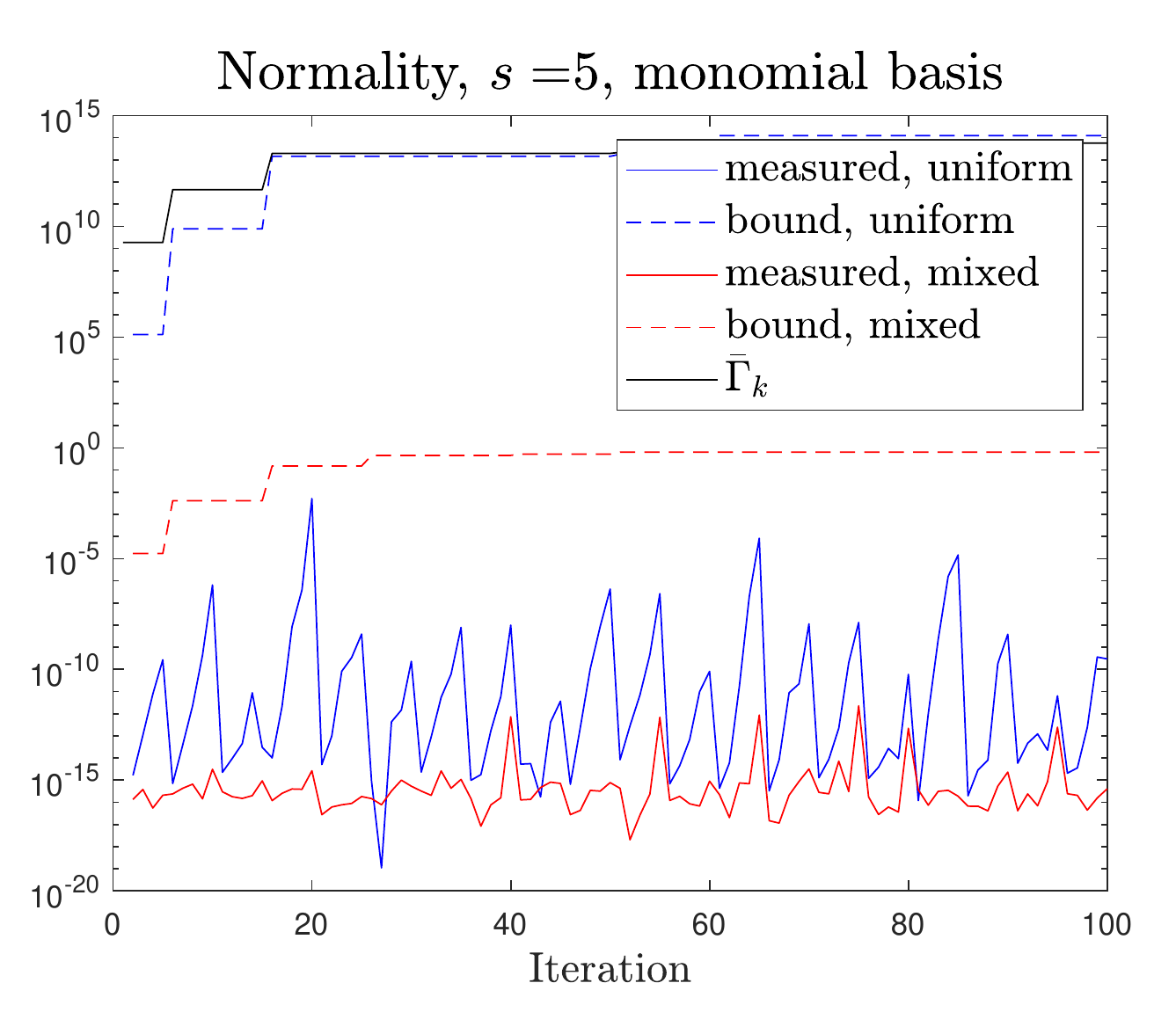}
    \includegraphics[trim={4cm 8cm 4cm 8cm},clip,width=6cm]{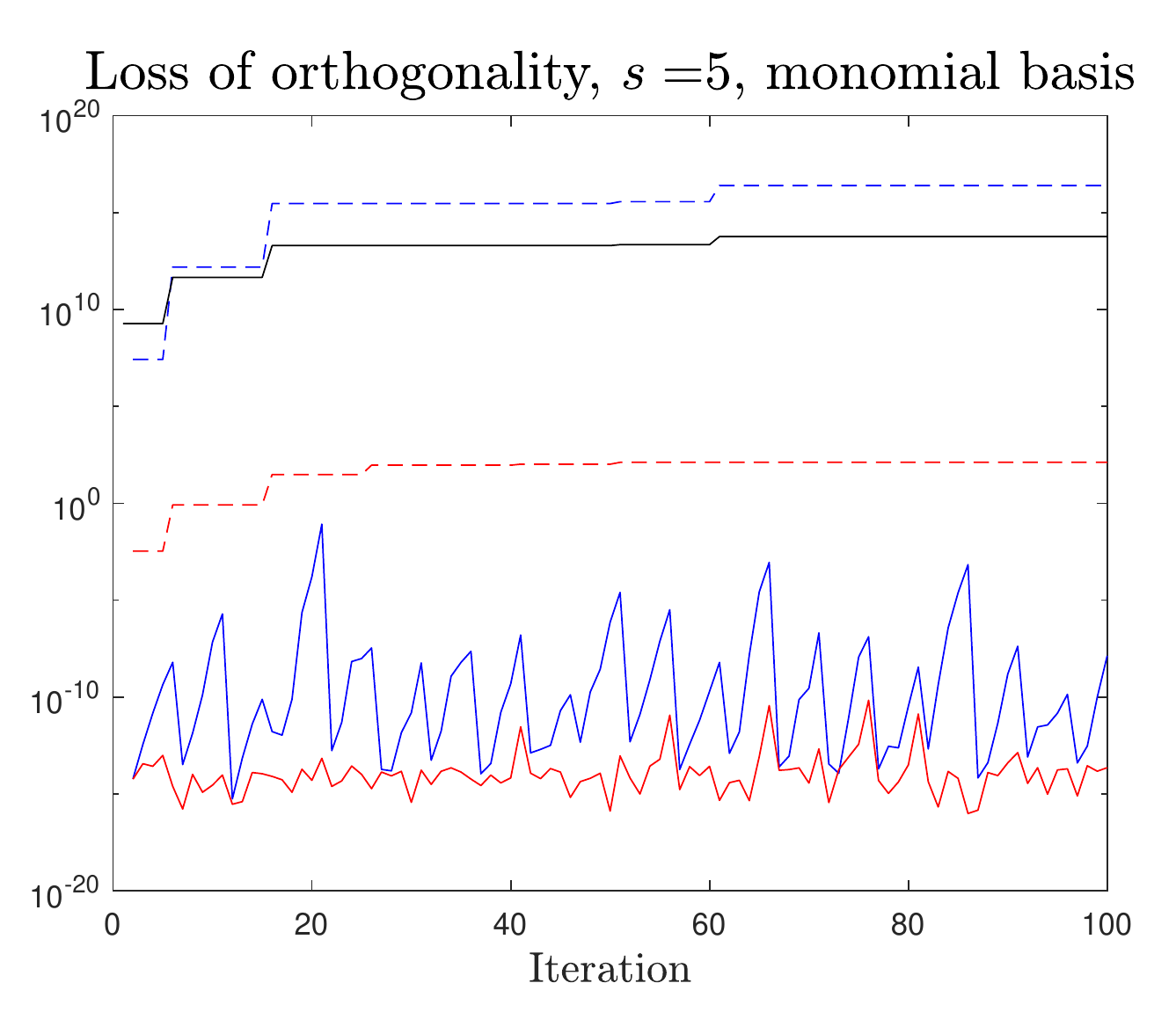}\\
    \includegraphics[trim={4cm 8cm 4cm 8cm},clip,width=6cm]{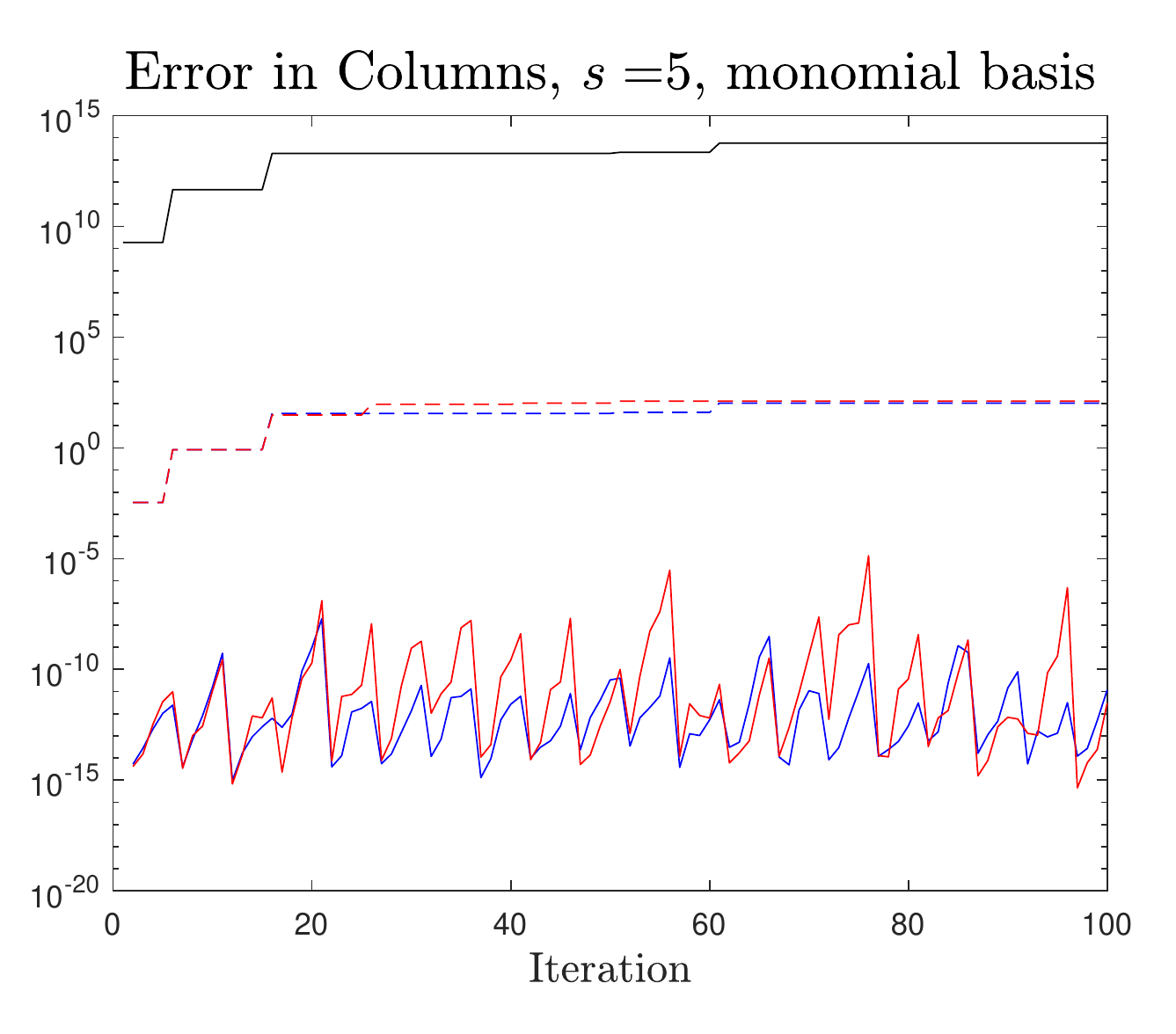}
    \includegraphics[trim={4cm 8cm 4cm 8cm},clip,width=6cm]{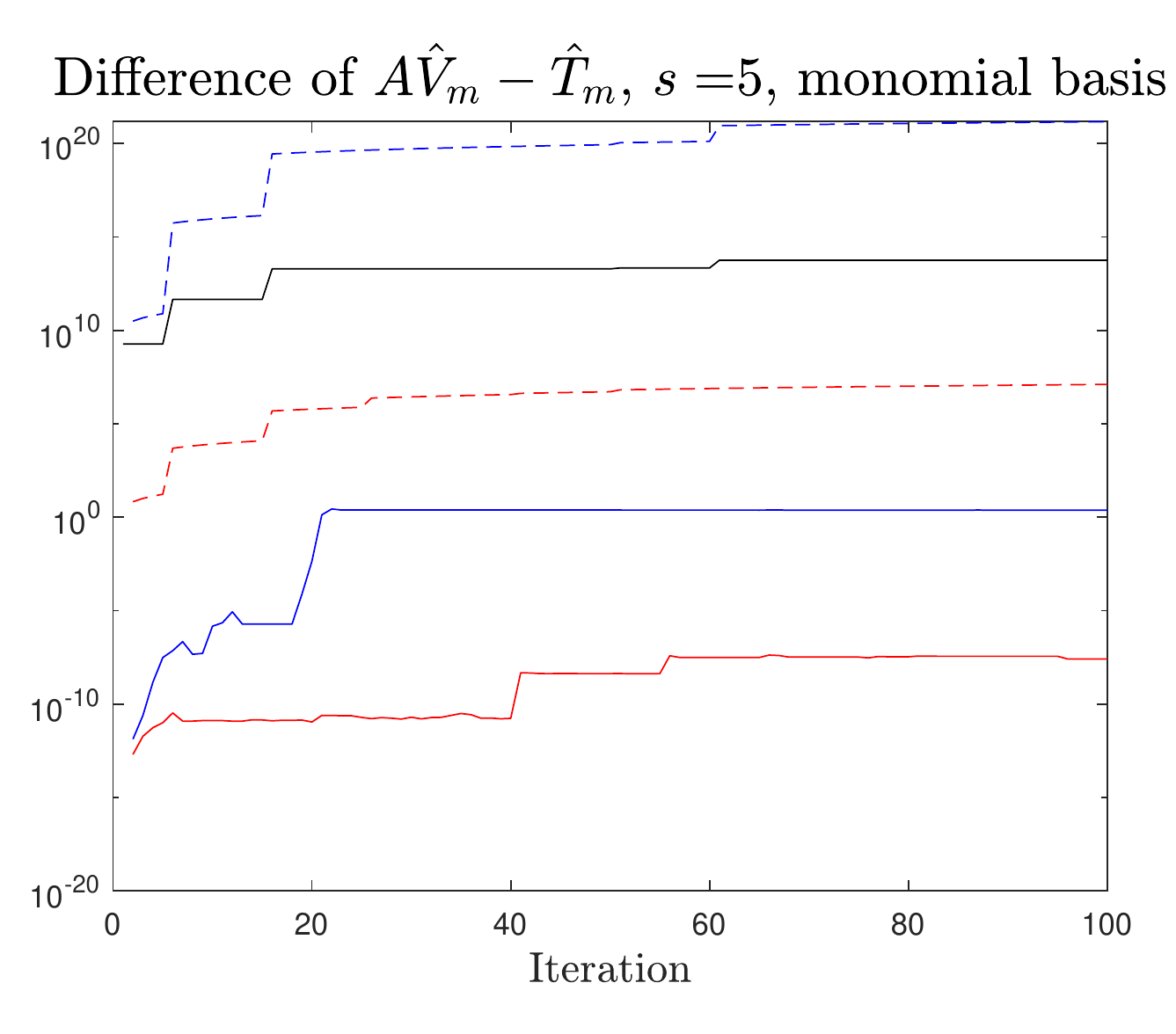}
    \caption{Comparison of uniform and mixed precision $s$-step Lanczos for the diagonal matrix problem \eqref{eq:strakos} with $n=100$, $\lambda_1=10^{2}$, $\lambda_n=10^{-3}$, $\rho=0.65$, using $s=5$ and a monomial basis.}
    \label{fig:strakoss5}
\end{figure}

We first test a problem commonly used in studies of the behavior of CG, the diagonal matrix with entries 
\begin{equation}
\lambda_i = \lambda_1 + \left(\frac{i-1}{n-1} \right)(\lambda_n - \lambda_1)\rho^{n-i}, \quad i = 1,\ldots,n,
\label{eq:strakos}
\end{equation}
with parameters $n=100$, $\lambda_1=10^{-3}$, $\lambda_n=10^2$, and $\rho=0.65$. The bounds and computed quantities for both uniform precision and mixed precision $s$-step Lanczos using $s=5$ with a monomial basis are shown in Figure \ref{fig:strakoss5}. 

A few things are clear from these plots. First, the loss of orthogonality, the loss of normality, and the bound~\eqref{mainthm:cols} are much improved by the use of mixed precision. It is also clear that the bounds for both the uniform precision and mixed precision approaches can be large overestimates of the actual quantities. We have tried to stress this in the analysis. The bounds themselves should not be taken as estimates of the errors; rather it is the structure of the accumulation and amplification of errors that tells us something meaningful about the numerical behavior. In this case, it is the fact that using mixed precision, we have linear rather than quadratic dependence on the quantity $\bar\Gamma_k$, and this is reflected in the computed quantities. We note that the computed quantities in \eqref{mainthm:deltav} (the lower left plots) are not necessarily smaller for the mixed precision approach than for the uniform precision approach. We expect this behavior, since the use of mixed precision does not improve the bound \eqref{mainthm:deltav}, which depends only on $\bar{\Gamma}_k$ even in the uniform precision case.  

\begin{figure}
    \centering
    \includegraphics[trim={4cm 8cm 4cm 8cm},clip,width=6cm]{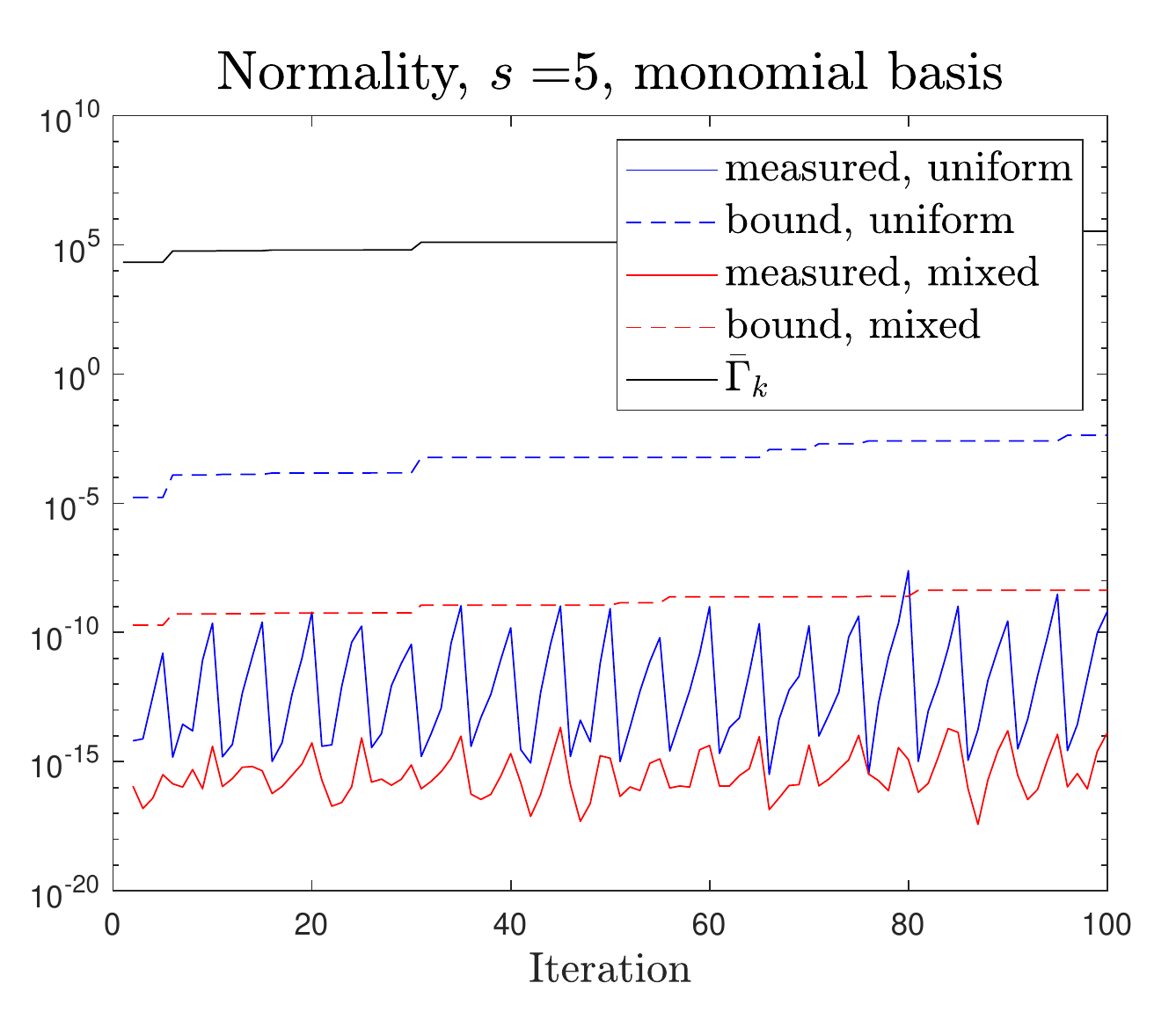}
    \includegraphics[trim={4cm 8cm 4cm 8cm},clip,width=6cm]{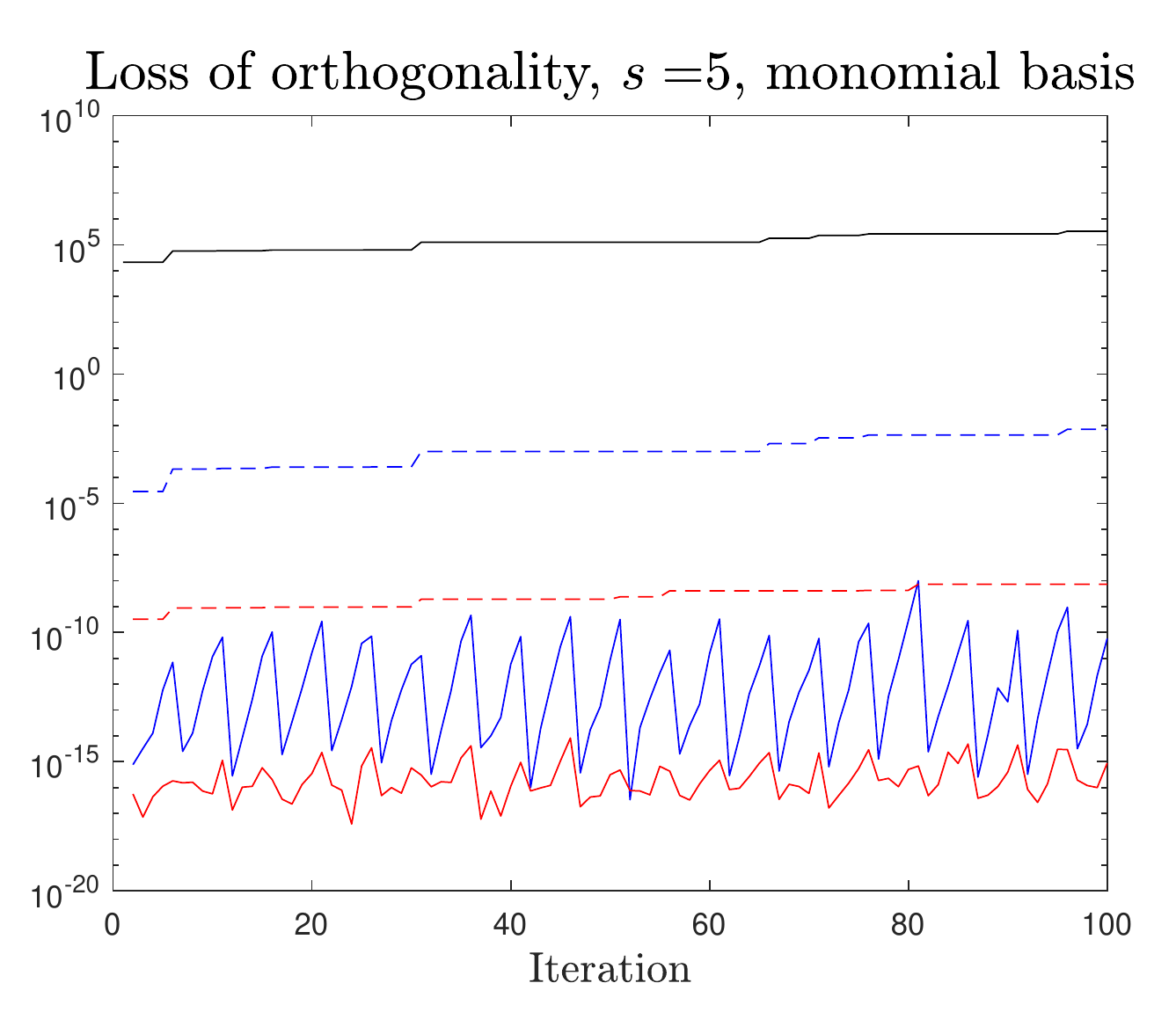}\\
    \includegraphics[trim={4cm 8cm 4cm 8cm},clip,width=6cm]{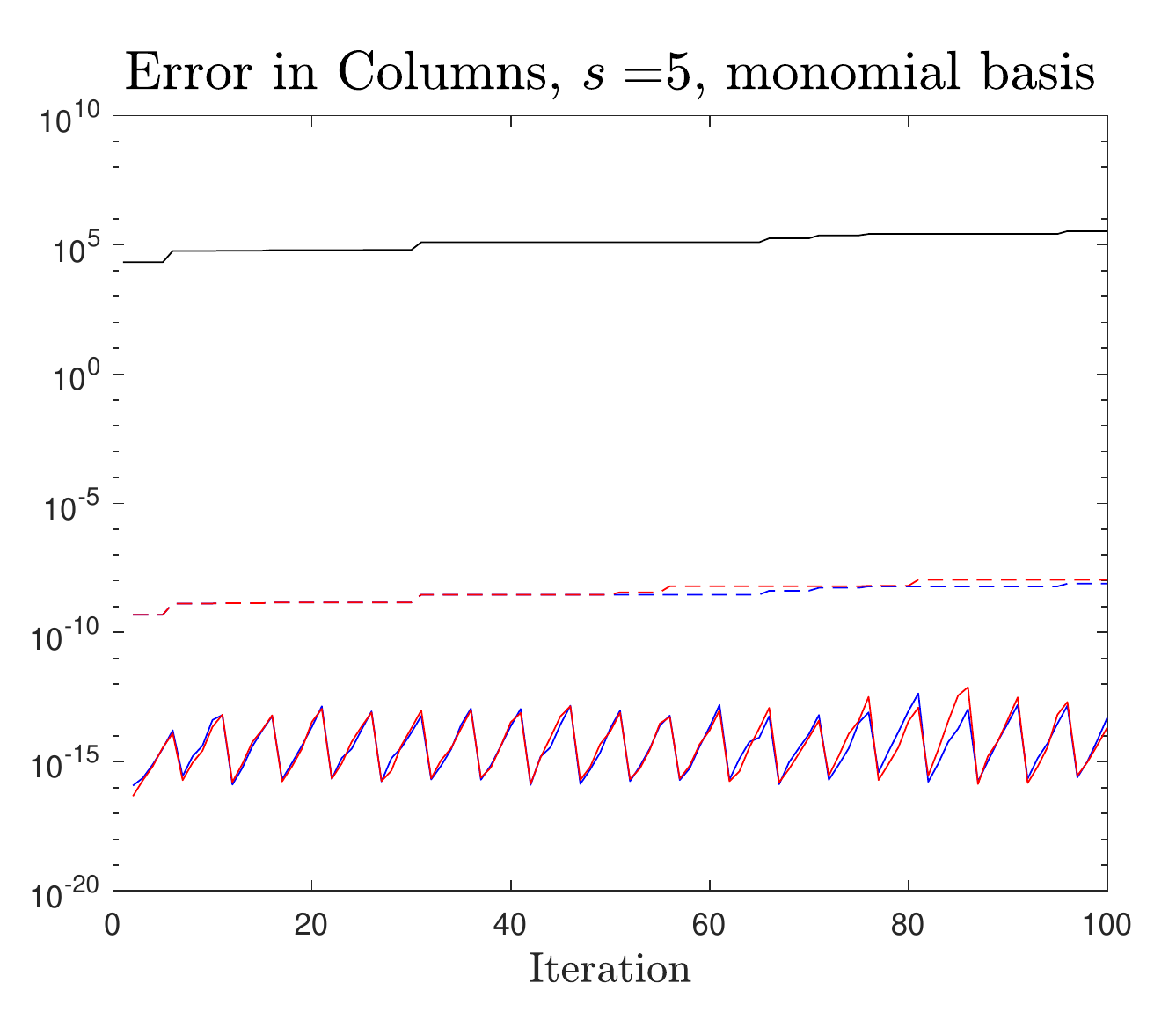}
    \includegraphics[trim={4cm 8cm 4cm 8cm},clip,width=6cm]{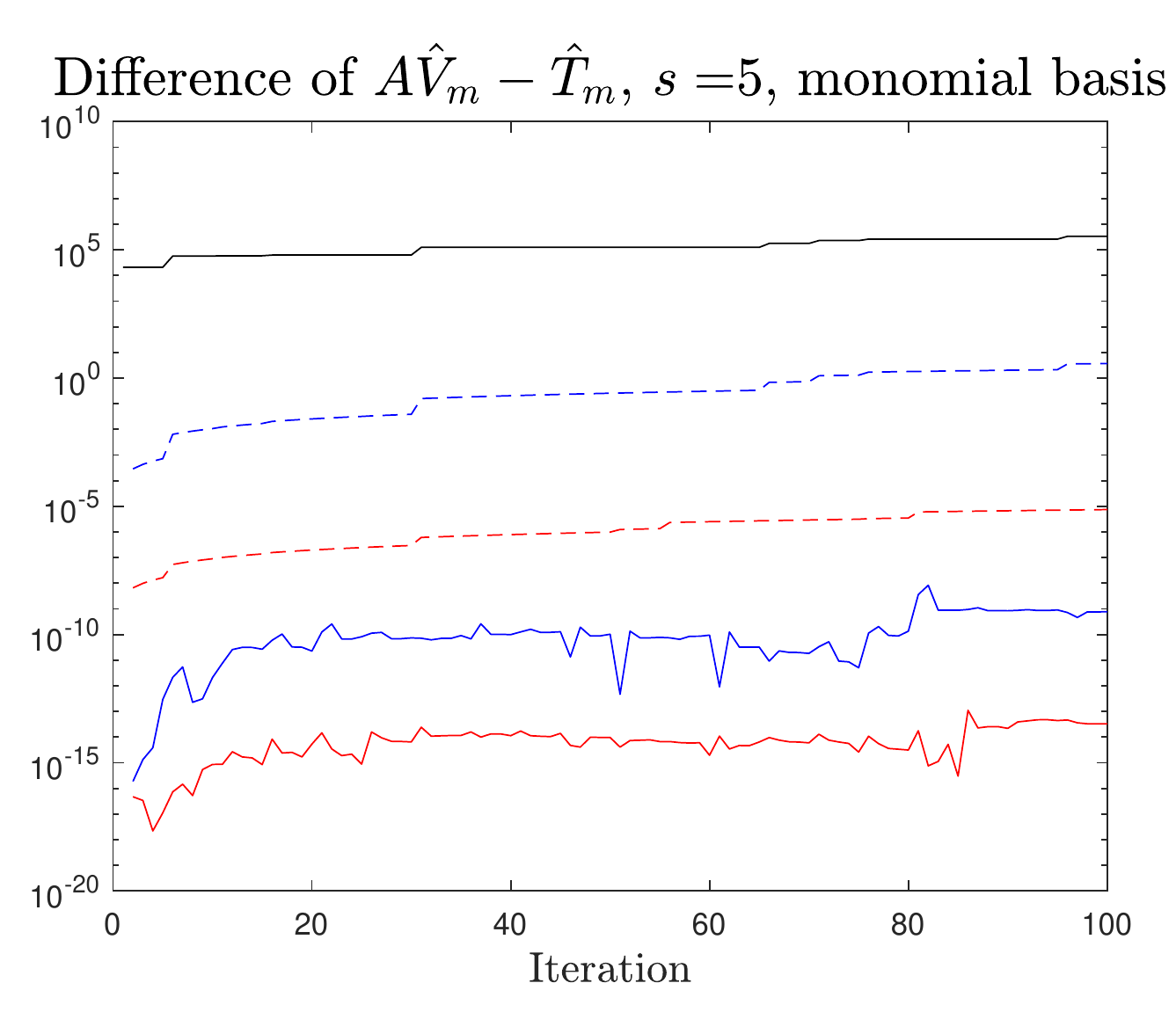}
    \caption{Comparison of uniform and mixed precision $s$-step Lanczos for the \texttt{nos4} problem from SuiteSparse, with $s=5$ and a monomial basis.}
    \label{fig:nos4s5}
\end{figure}

\begin{figure}
    \centering
    \includegraphics[trim={4cm 8cm 4cm 8cm},clip,width=6cm]{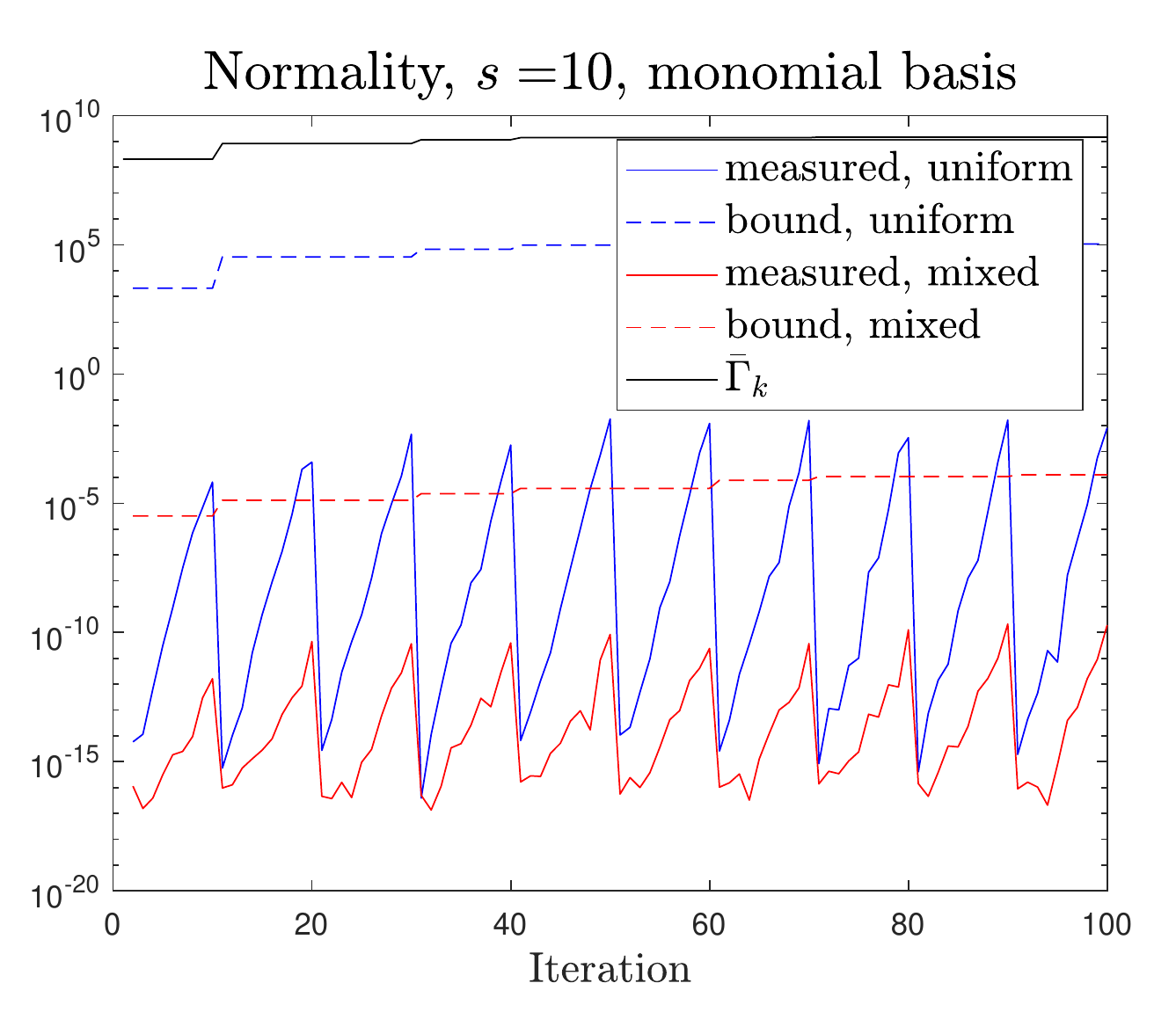}
    \includegraphics[trim={4cm 8cm 4cm 8cm},clip,width=6cm]{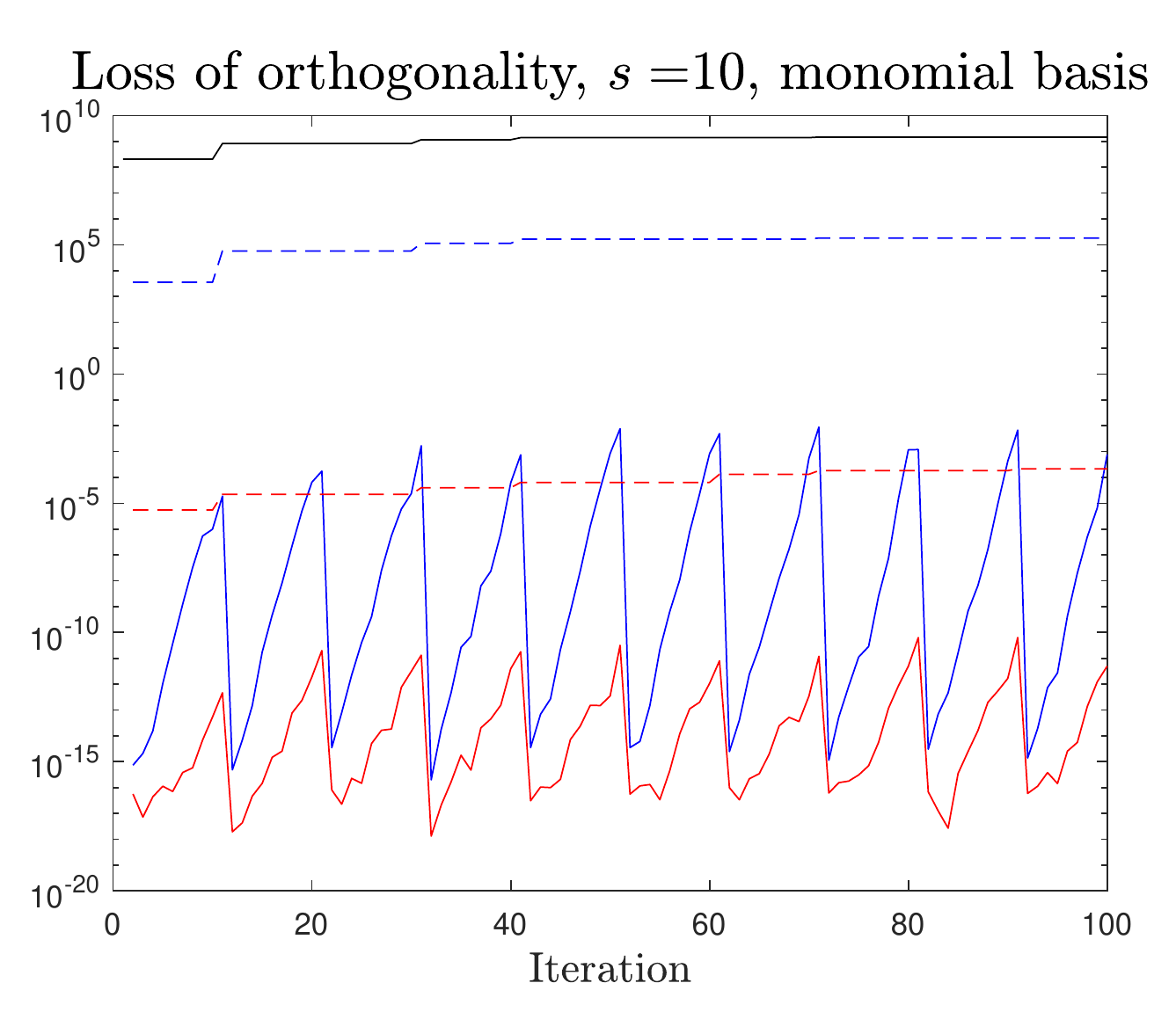}\\
    \includegraphics[trim={4cm 8cm 4cm 8cm},clip,width=6cm]{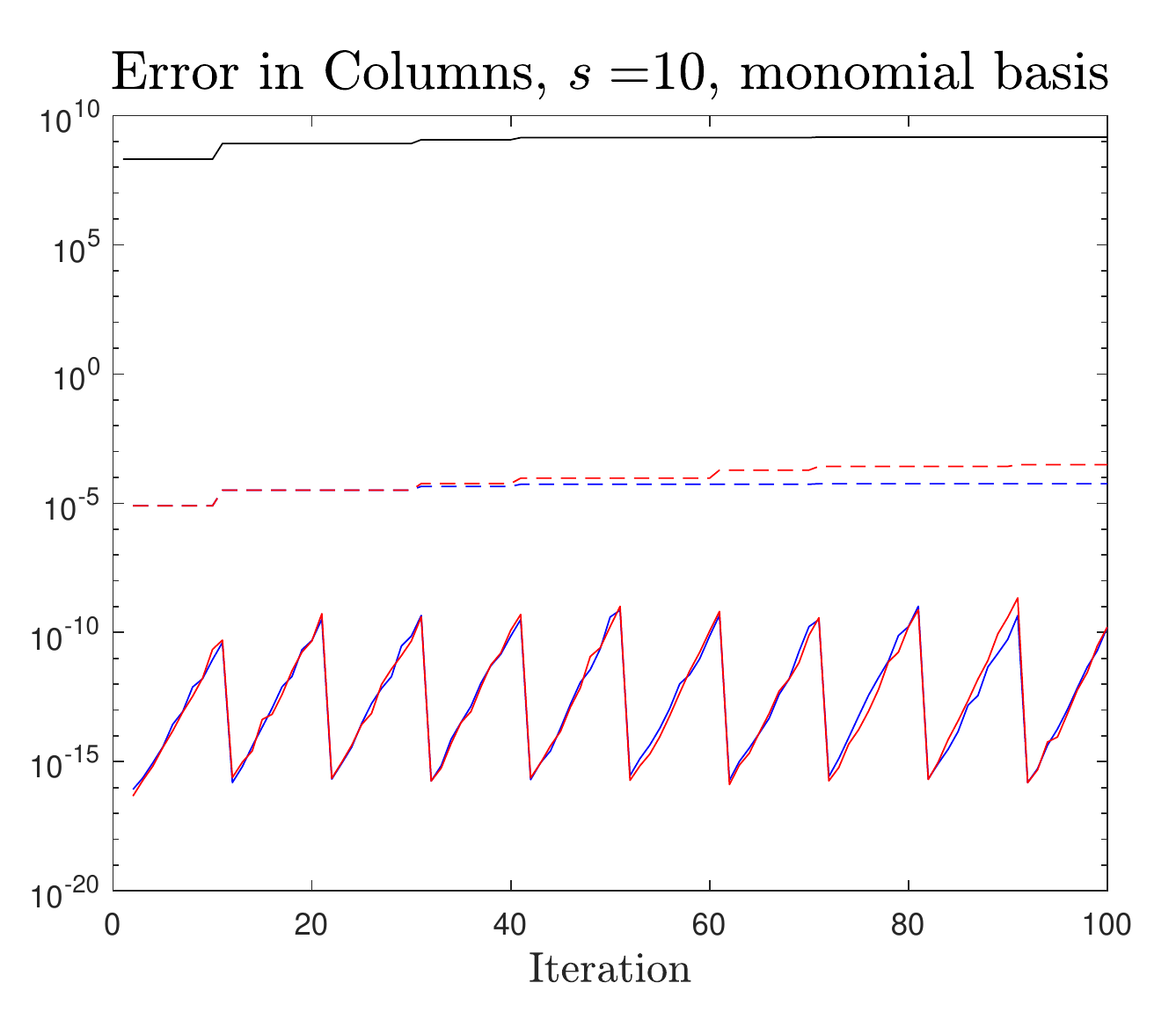}
    \includegraphics[trim={4cm 8cm 4cm 8cm},clip,width=6cm]{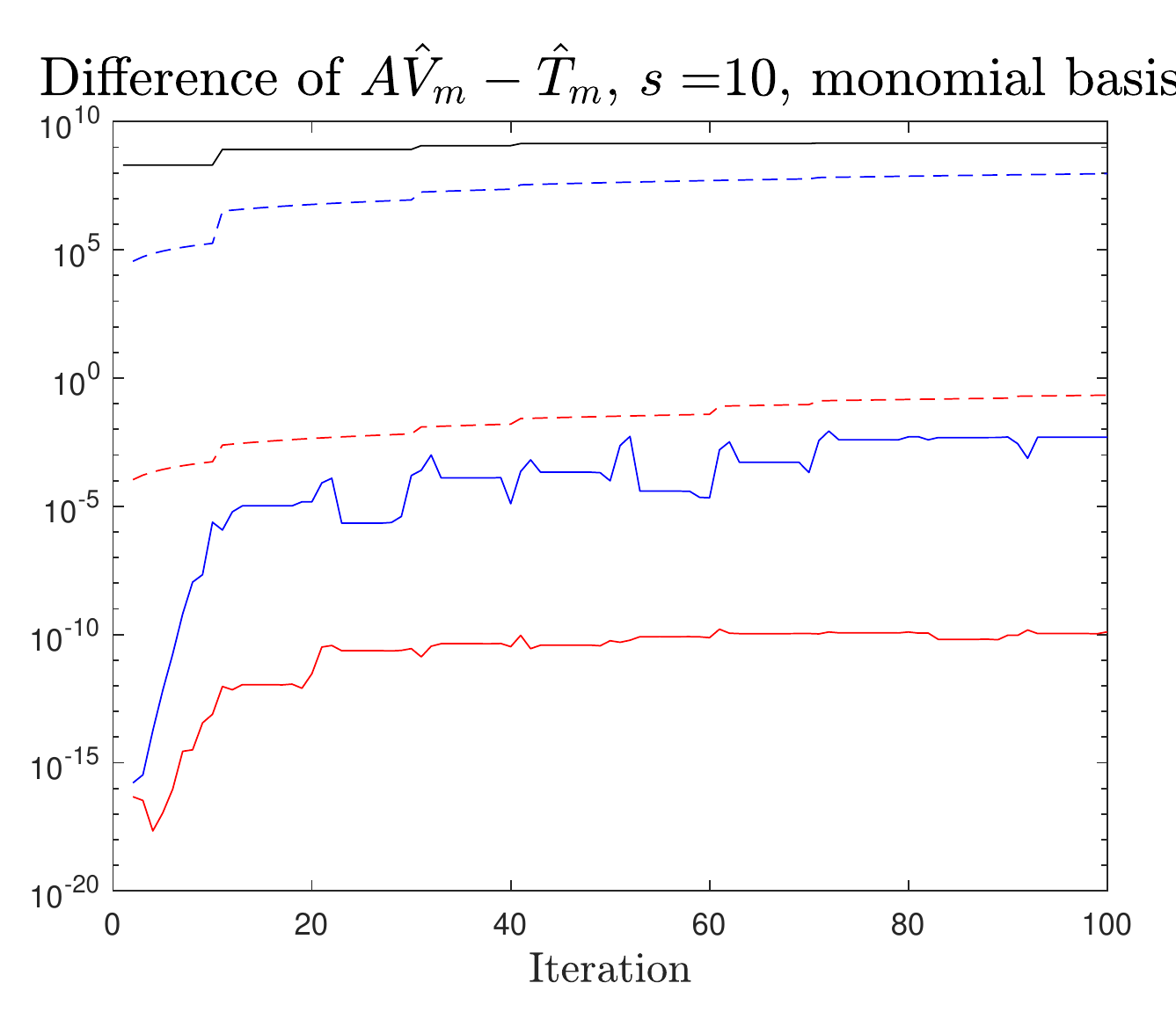}
    \caption{Comparison of uniform and mixed precision $s$-step Lanczos for the \texttt{nos4} problem from SuiteSparse, with $s=10$ and a monomial basis.}
    \label{fig:nos4s10}
\end{figure}

\begin{figure}
    \centering
    \includegraphics[trim={4cm 8cm 4cm 8cm},clip,width=6cm]{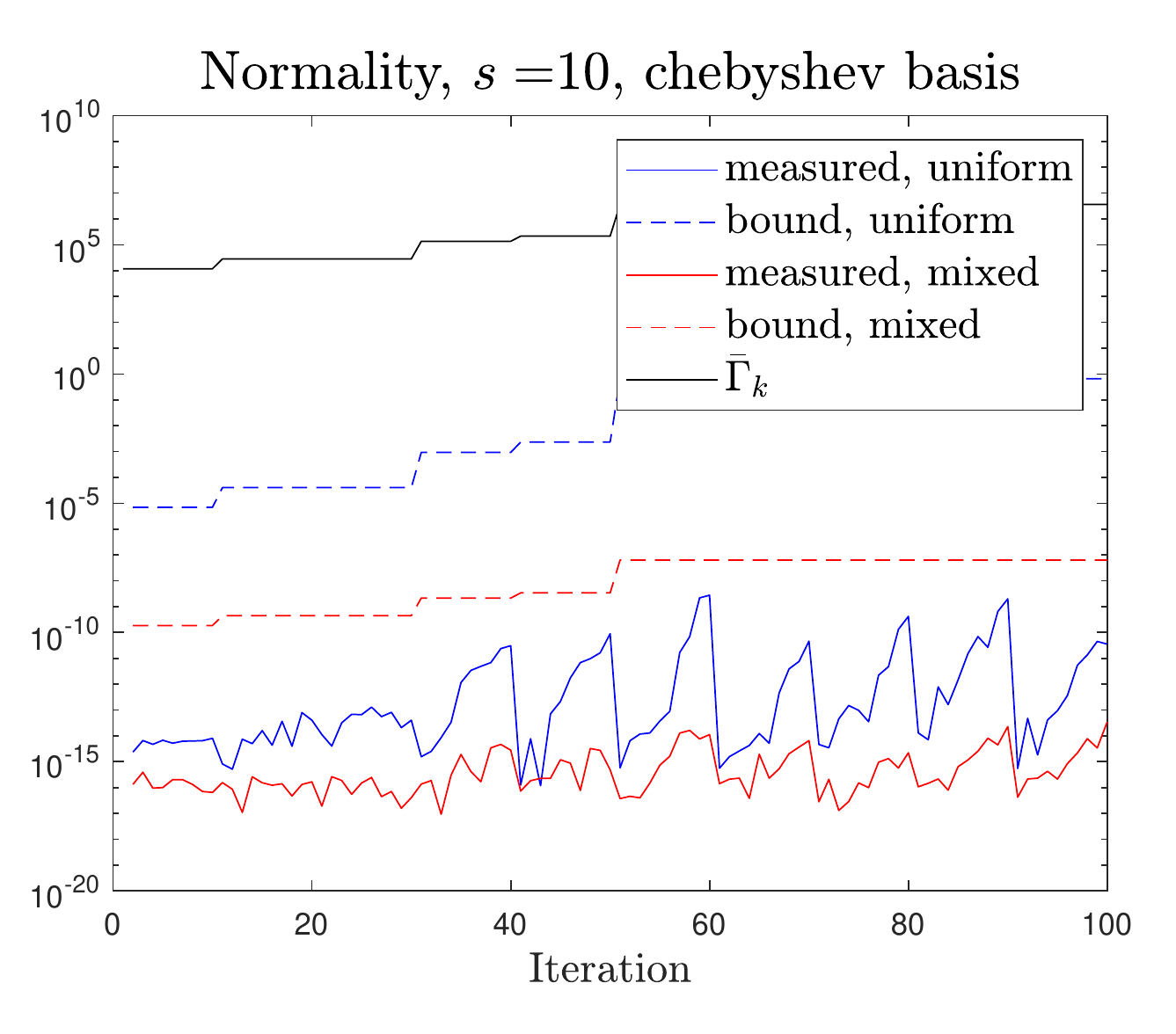}
    \includegraphics[trim={4cm 8cm 4cm 8cm},clip,width=6cm]{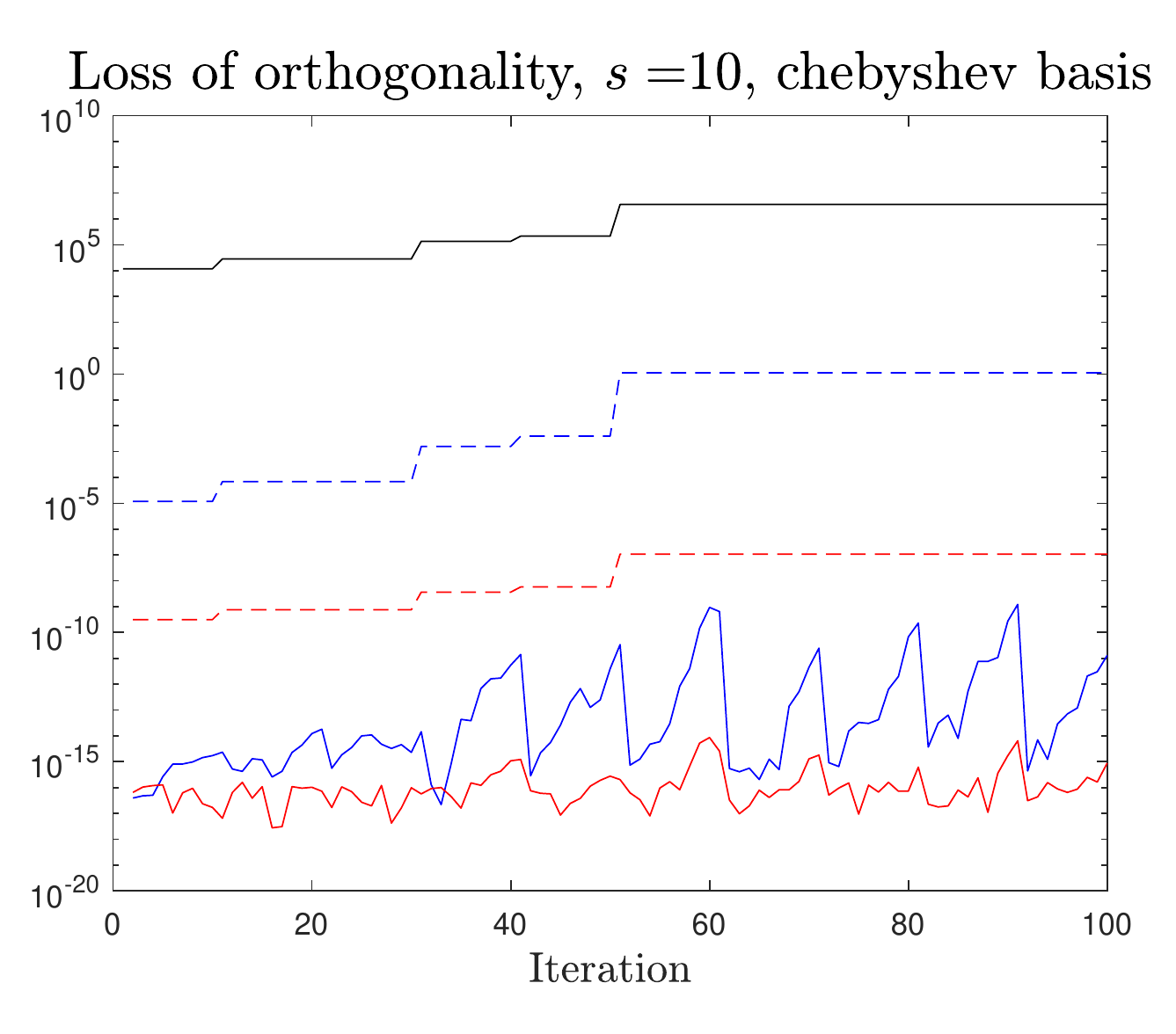}\\
    \includegraphics[trim={4cm 8cm 4cm 8cm},clip,width=6cm]{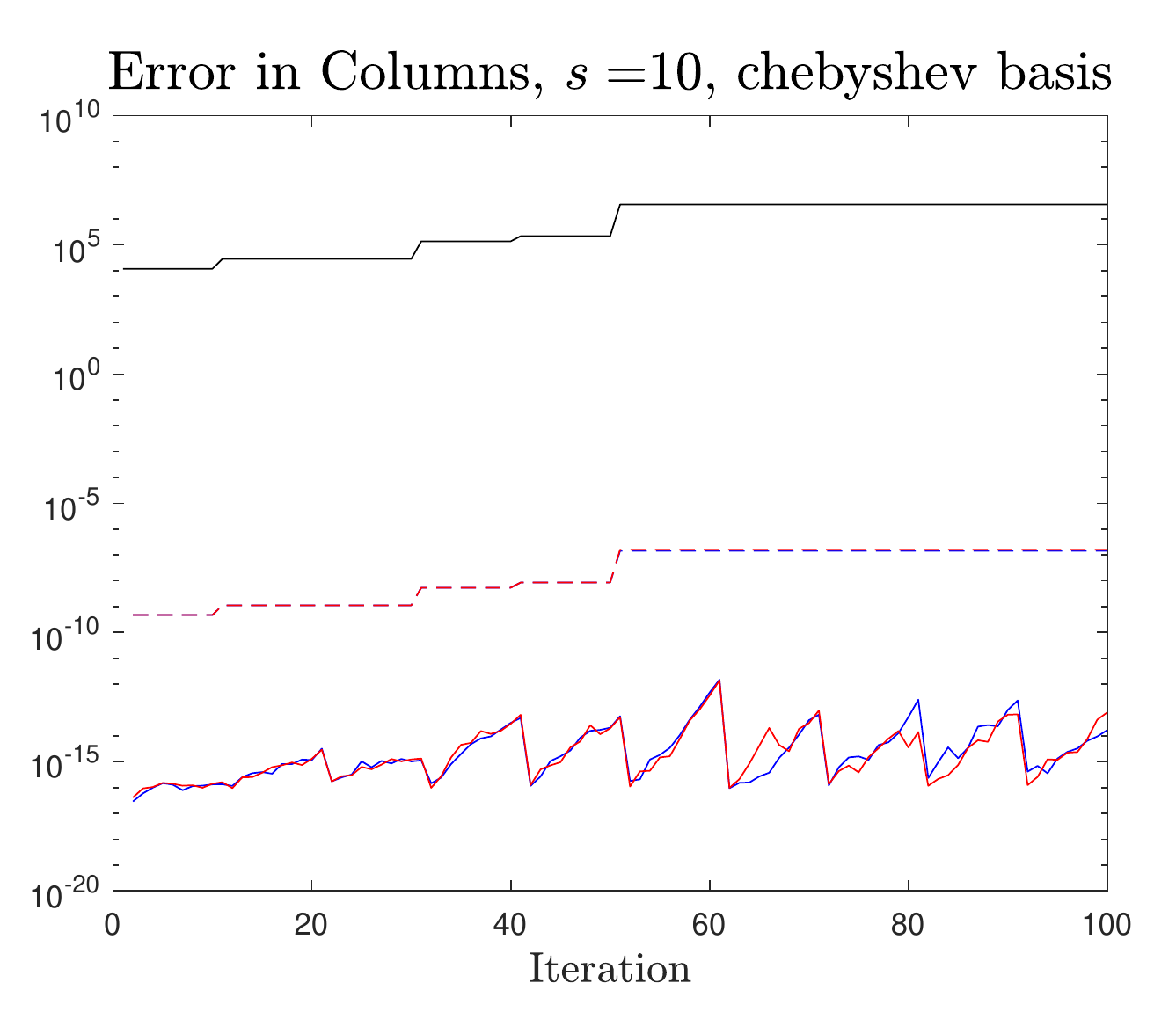}
    \includegraphics[trim={4cm 8cm 4cm 8cm},clip,width=6cm]{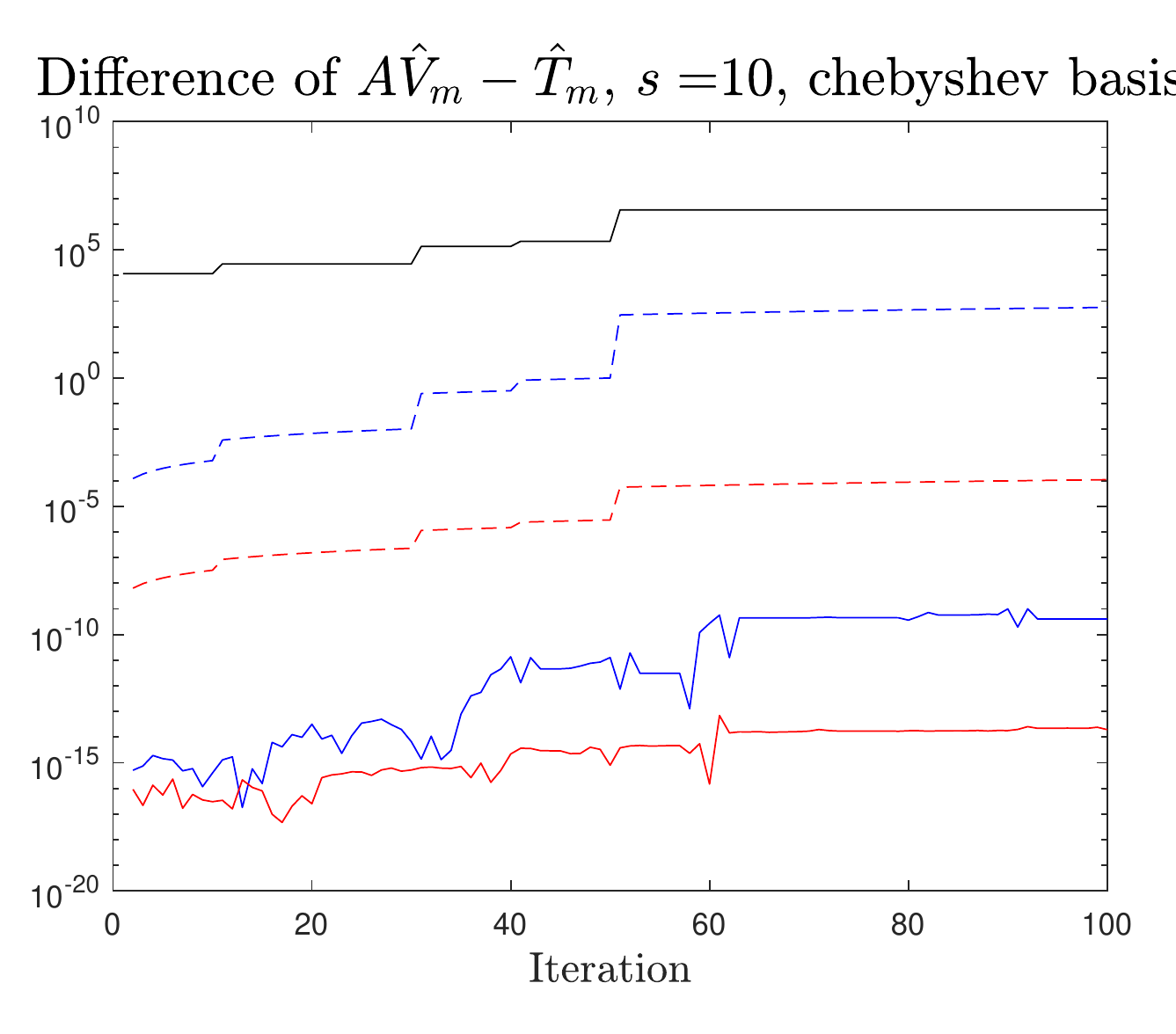}
    \caption{Comparison of uniform and mixed precision $s$-step Lanczos for the \texttt{nos4} problem from SuiteSparse, with $s=10$ and a Chebyshev basis.}
    \label{fig:nos4s10c}
\end{figure}

We next perform tests on the matrix \texttt{nos4} from the SuiteSparse collection \cite{davis2011university}. This matrix has dimension $n=100$ with condition number around $10^3$. Here we use $s=5$ with a monomial basis in Figure \ref{fig:nos4s5}, $s=10$ with a monomial basis in Figure \ref{fig:nos4s10}, and $s=10$ with a Chebyshev basis in Figure \ref{fig:nos4s10c}. Similar conclusions can be drawn here. It is clear that there is benefit to the mixed precision approach, especially as the quantity $\bar\Gamma_k$ grows large. In the $s=10$ case with the monomial basis, orthogonality and normality are nearly lost for the uniform precision approach, whereas the mixed precision approach dampens this effect significantly. As can be seen in Figure \ref{fig:nos4s10c}, using a more well-conditioned basis such as the Chebyshev basis results in lower $\bar\Gamma_k$ and thus improved behavior in both cases, although the benefit of mixed precision is still present. 

\section{Accuracy and Convergence of Eigenvalues and Eigenvectors}
\label{sec:evals}

In previous work \cite{carson2015accuracy}, we have demonstrated that bounds on accuracy and convergence for the finite precision classical Lanczos algorithm given by Paige \cite{paige1980accuracy} also apply to the uniform precision $s$-step Lanczos case assuming a bound on the maximum condition number of the precomputed $s$-step Krylov bases. The constraint in this case is that 
\[
\bar{\Gamma}_k < \big(24\eps (n+11s+15) \big)^{-1/2} = O(1/\sqrt{n\eps}).
\]
This comes from Paige's assumption that $\eps_0 < 1/12$ (see \cite[Eqn. (2.16)]{paige1980accuracy} and the value $\eps_0 = 2\eps(n+11s+15)\bar{\Gamma}_k^2$ (see \cite[Eqn. (4.10)]{carson2015accuracy}) for $s$-step Lanczos.

For the mixed precision $s$-step case, as in the uniform precision $s$-step case, we can also say that Paige's results apply as long as a certain constraint is met. In this case, the constraint is much looser; from \eqref{eq:e0e1}, we now require only that 
\[
\bar{\Gamma}_k < \big(2\eps(6s+11) \big)^{-1} = O(1/\eps).
\]
For example, in double precision, we can't expect predictable behavior from the uniform precision $s$-step Lanczos algorithm unless $\bar{\Gamma}_k \lesssim 10^8$, but for the mixed precision case, we can expect Paige's bounds to hold as long as $\bar{\Gamma}_k \lesssim 10^{16}$. 

In practical terms, these results have two interpretations. First, for a fixed value of $s$, we expect the mixed precision $s$-step Lanczos method to produce accuracy and convergence of Ritz values closer to that of the classical algorithm. Second, this means that the mixed precision variant can handle more ill-conditioned bases, corresponding to higher values of $s$ than the uniform precision case (and thus has more opportunity to potentially reduce communication in each iteration).

For completeness, we state the results here, and note that the proofs of the corresponding theorems and the complete analysis can be found in the work of Paige \cite{paige1980accuracy}. We assume, as Paige does, that 
\begin{equation}
\hat{\beta}_{i+1} \neq 0 \hspace{2mm}\text{for}\hspace{2mm}i\in\{1,\ldots,m\}, \quad m(3\eps_0+2\eps_1)\leq 1, \text{ and } \eps_0 < \frac{1}{12}.
\label{eq:assumptions}
\end{equation}
We further introduce the quantity (as does Paige),
\begin{equation}
  \eps_2 \equiv \sqrt{2}\max\{6\eps_0,\eps_1\}.  
  \label{eq:eps2}
\end{equation}
Importantly, note that in the uniform precision $s$-step Lanczos case, $\eps_2$ will contain a term $\bar{\Gamma}_k^2$, and in the mixed precision $s$-step Lanczos case, $\eps_2$ will contain a term $\bar{\Gamma}_k$.

Proceeding verbatim from Paige's work \cite{paige1980accuracy}, we write the eigendecomposition of $\hat{T}_m$ as
\begin{equation}
\hat{T}_m Q^{(m)} = Q^{(m)}\hspace{2pt} \text{diag}\big(\mu_i^{(m)} \big),
\label{eq:Tedec}
\end{equation}
for $i\in\{1,\ldots,m\}$, where the orthonormal matrix $Q^{(m)}$ has $i^{\text{th}}$ column $q_i^{(m)}$ and $(\ell,i)$ element $\eta_{\ell,i}^{(m)}$, and the eigenvalues are ordered such that 
\begin{equation*}
\mu_1^{(m)} > \mu_2^{(m)} > \cdots > \mu_{m}^{(m)}. 
\end{equation*}
If $\mu_i^{(m)}$ is an approximation to an eigenvalue $\lambda_i$ of $A$, then the corresponding approximate eigenvector is $z_i^{(m)}$, the $i$th column of
\begin{equation}
Z^{(m)} \equiv \hat{V}_m Q^{(m)}.
\label{eq:Zdef}
\end{equation}

A crucial part of Paige's analysis is the definition of stabilization of the $r^{th}$ Ritz value at iteration, $t$, denoted $\mu_r^{(t)}$:
\begin{definition} \cite[Definition 1]{paige1980accuracy}
We say that an eigenvalue $\mu_r^{(t)}$ of $\hat{T}_t$ \emph{has stabilized to within $\delta_{t,r}$} if, for every $m>t$, we know there is an eigenvalue of $\hat{T}_m$ within
\[
\delta_{t,r} \equiv \hat{\beta}_{t+1} |\eta_{t,r}^{(t)} |  \geq \min_{i} |\mu_i^{(m)} - \mu_r^{(t)}|
\]
of $\mu_{r}^{(t)}$. We will say $\mu_r^{(t)}$ has \emph{stabilized} when we know it has stabilized to within $\gamma(m+1)^{\omega} \sigma \eps_2$ where $\gamma$ and $\omega$ are small positive constants. 
\end{definition}
Thus, after $t$ steps, $\mu_r^{(t)}$ has necessarily stabilized to within $\delta_{t,r}$. 

Here we point out a crucial aspect of the difference between classical, uniform precision $s$-step, and mixed precision $s$-step Lanczos. The primary difference among these three variants is how tight the constraints are by which we consider an eigenvalue to be stabilized. The larger the value of $\eps_2$, the looser the constraint on stabilization becomes, and thus the sooner an eigenvalue is considered to be stabilized. Thus, somewhat counterintuitively, for the uniform precision $s$-step Lanczos process where $\eps_2$ is expected to be largest, we expect ``stabilization'' to happer \textit{sooner} than in the other methods (but again, to within a larger interval around the true eigenvalues of $A$). In the classical Lanczos method, the smaller value of $\eps_2$ means that we are more discriminating in what we consider to be a stabilized eigenvalue, and thus stabilization will take longer. For the mixed precision $s$-step Lanczos case, we expect the value of $\eps_2$ to fall somewhere in the middle of the other two variants. The significance of this will become clear after the remaining results are stated.

We also use the follow definition of Paige:
\begin{definition}\cite[Definition 2]{paige1980accuracy}
We will say that an eigenpair $(\mu, z)$ \emph{represents an eigenpair of $A$ to within} $\delta$ if we know that $\| Az - \mu z \| / \| z \| \leq \delta$.
\end{definition}

We can now state a number of results of Paige, rewritten for the mixed precision $s$-step Lanczos case:
\medskip
\begin{lem} 
\label{lem:lem1}
Let $\hat{T}_m$ and $\hat{V}_m$ be the result of $m$ steps of the mixed precision $s$-step Lanczos algorithm with~\eqref{eq:e0e1} and~\eqref{eq:eps2}, and let $R_m$ be the strictly upper triangular matrix defined in~\eqref{eq:Rm}. Then for each eigenpair $(\mu_i^{(m)}, q_i^{(m)})$ of $\hat{T}_m$, there exists a pair of integers $(r,t)$ with $0 \leq r \leq t < m$ such that
\begin{equation*}
\delta_{t,r} \equiv \hat{\beta}_{t+1} |\eta_{t,r}^{(t)} | \leq \psi_{i,m} \quad \text{and} \quad |\mu_{i}^{(m)} - \mu_r^{(t)}|\leq \psi_{i,m},  
\end{equation*}
where
\begin{equation*}
\psi_{i,m} \equiv \frac{ m^2 \sigma \eps_2}{\big|\sqrt{3} \hspace{2pt} q_i^{(m)T} R_m q_i^{(m)}\big| }.
\end{equation*}
\end{lem}
\medskip
\begin{proof} See \cite[Lemma 3.1]{paige1980accuracy}.
\end{proof}
\medskip
\begin{thm} 
If, with the con\-di\-tions of Lemma~\ref{lem:lem1}, an eigenvalue $\mu_i^{(m)}$ of $\hat{T}_m$ produced by mixed precision $s$-step Lanczos is stabilized so that 
\begin{equation}
\delta_{m,i} \equiv \hat{\beta}_{m+1} |\eta_{m,i}^{(m)}| \leq \sqrt{3} m^2 \sigma \eps_2,
\label{eq:deltael}
\end{equation}
and $\eps_0 < 1/12$, then for some eigenvalue $\lambda_c$ of $A$, 
\begin{equation}
|\lambda_c - \mu_i^{(m)}| \leq (m+1)^3 \sigma \eps_2.
\label{eq:lambmu}
\end{equation}
\label{thm:thm2}
\end{thm}
\medskip
\begin{proof} See \cite[Theorem 3.1]{paige1980accuracy}.
\end{proof}
\medskip

In fact, if~\eqref{eq:deltael} holds then we have an eigenvalue with a superior error bound to~\eqref{eq:lambmu} and we also have a good eigenvector approximation, as shown in the following Corollary. 
\medskip
\begin{cor}
\label{cor:cor1}
If~\eqref{eq:deltael} holds, then for the final $(r,t)$ pair in Theorem~\ref{thm:thm2}, $(\mu_{r}^{(t)}, \hat{V}_t q_r^{(t)})$ is an exact eigenpair for a matrix within $6t^2\sigma \eps_2$ of $A$. 
\end{cor} 
\medskip
\begin{proof} See \cite[Corollary 3.1]{paige1980accuracy} and~\cite[Corollary 5.3]{carson2014accuracy}.
\end{proof}
\medskip

Theorem~\ref{thm:thm2} says that, assuming~\eqref{eq:assumptions} holds, if an eigenvalue of $\hat{T}_m$ has stabilized to within $\sqrt{3} m^2 \sigma \eps_2$, then it is within $(m+1)^3 \sigma \eps_2$ of an eigenvalue of $A$, regardless of how many other eigenvalues of $\hat{T}_m$ are close. Corollary~\ref{cor:cor1} says that in this case we have an eigenpair of a matrix within $6m^2\sigma\eps_2$ of $A$. We stress again that the tightness of these intervals depends heavily on the value of $\eps_2$, which, in the mixed precision $s$-step Lanczos case, contains a factor of $\bar{\Gamma}_k$. 

We can now state some results of Paige \cite{paige1980accuracy} regarding convergence extended to the mixed precision $s$-step Lanczos case.
\medskip
\begin{thm} 
For the mixed precision $s$-step Lanczos algorithm, if $n(3\eps_0 + \eps_1)\leq 1$ and $\eps_0<1/12$, then at least one eigenvalue of $\hat{T}_{n}$ must be within $(n+1)^3\sigma\eps_2$ of an eigenvalue of the $n\times n$ matrix $A$, and there exist $r\leq t \leq n$ such that $(\mu_r^{(t)}, z_r^{(t)})$ is an exact eigenpair of a matrix within $6t^2 \sigma \eps_2$ of $A$. 
\end{thm}
\medskip
\begin{proof} See \cite[Theorem 4.1]{paige1980accuracy} and~\cite[Theorem 6.1]{carson2014accuracy}.
\end{proof}
\medskip

In the classical Lanczos case, this says that we have at least one eigenvalue of $A$ with high accuracy by iteration $m=n$. In both uniform and mixed precision $s$-step Lanczos, it is still true that we will find at least one eigenvalue with some degree of accuracy by iteration $m=n$ as long as \eqref{eq:assumptions} holds, but here the limit on accuracy is determined by the size of $\bar{\Gamma}^2_{\lceil n/s\rceil}$ in the uniform precision case and  $\bar{\Gamma}_{\lceil n/s\rceil}$ in the mixed precision case, which are both contained in the $\eps_2$ term. \textit{Thus we can expect in general, eigenvalue estimates will be about a factor $\bar{\Gamma}_{\lceil n/s\rceil}$ more accurate in the mixed precision case versus the uniform precision case. }

\medskip
\begin{thm} 
For $m$ iterations of the mixed precision $s$-step Lanczos algorithm, with~\eqref{eq:e0e1},~\eqref{eq:eps2}, and $m$ such that
\begin{equation}
\delta_{\ell,i} \equiv \hat{\beta}_{\ell+1} |\eta_{\ell,i}^{(\ell)}| \geq \sqrt{3} m^2 \sigma \eps_2, \quad  1\leq i \leq \ell < m,
\label{eq:delt2}
\end{equation}
the $m$ Lanczos vectors (columns of $\hat{V}_m$) span a  Krylov subspace of a matrix within $(3m)^{1/2} \sigma \eps_2$ of $A$. 
\end{thm}
\medskip
\begin{proof} See \cite[Theorem 4.2]{paige1980accuracy}.
\end{proof}
\medskip

This is analogous to the result of Paige for classical Lanczos: until an eigenvalue of $\hat{T}_{m-1}$ has stabilized, i.e., while~\eqref{eq:delt2} holds, the computed vectors $\hat{v}_1,\ldots,\hat{v}_{m+1}$ correspond to an exact Krylov sequence for the matrix $A+\delta A_m$. As a result of this, if we assume that the $s$-step bases generated in each outer loop are conditioned such that~\eqref{eq:assumptions} holds, then the (uniform and mixed precision) $s$-step Lanczos algorithm can be thought of as a numerically stable way of computing a Krylov sequence, at least until the corresponding Krylov subspace contains an exact eigenvector of a matrix within $6m^2 \sigma \eps_2$ of $A$. The issue is that, since the $\eps_2$ term contains either a factor of $\bar{\Gamma}^2_m$ or $\bar{\Gamma}_m$ in the uniform and mixed precision algorithms, respectively, then we can find an exact eigenvector of a matrix within $6m^2 \sigma \eps_2$ of $A$ much sooner, since this interval can be much larger. Comparing the uniform precision and mixed precision $s$-step variants, this also says that \textit{we expect that the mixed precision variant will compute a longer numerically stable Krylov sequence than in the uniform precision case.} 

When $\hat{T}_m$ and $\hat{V}_m$ are used to solve the eigenproblem of $A$, we want the eigenvalues and eigenvectors of $\hat{T}_m$ to be close to those of 
\begin{equation}
\hat{V}_m^T A \hat{V}_m q = \mu \hat{V}_m^T \hat{V}_m q, \quad \text{where } q^Tq =1,\label{eq:ideal}
\end{equation}
as would be the case with classical Lanczos with full reorthogonalization. If~\eqref{eq:delt2} holds, then the range of $\hat{V}_m$ is close to what we expect if we used classical Lanczos with full reorthogonalization. Thus the eigenvalues of~\eqref{eq:ideal} would be close (how close of course depends on the value of $\eps_2$, which is a factor $\bar{\Gamma}_k$ smaller in the mixed precision case) to those that would have been obtained using full reorthogonalization. 

\medskip
\begin{thm} 
\hspace{-.6pt} If $\hat{V}_m$ comes from the mixed precision $s$-step Lanc\-zos algorithm with~\eqref{eq:e0e1} and~\eqref{eq:eps2}, and~\eqref{eq:delt2} holds, then for any $\mu$ and $q$ which satisfy~\eqref{eq:ideal}, 
$(\mu,\hspace{-.6pt}\hat{V}_m q)$ is an exact eigenpair for a matrix within $\big(2\delta + 2 m^{1/2} \sigma \eps_2 \big) $ of $A$, where
\begin{equation*}
\eta \equiv e_{m}^T q, \quad \delta \equiv \hat{\beta}_{m+1} |\eta|. 
\end{equation*}
\end{thm}
\medskip
\begin{proof} See \cite[Theorem 4.3]{paige1980accuracy} and~\cite[Theorem 6.3]{carson2014accuracy}.
\end{proof}
\medskip

Following Paige \cite{paige1980accuracy}, it can also be shown that for each $\mu_i^{(m)}$ of $\hat{T}_m$,
\begin{equation*}
\min_{\mu \hspace{1mm}\text{in \eqref{eq:ideal}}} |\mu - \mu_i^{(m)}| \leq 2m^{1/2} \sigma \eps_2 + \frac{\sqrt{3}\delta_{m,i}}{\sqrt{m}}.
\end{equation*}
This means that when $(\mu_i^{(m)}, \hat{V}_m q_i^{(m)})$ represents an eigenpair of $A$ to within about $\delta_{m,i}$, there is a $\mu$ of~\eqref{eq:ideal} within about $\delta_{m,i}$ of $\mu_i^{(m)}$, assuming $m \geq 3$. 

Thus, assuming no breakdown occurs and the size of $\bar\Gamma_k$ satisfies~\eqref{eq:assumptions}, these results say the same thing for the mixed precision $s$-step Lanczos case as in the uniform precision $s$-step Lanczos and  classical Lanczos cases: \emph{until an eigenvalue has stabilized, the mixed precision $s$-step Lanczos algorithm behaves very much like the error-free Lanczos process, or the Lanczos algorithm with reorthogonalization}. 

Again, we stress that ``stabilization'', in the terms of Paige, occurs sooner the larger the value of $\eps_2$. This phenomenon occurs due to the fact that our constraint on what is considered to be ``stabilized'' is looser the larger the value of $\eps_2$. Thus in the uniform precision $s$-step Lanczos case, in which $\eps_2$ depends on $\bar{\Gamma}_k^2$, we expect stabilization will happen soonest, which here we see means that we expect faster deviation from the exact Lanczos process. In the mixed precision $s$-step Lanczos case, $\eps_2$ depends only on  $\bar{\Gamma}_k$, and thus we expect that the mixed precision algorithm will follow the exact Lanczos process for a greater number of iterations (but likely not as long as the classical Lanczos algorithm).

\section{A mixed precision $s$-step conjugate gradient algorithm}
\label{sec:cg}

The conjugate gradient (CG) method for solving linear systems is based on an underlying Lanczos process. Similarly, the $s$-step CG algorithm is based on an underlying $s$-step Lanczos algorithm. We therefore expect that the improved eigenvalue accuracy and orthogonality obtained by the use of the mixed precision approach in the $s$-step Lanczos algorithm will lead to improvements in convergence behavior in a corresponding mixed precision $s$-step CG algorithm. In this section, we present numerical experiments that support this conjecture. We note that we do not expect that the proposed use of extended precision in the Gram matrix computation will improve the maximum attainable accuracy in $s$-step CG, as this depends largely on the accuracy with which SpMVs are computed; see, e.g.,~\cite{carson2014residual} for bounds on the maximum attainable accuracy for $s$-step CG.

Details on the $s$-step CG algorithm can be found in, e.g.,~\cite{carson2015communication}, which also contains historical references. To keep the text simple, we do not include the algorithm here, but simply note that as in the mixed precision $s$-step Lanczos case, we use double the working precision in computing and applying the Gram matrix $G_k$ and working precision elsewhere. 

 We note that it is of future interest to extend the results of Greenbaum \cite{greenbaum1989behavior} for the classical CG algorithm to both the uniform and mixed precision $s$-step cases. Greenbaum's results say that finite precision classical CG behaves like exact CG applied to a larger matrix whose eigenvalues are in tight clusters around the eigenvalues of $A$. We expect that extending Greenbaum's analysis to the $s$-step variants will give a similar result, except the clusters may not be as tight; we anticipate that the cluster radius will contain a factor of $\bar{\Gamma}^2_k$ in the uniform precision case and a factor of $\bar{\Gamma}_k$ in the mixed precision case. 

We test the same two problems as in Section~\ref{sec:num}, the diagonal problem defined in \eqref{eq:strakos} (with the same parameters) and the matrix \texttt{nos4} from SuiteSparse \cite{davis2011university}, and use the same experimental setup. We additionally show results for the matrix \texttt{lund\_b}, also from SuiteSparse. For all matrices, we construct a right-hand side that has equal components in the eigenbasis of $A$ and unit $2$-norm (which presents a difficult case for CG) and take the zero vector as the initial approximate solution. 

For the both problems, we compare classical CG in double precision, uniform precision $s$-step CG in double precision, and mixed precision $s$-step CG in double/quad precision. In all plots, the y-axis shows the relative error in the $A$-norm, where the ``true solution'' was computed using MATLAB's backslash in quadruple precision via Advanpix \cite{ad06}. 

\begin{figure}
    \centering
    \includegraphics[trim={4cm 8cm 4cm 8cm},clip,width=6.0cm]{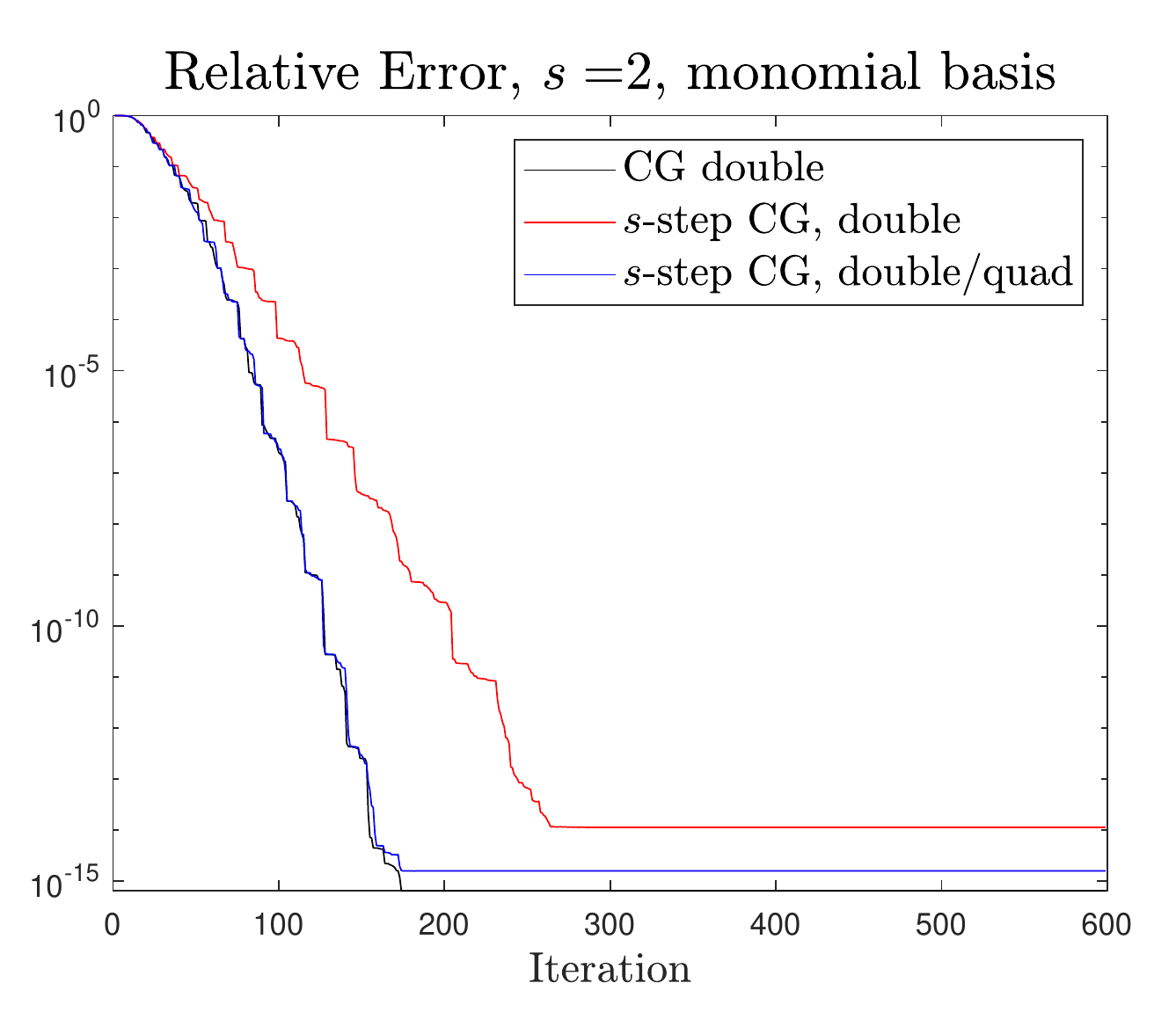}
    \includegraphics[trim={4cm 8cm 4cm 8cm},clip,width=6.0cm]{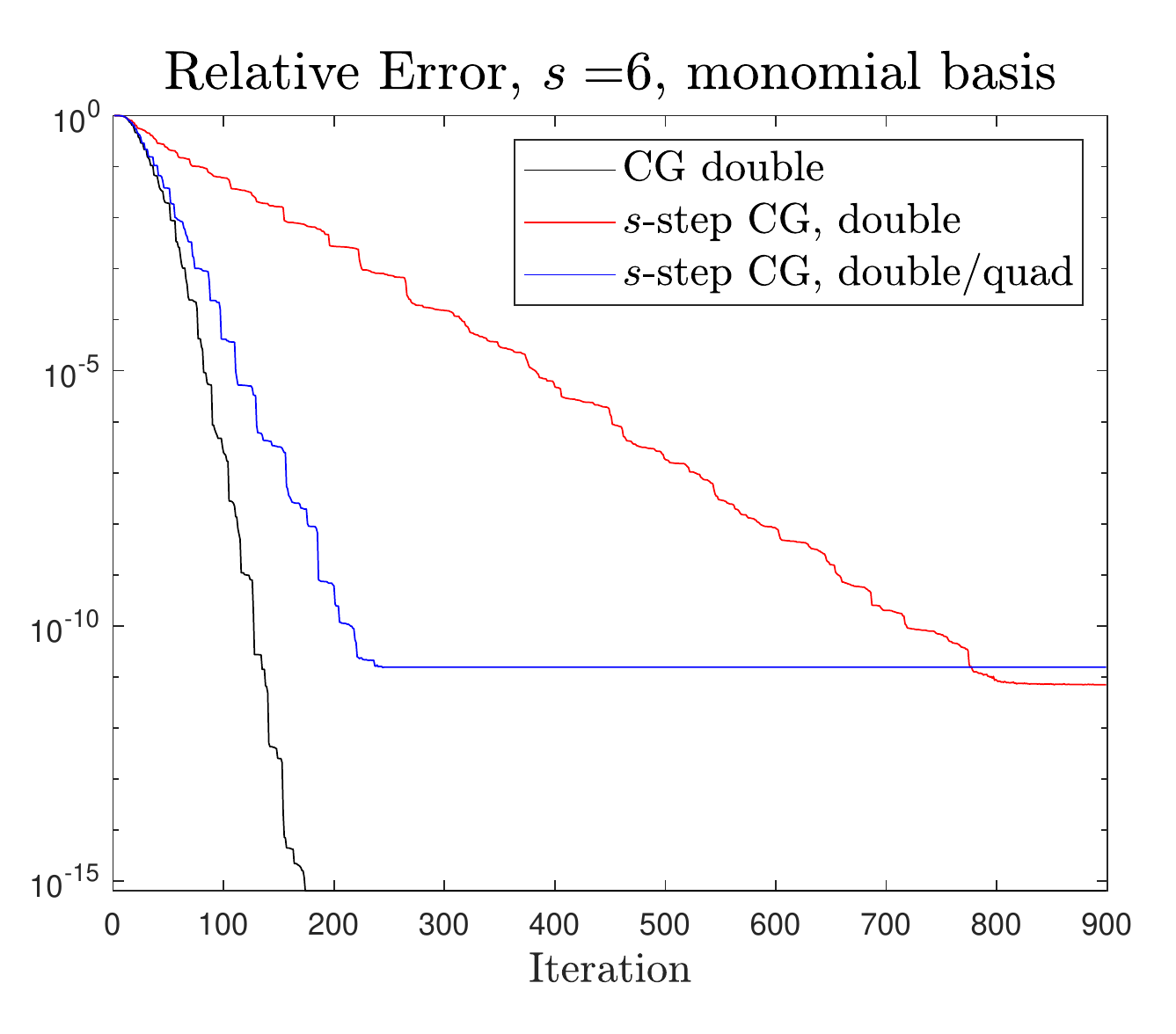}\\
    \includegraphics[trim={4cm 8cm 4cm 8cm},clip,width=6cm]{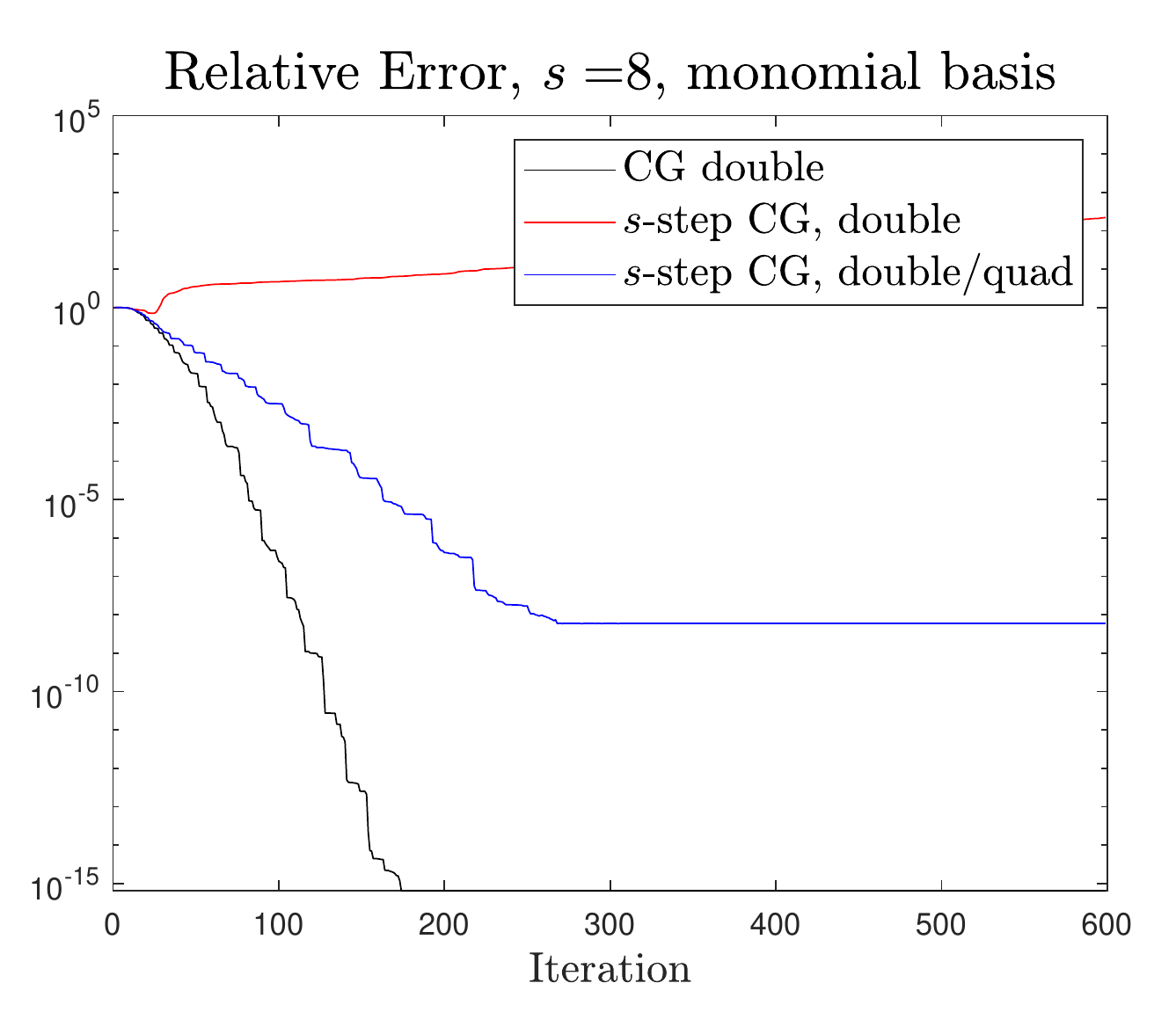}
    \includegraphics[trim={4cm 8cm 4cm 8cm},clip,width=6cm]{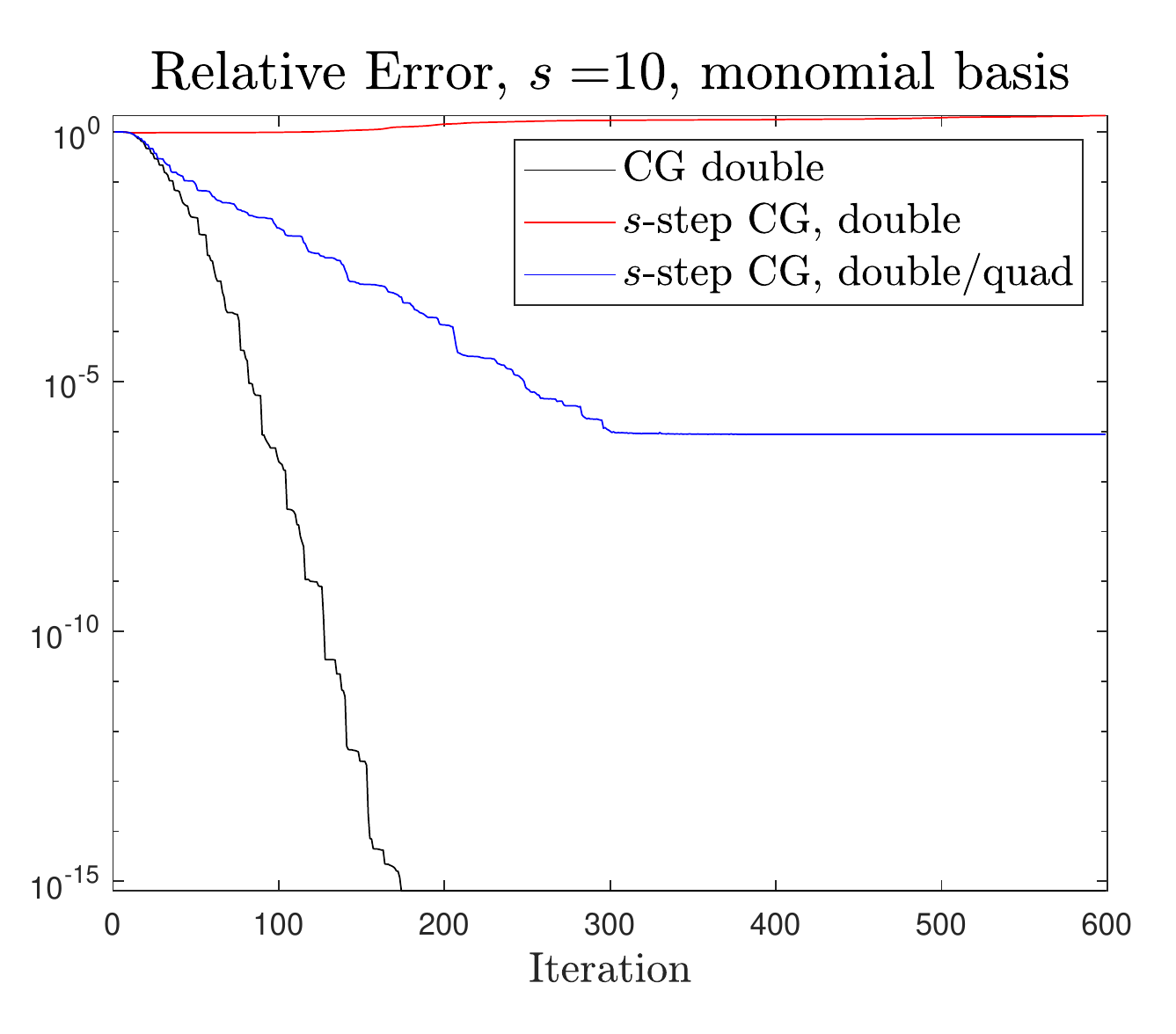}
    \caption{Comparison of uniform and mixed precision $s$-step CG for the diagonal matrix problem \eqref{eq:strakos} with $n=100$, $\lambda_1=10^{2}$, $\lambda_n=10^{-3}$, $\rho=0.65$, using $s\in\{2,6,8,10\}$ and a monomial basis, in terms of the relative error in the $A$-norm.}
    \label{fig:strakoscg}
\end{figure}

The results for the diagonal problem are shown in Figure \ref{fig:strakoscg}. Here we test $s$ values 2, 6, 8, and 10 and use the monomial basis. The advantage of the mixed precision approach from a numerical behavior standpoint is clear. For $s=2,6$ the convergence behavior of the mixed precision $s$-step algorithm is much closer to that of the classical algorithm. Note the extended x-axis in the upper right plot. For $s>6$, the uniform precision algorithm no longer converges. The mixed precision algorithm does still converge in this case, but the number of iterations until convergence increasingly deviates from the classical algorithm as $s$ increases. We also note here that, as the theory suggests, the mixed precision approach does not help improve the maximum attainable accuracy, which decreases with increasing $s$; see the bounds in \cite{carson2014residual}. 

\begin{figure}
    \centering
    \includegraphics[trim={4cm 8cm 4cm 8cm},clip,width=6.0cm]{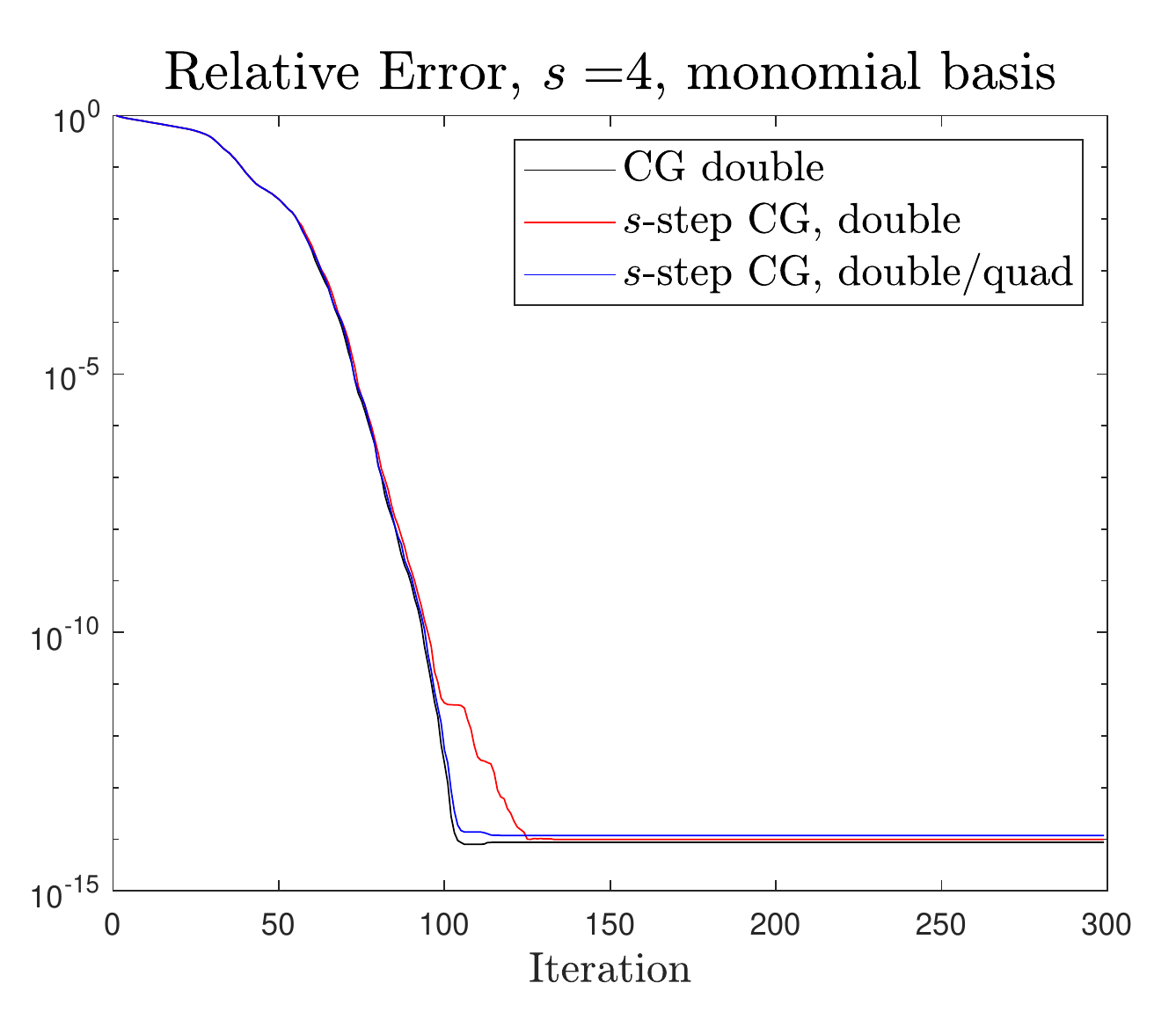}
    \includegraphics[trim={4cm 8cm 4cm 8cm},clip,width=6.0cm]{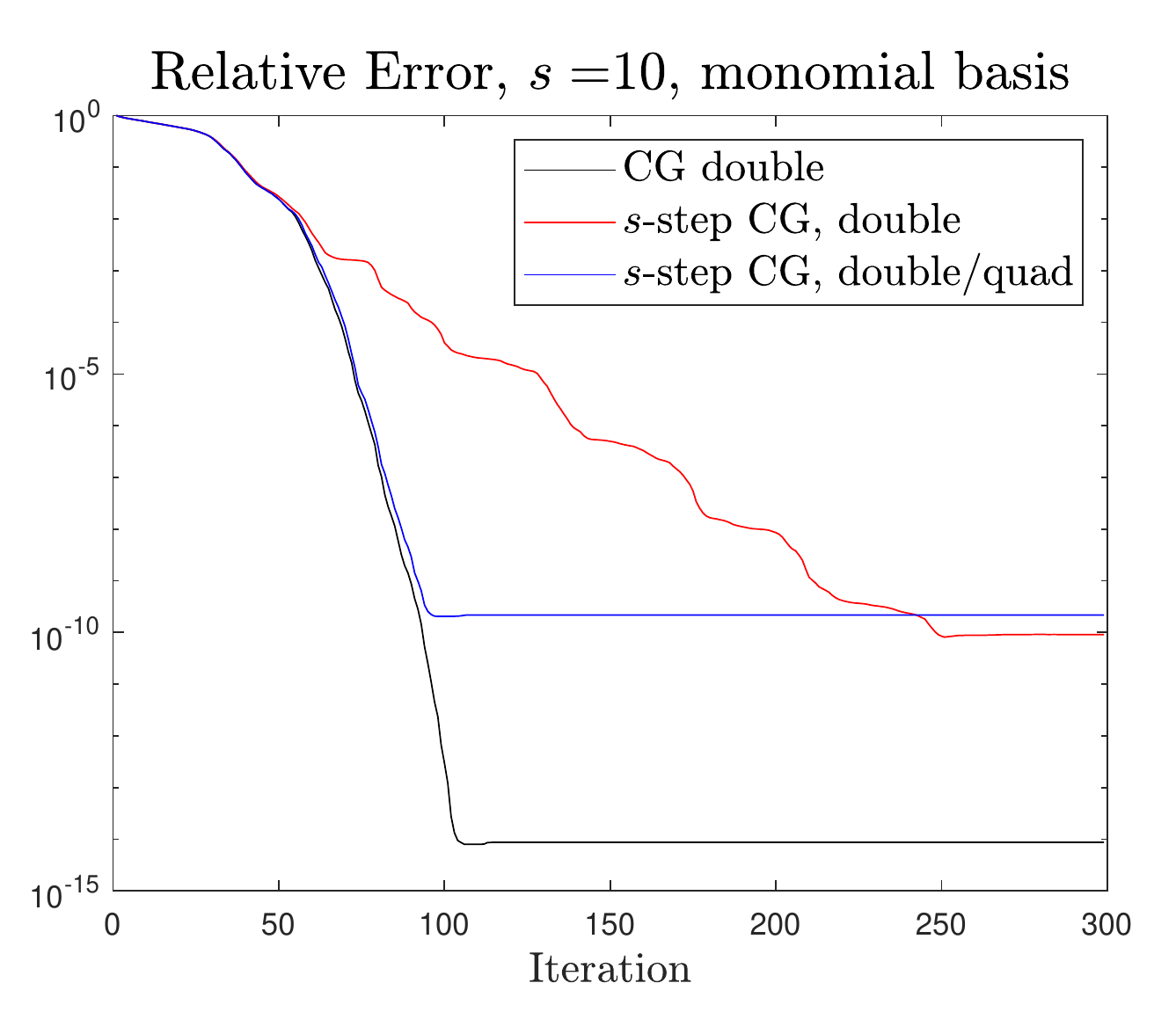}\\
    \includegraphics[trim={4cm 8cm 4cm 8cm},clip,width=6cm]{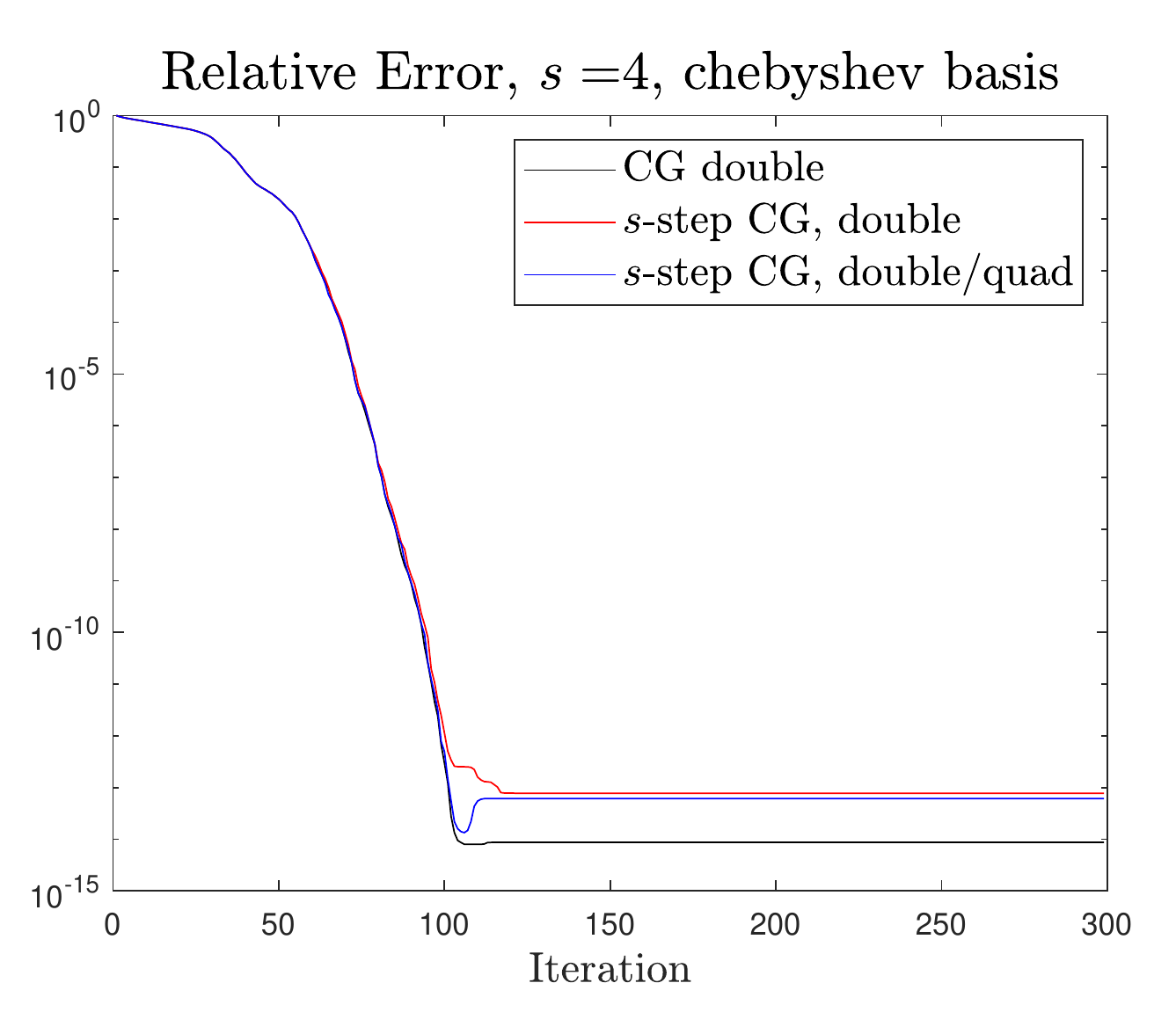}
    \includegraphics[trim={4cm 8cm 4cm 8cm},clip,width=6cm]{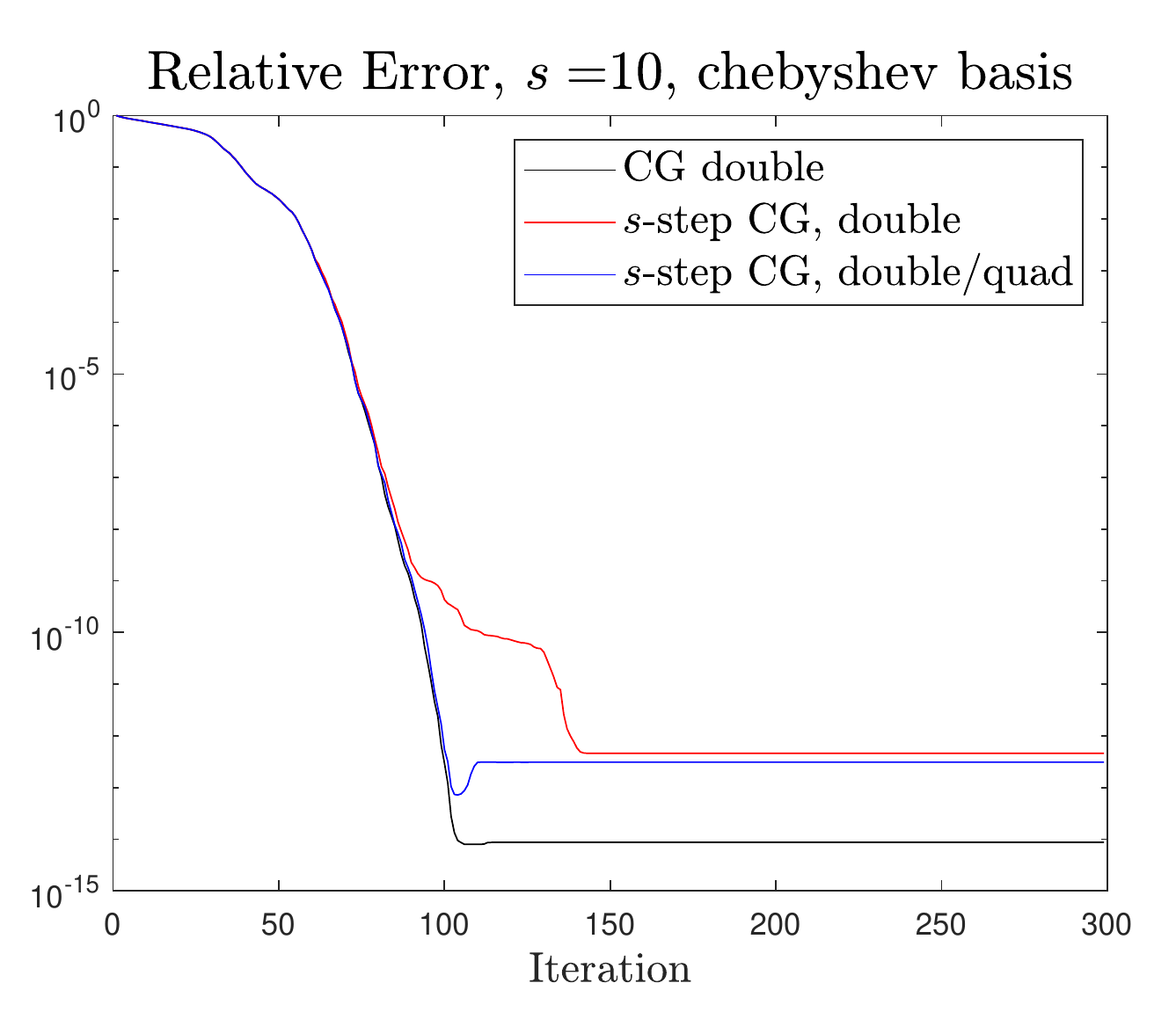}
    \caption{Comparison of uniform and mixed precision $s$-step CG for the \texttt{nos4} problem from SuiteSparse, with $s=4$ (left) and $s=10$ (right) with monomial (top) and Chebyshev (bottom) bases, in terms of the relative error in the $A$-norm.}
    \label{fig:nos4cg}
\end{figure}

The results for the \texttt{nos4} problem are shown in Figure \ref{fig:nos4cg}. Here we try 2 different $s$ values, 4 (left) and 10 (right) and two different polynomial bases, monomial (top) and Chebyshev (bottom). In both cases, for small $s$ values, the behavior of the uniform precision $s$-step CG algorithm is already quite close to that of the classical algorithm, so the extended precision does not significantly change things. For larger $s$ values, however, one can see the benefit of the mixed precision approach, in particular in the case of the monomial basis, where we expect $\bar\Gamma_k$ to be larger. Even for larger $s$ values, the mixed precision $s$-step CG algorithm behaves very similarly to the classical algorithm in terms of convergence. Again, we note the limitations on maximum attainable accuracy, which is linked to the value of $\bar\Gamma_k$ in both uniform and mixed precision settings. 

\begin{figure}
    \centering
    \includegraphics[trim={4cm 8cm 4cm 8cm},clip,width=6.0cm]{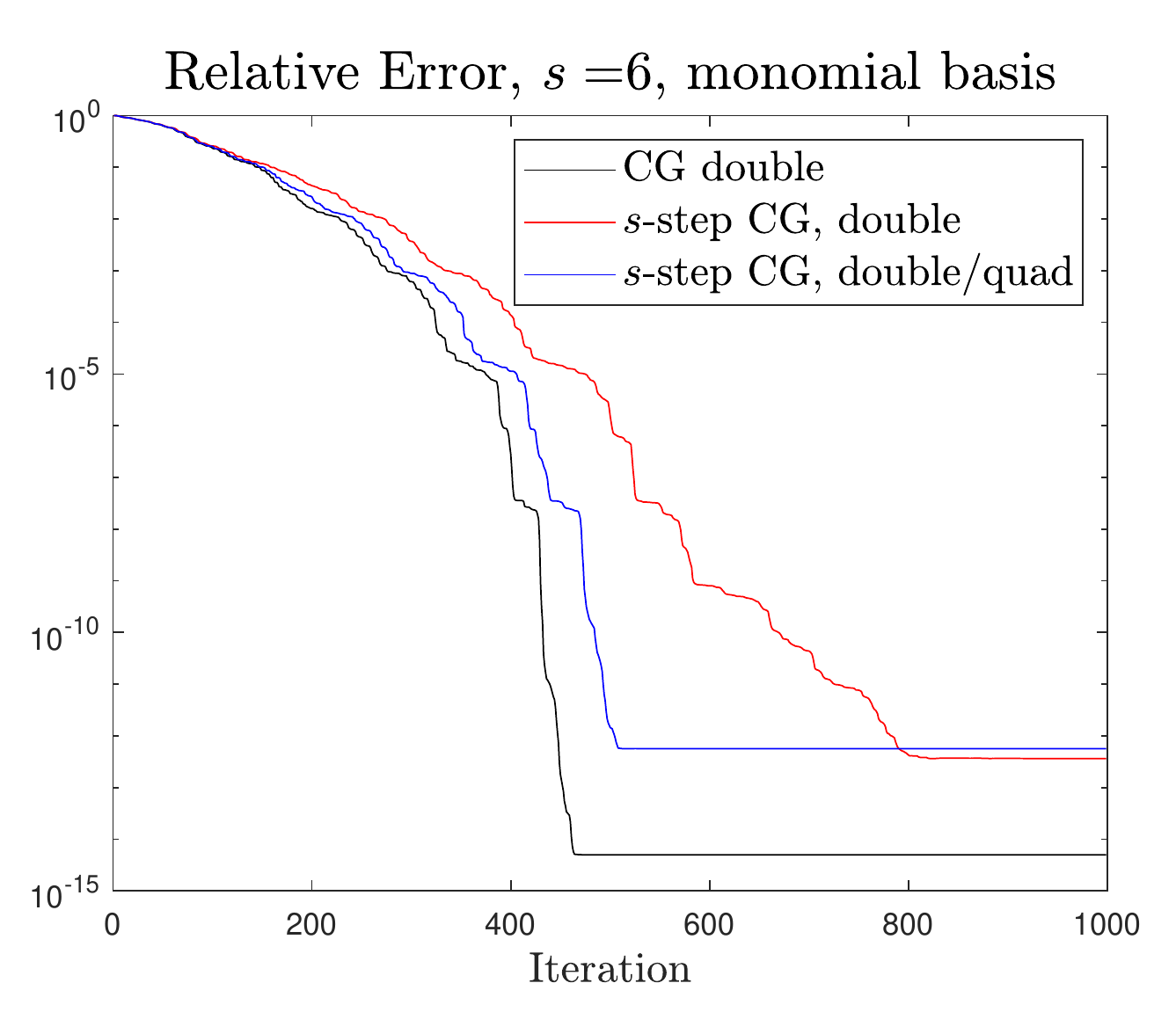}
    \includegraphics[trim={4cm 8cm 4cm 8cm},clip,width=6.0cm]{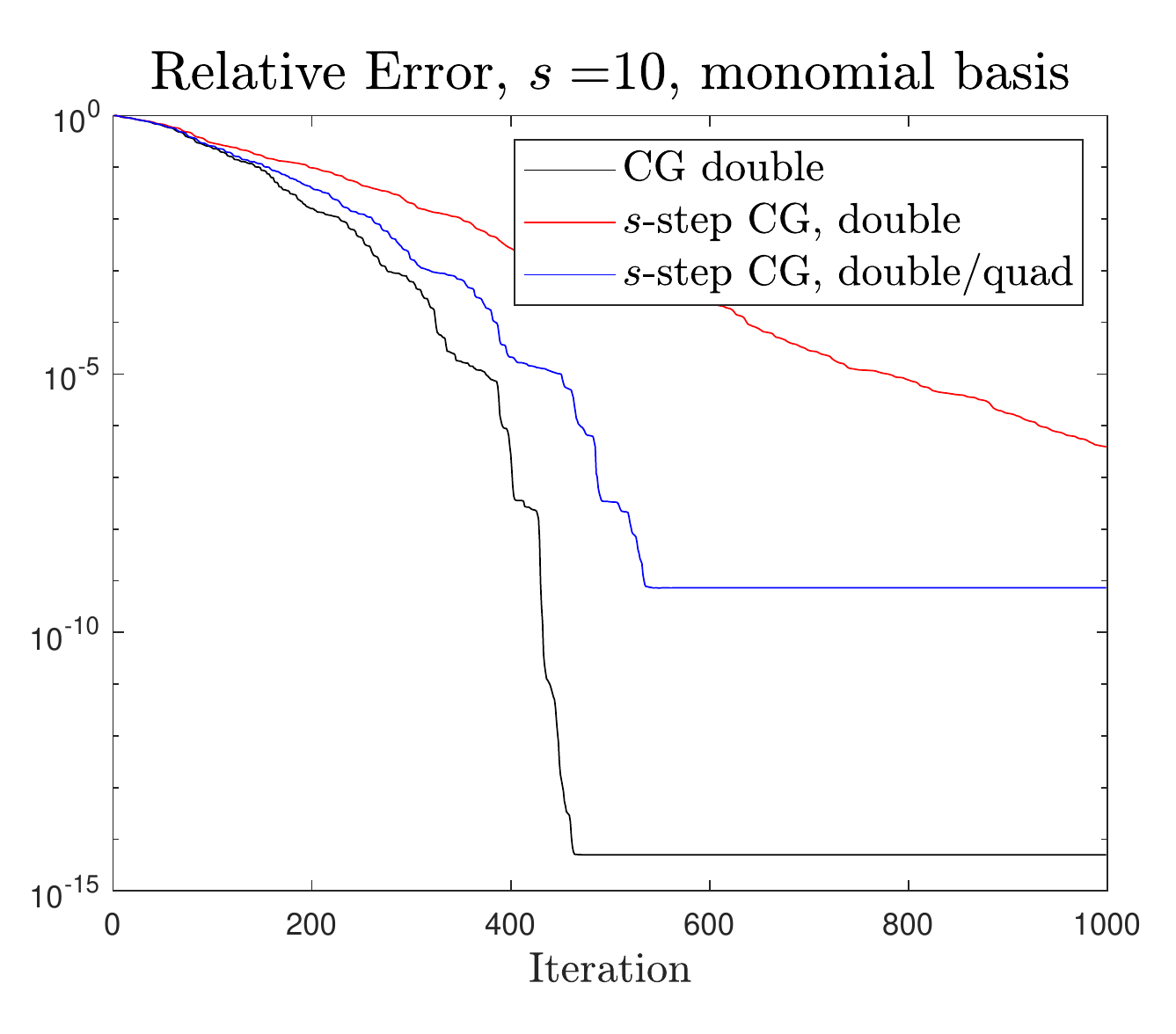}\\
    \includegraphics[trim={4cm 8cm 4cm 8cm},clip,width=6cm]{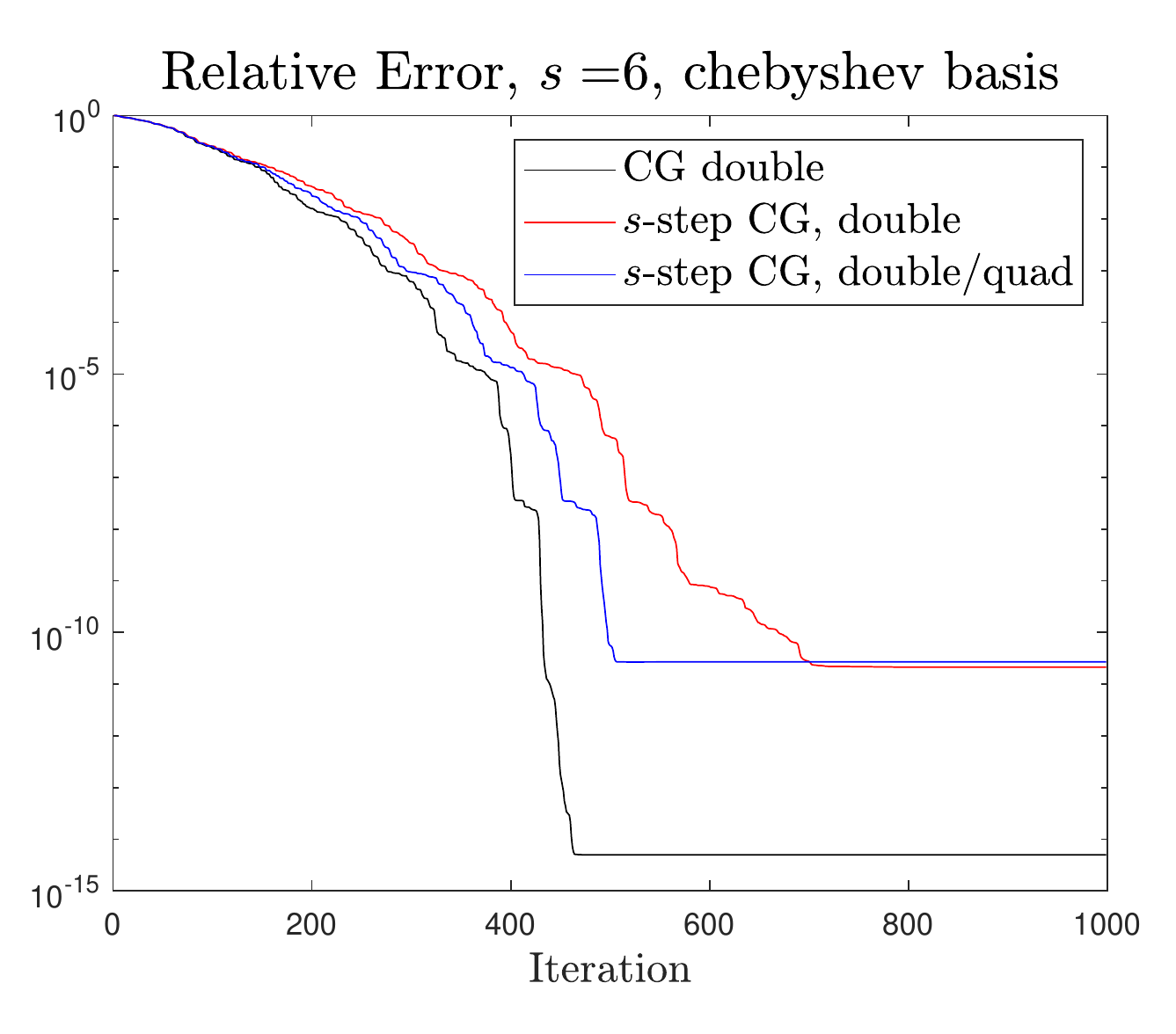}
    \includegraphics[trim={4cm 8cm 4cm 8cm},clip,width=6cm]{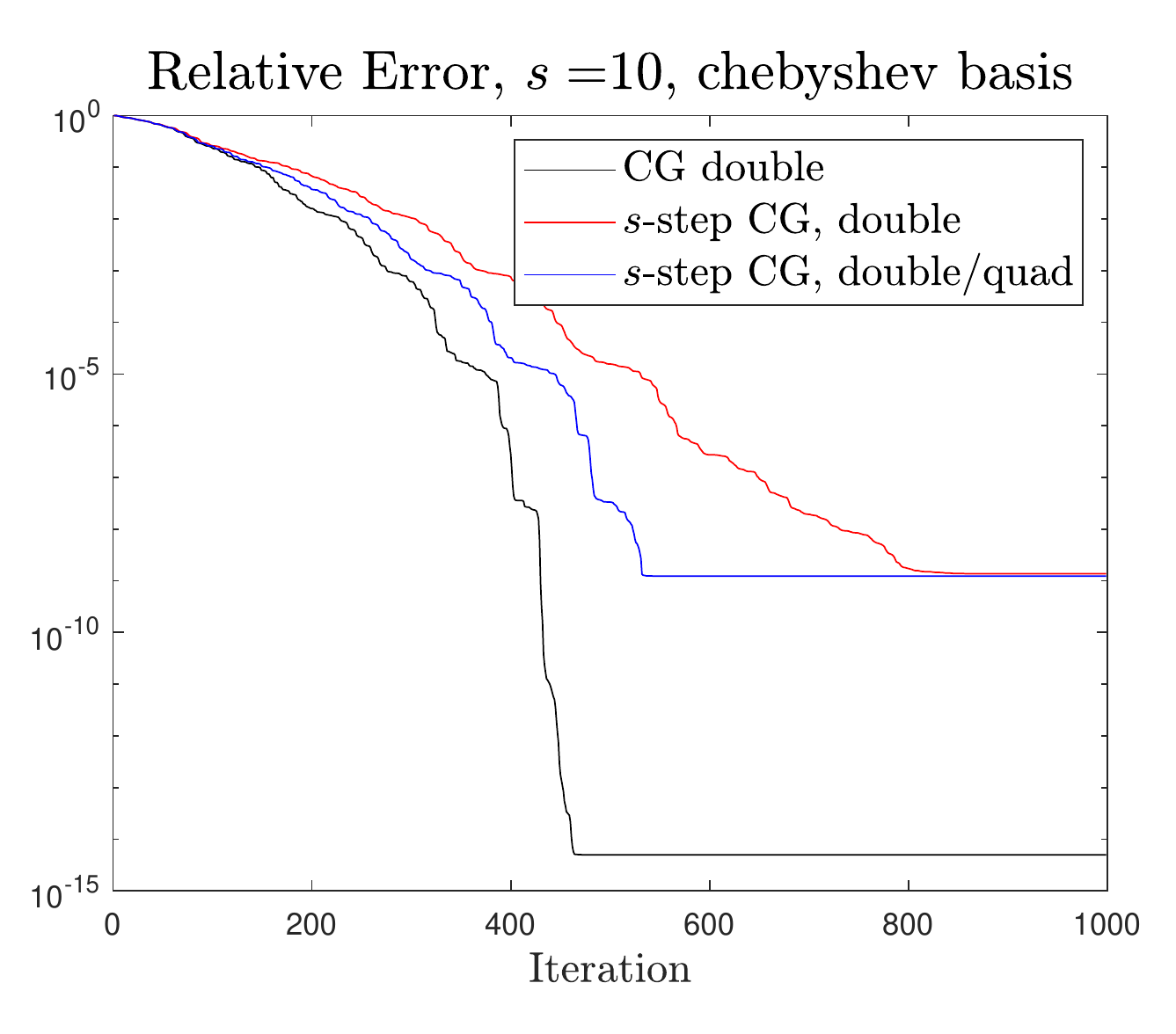}
    \caption{Comparison of uniform and mixed precision $s$-step CG for the \texttt{lund\_b} problem from SuiteSparse, with $s=6$ (left) and $s=10$ (right) with monomial (top) and Chebyshev (bottom) bases, in terms of the relative error in the $A$-norm.}
    \label{fig:lundbcg}
\end{figure}

The results for the \texttt{lund\_b} problem are shown in Figure \ref{fig:lundbcg}. For this problem $n=147$ and the condition number is about $3\cdot 10^4$. We test $s=6$ (left) and $s=10$ (right) and monomial (top) and Chebyshev (bottom) bases. Here we see that indeed, convergence of the mixed precision approach is much improved versus the uniform precision case. In this particular case, in contrast to \texttt{nos4}, it is clear that using an improved basis is not enough to correct the instabilities of the $s$-step approach. 

We argue that because the very selective use of extra precision is unlikely to cause a significant performance overhead, the techniques of using more well-conditioned polynomial bases and the use of mixed precision should be considered orthogonal. Both techniques should be used
in combination for the best behavior in practice as they improve numerical behavior through different means; using a more well-conditioned polynomial basis constructed using Chebyshev or Newton polynomials will reduce the value of $\bar{\Gamma}_k$ itself, and the use of mixed precision will reduce the dependence of the loss of orthogonality (and other quantities) from $\bar{\Gamma}_k^2$ to $\bar{\Gamma}_k$.

\section{Conclusions}
\label{sec:conclusion}

In this work, we have developed mixed precision variants of the $s$-step Lanczos and $s$-step CG algorithms. The insight is that by using double the working precision in a few select computations, namely, the computation and application of the Gram matrix, the bounds on the loss of orthogonality and other important measures are reduced by a factor related to the conditioning of the $s$-step bases. 

We extend the results of Paige in \cite{paige1976error} and \cite{paige1980accuracy} to the mixed precision $s$-step Lanczos algorithm and discuss the interpretation of the theorems. Key insights into the notion of ``stabilization'' of approximate eigenvalues are discussed, which provides theoretical explanation of why deviation from the exact Lanczos procedure happens earlier in $s$-step variants. We provide a few small numerical experiments in MATLAB that confirm the behavior expected by the theory. 

While we have also developed a mixed precision variant of $s$-step CG, there is still remaining work to do extended the classical ``backward-like'' stability results of Greenbaum \cite{greenbaum1989behavior} to both the uniform precision and mixed precision $s$-step CG algorithms. Based on existing insights, we expect that the theory will say that, under some constraints, the $s$-step variants behaves like exact CG applied to a larger matrix whose eigenvalues are in tight clusters around the eigenvalues of $A$, where the cluster radius will contain a factor of $\bar{\Gamma}^2_k$ in the uniform precision case and a factor of $\bar{\Gamma}_k$ in the mixed precision case. Extending Greenbaum's analysis is largely technical in nature and remains future work. 

Another missing piece is experimental confirmation that the overhead of using double the working precision as described is not too high in large-scale parallel settings. We predict that extra cost per iteration will not be significant, especially in latency-bound settings, and that in any case, it is likely that the improved convergence behavior will lead to faster time-to-solution. 

As noted in the theory and confirmed in the experiments, the mixed precision approach can improve convergence behavior, but does not improve the maximum attainable accuracy, which is largely limited by the precision of the SpMVs, and in both uniform and mixed precision $s$-step algorithms, depends on $\bar{\Gamma}_k$. One possibility to also improve the maximum attainable accuracy is to combine the mixed precision approach with a residual replacement strategy, which has been developed for $s$-step CG in \cite{carson2014residual}.

\bibliographystyle{siam}
\bibliography{main}

\end{document}